\documentclass[sn-mathphys-num]{sn-jnl}

\usepackage{graphicx}%
\usepackage{multirow}%
\usepackage{amsmath,amssymb,amsfonts}%
\usepackage{amsthm}%
\usepackage{mathrsfs}%
\usepackage{appendix}%
\usepackage{xcolor}%
\usepackage{textcomp}%
\usepackage{manyfoot}%
\usepackage{booktabs}%
\usepackage[noend]{algorithmic}%
\usepackage{algorithm}%
\usepackage{listings}%
\usepackage{afterpage}

\theoremstyle{thmstyleone}%
\newtheorem{theorem}{Theorem}

\theoremstyle{thmstyletwo}%

\theoremstyle{thmstylethree}%

\usepackage{orcidlink}
\usepackage{enumerate}
\usepackage{accents}
\newcommand\munderbar[1]{%
\underaccent{\bar}{#1}}
\usepackage{mathdots}
\usepackage{booktabs}
\usepackage{mathtools}

\usepackage{pgfplots}
\pgfplotsset{compat=newest} 

\usepackage{booktabs}\let\cline\cmidrule
\usepackage{lscape}  

\usepackage{tikz}
\usepackage{tikz-3dplot}

\usepackage{subcaption}

\newtheorem{proposition}{Proposition}
\newtheorem{lemma}{Lemma}
\newtheorem{remark}{Remark}

\usepackage{amssymb}

\usepackage{geometry} 
\usepackage{cleveref}
\crefname{enumi}{}{\unskip}

\usepackage{verbatim}

\usepackage{appendix}
\usepackage{titlesec}

\begin{document}

\title[Gas Transportation and Storage]{Global Optimization of Gas Transportation and Storage: Convex Hull Characterizations and Relaxations}

\author[1]{\fnm{Bahar Cennet} \sur{Okumuşoğlu} \orcidlink{0000-0002-5036-5687}
}\email{bahar.okumusoglu@tum.de}

\author[2]{\fnm{Burak} \sur{Kocuk} \orcidlink{0000-0002-4218-1116
}}\email{burak.kocuk@sabanciuniv.edu}


\affil[1]{\orgdiv{TUM School of Management}, \orgname{Technical University of Munich}, \orgaddress{\postcode{80333}, \city{Munich}, \country{Germany}}}

\affil[2]{\orgdiv{Industrial Engineering Program}, \orgname{Sabancı University}, \orgaddress{ \postcode{34956}, \city{Istanbul}, \country{Turkey}}}


\abstract{Gas transportation and storage has become one of the most relevant and important optimization problems in energy systems. This problem inherently includes highly nonlinear and nonconvex aspects due to gas physics, and discrete aspects due to the control decisions of active network elements. Obtaining even locally optimal solutions for this problem presents significant mathematical and computational challenges for system operators. In this paper, we formulate the gas transportation and storage problem as a nonconvex mixed-integer nonlinear program (MINLP) through disjunctions on the flow directions. Moreover, we study the nonconvex sets induced by gas physics and propose mixed-integer second-order cone programming relaxations for the nonconvex MINLP problem. The proposed relaxations are based on the convex hull representations of two nonconvex sets: Firstly, we give the convex hull representation of the nonconvex set for pipes and show that it is second-order cone representable. Secondly, we also give a complete characterization of the extreme points of the nonconvex set for compressors and show that the convex hull of the extreme points is power cone representable. Moreover, for practical applications, we propose a second-order cone outer-approximation for the nonconvex set for compressors. To obtain (near) globally optimal solutions, we develop an algorithmic framework based on our convex hull results. We evaluate our framework through extensive computational experiments on various GasLib networks in comparison with the convex relaxations from the literature and a state-of-the-art global solver. Our results highlight the computational efficiency and convergence performance of our convex relaxation method compared to other methods. Moreover, our method also consistently provides (near) global solutions as well as high-quality warm-starting points for local solvers.}




\maketitle

\section{Introduction}
Gas transportation is one of the most important optimization problems in gas network systems. It determines the optimal nodal pressures and gas flows through pipelines while balancing the gas supply and demand with minimum operational costs. 
In addition, the recent technological developments in the gas industry have highlighted the significant role of gas storage \cite{IEA2023}. In this respect,  considering storage sites in gas transportation optimization problem has become highly relevant and significant to the system operators. The gas stores provide considerable  security and flexibility for the network operations as they are the most responsive elements in case of supply disruptions and demand fluctuations. However, 
there remain certain mathematical and computational challenges in the gas transportation and storage problem summarized as follows: 
\begin{itemize}
    \item \textbf{Nonlinear and nonconvex relations}: The pressure loss in pipes and the fuel consumption of compressors are described by highly challenging nonconvex and nonlinear equations governed by the laws of physics. 
    \item \textbf{Discrete decisions}: The control decisions of active network elements and bidirectional gas flows require disjunctive formulations and binary decision variables. 
    \item \textbf{Multi-period structure}: Gas injection and withdrawal rates for stores and switching decisions of active network elements necessitate intertemporal constraints. 
\end{itemize}
The gas transportation and storage problem with the aforementioned nonconvex, nonlinear, and discrete aspects belongs to the class of nonconvex mixed-integer nonlinear programming (MINLP) problems. These problems are known to be very hard to solve to optimality as even their continuous relaxations involve nonconvexities.  


\subsection{Literature Review}

The gas networks are characterized by the existence of controllable network elements. These networks are called \textit{passive} if they only consist of controllable elements such as pipes and  \textit{active} if there also exist controllable elements such as compressors and (control) valves. 
There is a growing literature on the optimization problems considering different network characteristics as well as their computational complexities 
(for a recent review, see \cite{Rios2015} and references therein). 
The gas transportation problem in \cite{Gross2019} can be efficiently solved for single-source, single-sink passive networks, and gas nominations can be decided in polynomial time in single-cycle passive networks \cite{labbe2021polynomial}. 
In the existence of active network elements such as compressors and (control) valves, gas network optimization problems become NP-hard. For example, the work by \cite{Humpola2014} shows that the topology optimization problem in an active gas network is strongly NP-hard. The optimization problems in active gas networks are also known to be practically intractable within the scope of current state-of-the-art global solvers. Therefore, obtaining globally optimal or even local solutions for optimization problems in active networks is a very challenging task. 

Different solution methodologies are applied to solve the gas network optimization problems, which can be categorized into three groups: nonlinear optimization methods, approximation-based approaches, and conic programming relaxations. In nonlinear optimization methods, governing equations can be tackled in their original forms. These methods include sequential linear and quadratic programming \citep{Wolf2000, Ehrhardt2005}, interior point methods \citep{Steinbach2007, Schmidt2015} and primal-relaxed dual decomposition methods \citep{wu2007gas}. However, these studies do not include active network elements or assume that their control decisions are known a priori. While being able to find locally optimal solutions, these methods provide only a primal bound on the optimal objective value without any optimality guarantees for nonconvex nonlinear problems. Regularization schemes are also applied to nonlinear models with complementarity conditions for discrete aspects \citep{Schmidt2013, Pfetsch2015}. Although the solution approach in \cite{Pfetsch2015} is promising to obtain near-global solutions, the authors observe some infeasibility and numerical inaccuracy issues. Also, note that such issues may arise in local methods without acceptable warm-starting points.


Approximation-based approaches are widely applied to address nonlinear and nonconvex functions in the literature \citep{Zhang1996, Wu2000, Andre2009, Babonneau2012, Fugenschuh2013, Wu2017Mono, Ordoudis2019}. These approaches are mainly based on piecewise linearizations. Bilinear terms in the nonconvex equations are typically replaced by McCormick inequalities (see, e.g., \citep{Borraz2016SOCP, Wu2017Mono}). Although this approach is commonly used in global optimization, it might result in weak formulations when continuous variables have large bounds. In fact, this is the case for nodal pressure and gas flow bounds in gas network optimization problems. Piecewise linear functions are more successful in approximating inherent nonlinearities, which requires the additional of binary variables \citep{Martin2006, Zheng2010, Pfetsch2015, Wu2017Mono, Wang2018}. These piecewise linear approximations lead to mixed-integer linear programming (MILP) models, which are very convenient for the current state-of-the-art commercial solvers. High-resolution solutions produced by these approximations still require many binary variables and MILP formulations might become prohibitively expensive, especially in the presence of discrete control decisions.

Conic programming relaxations have recently attracted attention to deal with nonconvex and nonlinear aspects of gas physics. A semidefinite programming (SDP) relaxation is used to solve the nonconvex equations in a small tree network \cite{Ohja2017SDP}. Although SDP relaxations are usually tight, they become computationally more expensive as the network size increases. 
As an alternative, second-order cone programming (SOCP) relaxations have become popular among researchers due to their scalability and exactness under certain assumptions \citep{Borraz2016SOCP, BorrazSanchez2016-2, Wu2017Mono, Wen2018, He2018, Singh2019, Schwele2019SOCP}. In radial gas networks with known flow directions, pressure loss equations can be relaxed into second-order cone constraints \citep{He2018}. However, in the more realistic case of arbitrary flow directions, the situation becomes more challenging and necessitates disjunctive constraints. Bidirectional flows are commonly reformulated with flow direction (binary) variables and lead to mixed-integer second-order cone programming (MISOCP) models through big-M formulations \citep{Singh2019, Schwele2019SOCP} and McCormick inequalities \citep{Borraz2016SOCP, Schwele2019SOCP, li2024misocp}. Extended formulations built on perspective functions are also used (see, e.g., \citep{Borraz2016SOCP, BorrazSanchez2016-2, Wen2018, li2024misocp}). While the proposed relaxations in previous works are computationally tractable, they do not lead to convex hull formulations.

Next, we address the existing gaps in the literature. The optimization problems in gas networks include highly detailed nonlinear, nonconvex, and discrete aspects as mentioned above. Several studies use simplified models and focus on different aspects of these problems. For example, feasible sets for pipes induced by nonconvex pressure equations are commonly approximated by mixed-integer conic sets \citep{Singh2019, Schwele2019SOCP, Borraz2016SOCP, He2018,  BorrazSanchez2016-2, Wen2018}. The novel work in \cite{Borraz2016SOCP} proposes stronger formulations than previous works under arbitrary gas flows. The authors also emphasize the benefits of MISOCP models. However, the proposed mixed-integer conic sets do not give the convex hull representations of the nonconvex sets for pipes, and the proposed formulations produce weak SOCP relaxations. To approximate the operating domains for compressors simplified models (hereafter referred to as \textit{the feasible set for compressors}) are widely used in the literature \citep{Pfetsch2015, burlacu2019, BorrazSanchez2016-2}. However, the optimization models in active networks also assume idealized compressors, i.e., no fuel consumption \citep{Borraz2016SOCP, Schwele2019SOCP, Singh2019}, or use outer-approximations to preserve convexity \citep{BorrazSanchez2016-2, Wen2018}. 
To the best of our knowledge, convex hull descriptions and tight outer-approximations of the feasible sets for these network elements have not been studied in the literature. This paper fills these literature gaps as discussed in the next section.



\subsection{Our Contributions}

In this paper, we aim to find global solutions produced by mixed-integer conic relaxations for the general class of optimization problems in gas networks, including the aforementioned challenges. For this purpose, we study the gas transportation problem and extend it to a multi-period setting with gas stores. We also include the control decisions of active network elements such as compressors and (control) valves. We formulate the gas transportation and storage problem as a nonconvex MINLP. This problem includes 
highly nonlinear and nonconvex characteristics of the underlying gas physics, namely, two nonconvex feasible sets induced by the nonconvex pressure loss equations for pipes and the fuel consumption equations for compressors. We propose an MISOCP relaxation for this problem. In particular, our MISOCP relaxation is based on the exact convex hull representation of the feasible set for pipes and the outer-approximation of the convex hull of the feasible set for compressors.

Next, we elaborate on our main contributions: 
\begin{itemize}
    \item \textbf{Convexification of the feasible set for pipes}:  
    We focus on the nonconvex feasible set of a single pipe. 
    In Theorem~\ref{theorem:convexHull_pipes}, we show that the convex hull of this set is second-order cone representable (SOCr). We also give an extended formulation for the convex hull of this feasible set as Theorem~\ref{theorem:convexHull_pipes} is constructive. To the best of our knowledge, this is the first convex hull formulation of the feasible set for pipes. Moreover, our approach can be extended for those optimization problem that include nonconvex sets with similar characteristics (e.g., water supply networks).
    \item \textbf{Convexification of the feasible set for compressors}: 
    We later focus on the nonconvex feasible set of a single compressor. We give a complete characterization of the extreme points of this nonconvex set. Moreover, we show that the convex hull of the extreme points is power cone representable (POWr) in Theorem~\ref{theorem:convexHull_compressors_union}. To the best of our knowledge, this is the first convex hull characterization of the feasible set for compressors. Moreover, we also propose a mixed-integer second-order cone outer-approximation for this nonconvex set in Theorem~\ref{theorem:compressorOuterApprox} since it is not quite practical to solve mixed-integer power cone programs directly using the current state-of-the-art solvers. 
    \item \textbf{Algorithmic developments}: 
    To obtain global or near-global optimal solutions of the problem considered in this paper, we present an algorithmic framework based on our main convex hull results, in particular, Theorems~\ref{theorem:convexHull_pipes}, ~\ref{theorem:convexHull_compressors_union} and \ref{theorem:compressorOuterApprox}.
\end{itemize}

In order to test the empirical efficacy of our theoretical and algorithmic developments, we conduct an extensive computational study and compare the proposed framework with the state-of-the-art global solver. We also compare our relaxations with different conic relaxations from the literature under different objective functions. Our computational results show that the proposed framework significantly outperforms the global solver in terms of computational performance and convergence. Our results also illustrate the strength of our conic relaxations in comparison with the relaxations from the literature. Moreover, our approaches offer numerically consistent global and near-global solutions as well as high-quality warm-starting points for local solvers.

The rest of the paper is organized as follows: In Section \ref{section:problemFormulation}, we formally introduce the multi-period gas transportation and storage problem. 
We present our main results on the (nonconvex) feasible sets for pipes and compressors in Section~\ref{section:mainResults}. The proposed algorithmic framework is presented in Section~\ref{section:solutionAlgorithm}. The numerical results for our computational experiments on various GasLib instances follow in Section~\ref{section:computationalStudy}. Finally, we conclude our paper with further remarks and future research directions in Section~\ref{section:conclusion}.

\section{Problem Formulation} \label{section:problemFormulation}
We describe the problem setting in Section~\ref{ss:ProblemSetting}, and explain the gas network modeling along with decision variables and problem parameters in Section~\ref{ss:GasNetworkComponents}. Finally, we present the MINLP formulation of the problem in Section~\ref{ss:MINLPFormulation}.
\subsection{Problem Setting} \label{ss:ProblemSetting}
\begin{figure}[h]
\centering
 \includegraphics[scale = 0.38]{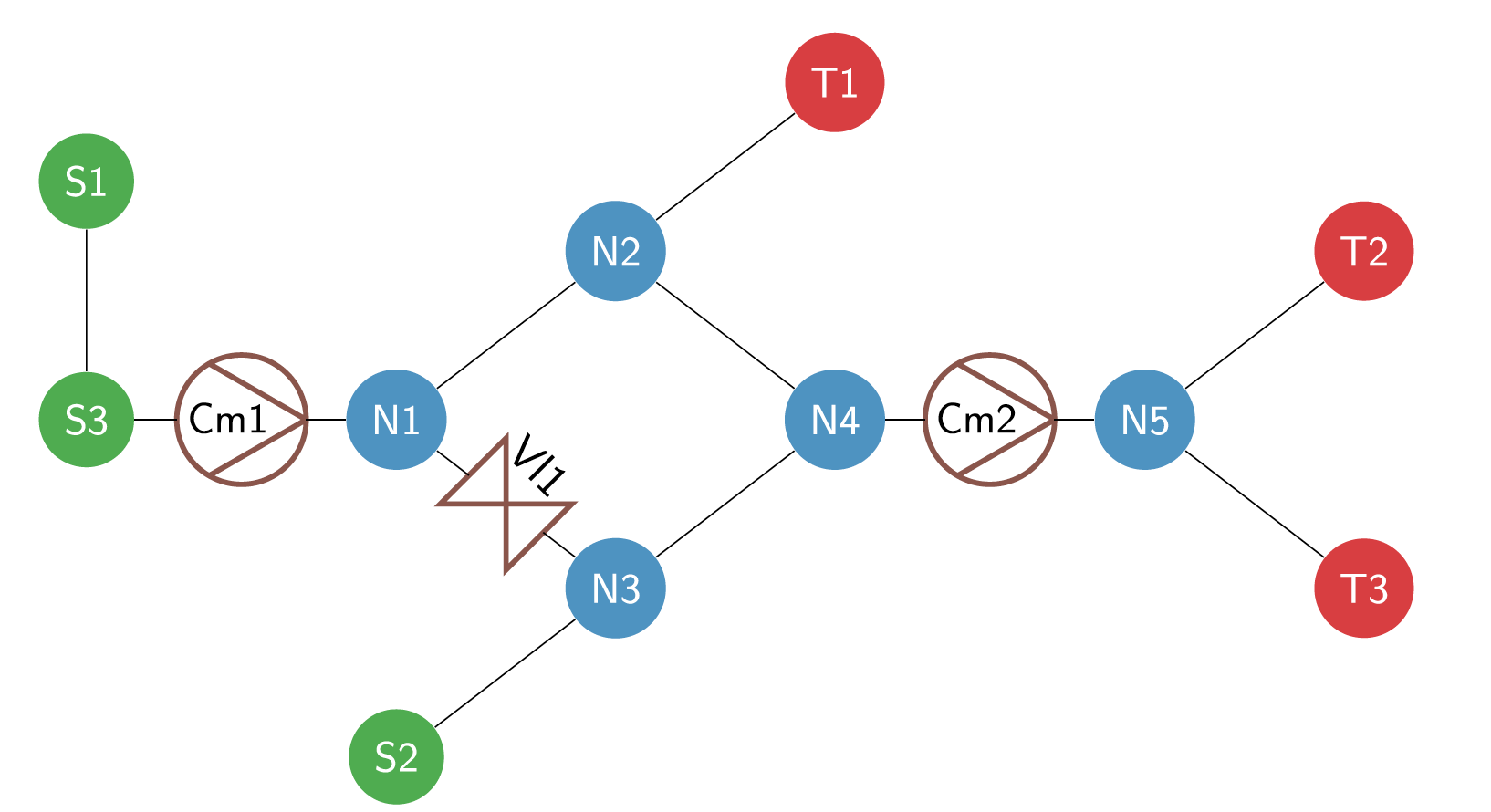}
 \caption{The GasLib-11 instance (figure is taken from \citep{burlacu2019}). Here, $\mathcal{N}_{source} = \{\text{S1, S2, S3}\}$, $\mathcal{N}_{sink} = \{\text{T1, T2, T3}\}$, $\mathcal{N}_{inner} = \{\text{N1, N2, N3, N4, N5}\}$, $\mathcal{C} = \{\text{Cm1, Cm2}\}$, $\mathcal{R} = \{\text{Vl1}\}$, $\mathcal{P}$ is the set of unlabeled~arcs.}
 \label{figure:gaslib11}
\end{figure}

In this paper, we consider an isothermal and stationary natural gas network, which is a common industry practice \citep{Sanchez2009Tabu}. We denote this network by $\mathcal{G} = (\mathcal{N}, \mathcal{E})$, where $\mathcal{N}$ denotes the set of nodes, and  $\mathcal{E}$ denotes the set of edges. Figure \ref{figure:gaslib11} shows an example of a small-scale gas network (taken from \cite{burlacu2019}). The set $\mathcal{N}$ consists of the set $\mathcal{N}_{source}$ of source nodes, the set $\mathcal{N}_{inner}$ of inner nodes and the set $\mathcal{N}_{sink}$ of sink nodes. 
We let $\mathcal{E} = \mathcal{P} \cup \mathcal{A}$, where $\mathcal{P}$ and $\mathcal{A}$ denote the set of passive and active elements, respectively. In our notation, set~$\mathcal{P}$ consists of the set of pipes, resistors and connections (i.e., short pipes). For simplicity, we refer to the elements of this set as pipes when necessary. We split $\mathcal{A}$ into the set $\mathcal{C}$ of compressors, the set $\mathcal{V}$ of control valves, and the set $\mathcal{R}$ of regular valves. For each edge $(i,j) \in \mathcal{E}$, we assume that the direction of the gas flow is arbitrary, e.g., it is positive if gas flows from node $i$ to node $j$, and negative otherwise. We restrict ourselves to a simplified compressor station, which consists of a bidirectional single compressor unit (i.e., compressor). We assume that all nodes are located horizontally (i.e., arcs are at the same height), and all arcs are cylindrical. For each node $i \in \mathcal{N}$, we define $\Delta^+(i) :=\{j \in \mathcal{N}: (i,j) \in \mathcal{A}\}$ as the set of outgoing neighbors and $\Delta^-(i):=\{j \in \mathcal{N}: (j,i) \in \mathcal{A}\}$ as the incoming neighbors of node $i$. 
Also, we respectively define $\Delta_{\mathcal{C}}^+(i) = \{j \in \mathcal{N}, (i,j) \in \mathcal{C}\}$ and $\Delta_{\mathcal{C}}^-(i)  = \{j \in \mathcal{N}, (j,i) \in \mathcal{C}\}$ as the sets of outgoing and incoming neighbors of node $i$, which are connected to this node via a compressor. 
We denote the set of gas stores linked to node~$i$ as $\mathcal{S}(i)$. 
Finally, set $\mathcal{T}$ represents the set of periods within the planning horizon.

\subsection{Gas Network Modeling} \label{ss:GasNetworkComponents}

In this section, we explain the physical and operational constraints of 
network elements. For more detailed explanations on the gas network modeling, we refer the reader to the book \cite{Koch2015}. 
For period $t \in \mathcal{T}$, each node $i \in \mathcal{N}$ is associated with a nonnegative pressure variable $p_{it}$, which is restricted between pressure bounds denoted by $\munderbar{p}_i$ and $\bar{p}_i$, respectively. Moreover, the gas load parameter at node $i \in \mathcal{N}$ is denoted by $q_{it}$ for period $t \in \mathcal{T}$. Note that a positive gas load (e.g., supply) is assigned for each $i \in \mathcal{N}_{source}$ and a negative gas load (e.g., demand) is assigned for each $i \in \mathcal{N}_{sink}$. We also let $q_{it} = 0$ for each $i \in \mathcal{N}_{inner}$. We denote the minimum and maximum store levels by $\munderbar{s}_j$ and $\bar{s}_j$, respectively. For each period $t \in \mathcal{T}$, the nonnegative variable $s_{kt} \in [\munderbar{s}_k, \bar{s}_k]$ represents the amount of gas (the store level) at store $k \in \mathcal{S}(i)$ linked to node $i \in \mathcal{N}$, whereas the gas injection and withdrawal amounts at store $k$ with predetermined rates $\eta_k^+$ and $\eta_k^-$ are denoted by $s_{kt}^+$ and $s_{kt}^-$, respectively. The gas withdrawals from store $k$ are associated with a positive supply cost $C^{wd}_k$. For each period $t \in \mathcal{T}$, the variable $f_{ijt}$ represents the gas flow through edge $(i,j) \in \mathcal{E}$. The flow variable $f_{ijt}$ is between allowable lower and upper flow limits denoted as $\munderbar{f}_{ij}$ and $\bar{f}_{ij}$, respectively. For each active element $(i,j) \in \mathcal{A}$, we associate a control variable $x_{ijt} \in \{0,1\}$, which takes the value $1$ if active element $(i,j)$ is in the operating state, and $0$ otherwise. For period $t \in \mathcal{T}$, each compressor $(i,j) \in \mathcal{C}$ is associated with binary start-up and shut-down variables $u_{ijt}$ and $v_{ijt}$, respectively. Whenever compressor $(i,j)$ is turned on from the off state, the positive start-up cost $C^{up}_{ij}$ is incurred.

Gas transportation in pipes is governed by the Euler equations, which are a system of hyperbolic partial differential equations (see, e.g., \cite{Osiadacz1987}). Finding a solution to this complex system is very hard. Following the common practice, we simplify this system under steady-state conditions 
to model the pressure losses in pipes due to internal friction. Under steady-state conditions, gas transportation through each pipe $(i,j) \in \mathcal{P}$ is governed by the well-known nonlinear and nonconvex Weymouth equation:
\begin{equation} \label{eq:weymouthPipe}
    p_{it}^2 - p_{jt}^2  = w_{ij} |f_{ijt}|f_{ijt} \qquad  (i,j) \in \mathcal{P}, t \in \mathcal{T}.
\end{equation} 
Here, the resistance coefficient $w_{ij}$ is often called as the Weymouth coefficient. We note that short pipes are artificial network elements that do not induce pressure losses, i.e., $w_{ij} = 0$. Also, resistors are modeled with relatively small resistance coefficients.

In gas networks, active elements control the gas flow and the pressures of their adjacent nodes. For example, compressors are used to increase the pressure of the incoming gas and transport it over long distances whereas control valves are used to reduce the pressure between its incident nodes. For an operating active element, the compression (pressure) ratio $p_{jt}/p_{it}$ belongs to the interval $[\munderbar{a}_{ij}, \bar{a}_{ij}]$ if the direction of the gas flow is positive, and $p_{it}/p_{jt} \in [\munderbar{a}_{ij}, \bar{a}_{ij}]$ otherwise. 
Whenever the active element is closed, nodal pressures at adjacent nodes are decoupled due to blocked gas flow, e.g., $f_{ij} = 0$. 
We summarize the operational conditions for active elements as follows: 
\begin{subequations}  \label{eq:activeElementConditions}
\begin{align}
\hspace{0.5em} 
x_{ijt} = 0  \quad &\implies \quad f_{ijt} = 0 &(&i,j) \in \mathcal{A}, t \in \mathcal{T} \label{eq:activeElementConditions_1}\\
\hspace{0.5em} 
f_{ijt} > 0, x_{ijt} = 1 \quad &\implies \quad\frac{p_{jt}}{p_{it}} \in [\munderbar{a}_{ij}, \bar{a}_{ij}] &(&i,j) \in \mathcal{A}, t \in \mathcal{T}\\
\hspace{0.5em} f_{ijt} < 0, x_{ijt} = 1 \quad &\implies \quad \frac{p_{it}}{p_{jt}} \in [\munderbar{a}_{ij}, \bar{a}_{ij}] &(&i,j) \in \mathcal{A}, t \in \mathcal{T}.
\end{align}
\end{subequations}

We assume that gas flows are bidirectional on each active element $(i,j) \in \mathcal{A}$, i.e.,  
$\munderbar{f}_{ij} < 0 < \bar{f}_{ij}$. 
Also, we let $1 < \munderbar{a}_{ij} < \bar{a}_{ij}$ for each compressor $(i,j) \in \mathcal{C}$, which results in a pressure increase at node $i$ (resp. node $j$) in case of a positive (resp. negative) gas flow. Similarly, we have $0 < \munderbar{a}_{ij} < \bar{a}_{ij} \le 1$ for each control valve $(i,j) \in \mathcal{V}$. We assume that there is no pressure loss in regular valves, e.g., $\munderbar{a}_{ij} =\bar{a}_{ij} = 1$ for $(i,j) \in \mathcal{R}$. 
 
We assume that each compressor $(i,j) \in \mathcal{C}$ 
associated with a special type of compressor drive (i.e., gas turbines), which provides the necessary power to compress the gas. Whenever compressor $(i,j)$ is operating, the turbines consume the gas directly taken from the network as their power source, and the positive fuel consumption cost denoted by $C^{fc}_{ij}$ is incurred. We denote the fuel consumption (i.e., gas loss) incurred by each compressor $(i,j)$ with the variable $l_{ijt}$ for each period $t \in \mathcal{T}$. The loss variables depend on the direction of the gas flow, some physical parameter $H_{ij}$ and isentropic exponent $\kappa$. 
If the gas flows from node~$i$ to node~$j$, 
the loss variable depends on the compression ratio $p_{jt}/p_{it}$, and $p_{it}/p_{jt}$ in case of a negative gas flow. 
We summarize these operational conditions for compressors as follows:
\begin{subequations}   \label{eq:CompressorConditions}
\begin{align}
\hspace{0.5em}  x_{ijt} = 0 \quad &\implies \quad f_{ijt} = 0, l_{ijt} = 0 &(&i,j) \in \mathcal{C}, t \in \mathcal{T} \label{eq:CompressorConditions_1}\\
\hspace{0.5em}
  f_{ijt} > 0, x_{ijt} = 1 \quad &\implies \quad 
l_{ijt} = H_{ij}((\frac{p_{jt}}{p_{it}})^{\kappa} -1) f_{ijt} &(&i,j) \in \mathcal{C}, t \in \mathcal{T} \label{eq:CompressorConditions-II}\\
\hspace{0.5em}  
f_{ijt} < 0, x_{ijt} = 1 \quad &\implies \quad 
l_{ijt} = H_{ij}((\frac{p_{it}}{p_{jt}})^{\kappa} -1) f_{ijt} 
&(&i,j) \in \mathcal{C}, t \in \mathcal{T} \label{eq:CompressorConditions-III}.
\end{align}
\end{subequations}
Finally, we consider the switching conditions for each active element $(i,j) \in \mathcal{A}$. The minimum up-time restrictions ensure that an active element $(i,j)$ remains open for at least $MU_{ij}$ many periods if it is activated. Similarly, the minimum down-time restrictions impose that the network element remains closed for at least $MD_{ij}$ many periods whenever it is turned down. The start-up and shut-down decisions are also coupled with switching decisions for each compressor $(i,j) \in \mathcal{C}$. We summarize the switching conditions for active elements as follows:
\begin{subequations}  \label{eq:activeTechnicalConditions}
\begin{align} 
&\hspace{0.5em}  x_{ijt} - x_{ij(t-1)} \le x_{ijt'}   \qquad &(&i,j) \in \mathcal{A}, t \in \mathcal{T} , t' \in \{t+1, \dots, t + MU_{ij} -1\}  \\
&\hspace{0.5em}  x_{ij(t-1)} - x_{ij} \le 1 - x_{ijt'}    &(&i,j) \in \mathcal{A}, t \in \mathcal{T}, t' \in \{t+1, \dots, t + MD_{ij} -1\}  \\
&\hspace{0.5em} x_{ijt}-x_{ij(t-1)} = u_{ijt} - v_{ijt}   &(&i,j) \in \mathcal{C}, t \in \mathcal{T}\\
&\hspace{0.5em} u_{ijt} + v_{ijt} \le 1  &(&i,j) \in \mathcal{C}, t \in \mathcal{T}.
\end{align}
\end{subequations}
We note that the switching decisions of compressors are analogous to commitment decisions in the unit commitment problem in power systems (see, e.g., \citep{Tuncer2023, Okumusoglu2024}). Finally, we assume that the initial state of each active element is off, i.e., (control) valves are closed and compressors are deactivated.
\subsection{MINLP Formulation}  \label{ss:MINLPFormulation} 
In this section, we present the MINLP formulation of the multi-period gas transportation and storage problem. 
Observe that one type of nonlinearity is in the form of $p_i^2 -p_j^2$ in equation \eqref{eq:weymouthPipe}. Also, note that these nonnegative pressure variables do not linearly appear elsewhere. Thus, we can easily capture this nonlinearity by introducing a squared pressure (potential) variable $\pi_{it} \in [\munderbar{\pi}_i, \bar{\pi}_i]$ for each node $i \in \mathcal{N}$ and each period $t \in \mathcal{T}$, where $\bar{\pi}_i = \bar{p}_i^2$ and $\munderbar{\pi}_i = \munderbar{p}_i^2$. 
This reformulation leads to the following nonconvex equations over two finite intervals: 
\begin{subequations}  \label{eq:weymouthPIpe}
\begin{align}
    &\hspace{0.5em} \pi_{it} - \pi_{jt}  = w_{ij} |f_{ijt}|f_{ijt} &(&i,j) \in \mathcal{P}, t \in \mathcal{T} \label{eq:weymouth1}\\
    &\hspace{0.5em} \munderbar{\alpha}_{ij} \le f_{ijt} \le \bar{\alpha}_{ij} &(&i,j) \in \mathcal{P}, t \in \mathcal{T} \label{eq:weymouth2}\\
    &\hspace{0.5em} \munderbar{\beta}_{ij} \le \pi_{it} - \pi_{jt} \le \bar{\beta}_{ij} &(&i,j) \in \mathcal{P}, t \in \mathcal{T}, \label{eq:weymouth3}
    \end{align}
\end{subequations}
where
\begin{align*}
    \munderbar{\alpha}_{ij} := &\max \bigg \{\munderbar{f}_{ij}, -\sqrt{\frac{\bar{\pi}_{j} - \munderbar{\pi}_i}{w_{ij}}} \bigg \} 
    &\bar{\alpha}_{ij} := &\min \bigg \{\bar{f}_{ij}, \sqrt{\frac{\bar{\pi}_{i} - \munderbar{\pi}_j}{w_{ij}}} \bigg \} \\
     \munderbar{\beta}_{ij} := &\max \bigg \{-w_{ij}(\munderbar{f}_{ij})^2, \munderbar{\pi}_{i} - \bar{\pi}_j\bigg \}
   &\bar{\beta}_{ij} := &\min \bigg \{w_{ij}(\bar{f}_{ij})^2, \bar{\pi}_{i} - \munderbar{\pi}_j\bigg \}.
\end{align*}
In our gas network modeling, we completely omit nonnegative pressure variables. 
They can be easily computed (if necessary) once the potential variables are obtained. For each active element $(i,j) \in \mathcal{A}$, we let $\munderbar{r}_{ij} = (\munderbar{a}_{ij})^2$, $\bar{r}_{ij} = (\bar{a}_{ij})^2$, and $\hat{\kappa} = \kappa / 2$ for the fractional exponent used in modeling gas losses induced by compressors to be consistent with this variable transformation.

Now, we are ready to present the mathematical programming formulation of the considered problem with the aforementioned variable and bound changes: 
\begin{subequations} \label{eq:MINLP}
\begin{align}
    \min &\hspace{0.5em} g(s, l, \pi) \label{eq:Gas_objFunction} \\
    \mathrm{s.t.} &\hspace{0.5em} \text{For each } t \in \mathcal{T}: \notag \\
        &\hspace{0.5em}  q_{it}  =  \sum_{k \in \mathcal{S}(i)} (s_{kt}^+ - s_{kt}^-) +  \sum_{k \in \Delta^+(i)} f_{ijt} - \sum_{j \in \Delta^-(i)} f_{jit}  + \sum_{j \in \Delta^{-}_\mathcal{C}(i)} l_{ijt}  -  \sum_{j \in \Delta^{+}_\mathcal{C}(i)} l_{jit} &i& \in \mathcal{N}
        \label{eq:MINLP_flowBalance-gas} \\
        &\hspace{0.5em} s_{kt} = s_{k(t-1)} + \eta_j^+ s_{kt}^+ -  \eta_j^- s_{kt}^-  &i& \in \mathcal{N}, k \in \mathcal{S}(i)
        \label{eq:MINLP_balance-store}\\
        &\hspace{0.5em} s_{k|\mathcal{T}|} = s_{k0}   &i& \in \mathcal{N}, k \in \mathcal{S}(i)
            \label{eq:MINLP_endInventory}\\
        &\hspace{0.5em} \munderbar{\pi}_{i} \le  \pi_{it} \le  \bar{\pi}_{i}    &i& \in \mathcal{N}
        \label{eq:pressureBounds_connectives}\\ 
        &\hspace{0.5em} \munderbar{f}_{ij} \le  f_{ijt} \le \bar{f}_{ij}   &(&i,j) \in \mathcal{E} \setminus \mathcal{P}
        \label{eq:gasFlowBounds_connectives}\\
        &\hspace{0.5em} f_{ijt} > 0, x_{ijt} = 1 \quad \implies \quad\frac{\pi_{jt}}{\pi_{it}} \in [\munderbar{r}_{ij}, \bar{r}_{ij}] &(&i,j) \in \mathcal{A}
        \label{eq:MINLP_activeElementRelative1}\\
        &\hspace{0.5em} f_{ijt} < 0, x_{ijt} = 1 \quad \implies \quad \frac{\pi_{it}}{\pi_{jt}} \in [\munderbar{r}_{ij}, \bar{r}_{ij}] &(&i,j) \in \mathcal{A}
        \label{eq:MINLP_activeElementRelative2} \\
        &\hspace{0.5em}  f_{ijt} > 0, x_{ijt} = 1 \quad \implies \quad l_{ijt} = H_{ij}((\frac{\pi_{jt}}{\pi_{it}})^{\hat{\kappa}} -1) f_{ijt} &(&i,j) \in \mathcal{C}
        \label{eq:MINLP_compressorLoss1}\\
        &\hspace{0.5em}  f_{ijt} < 0, x_{ijt} = 1 \quad \implies \quad l_{ijt} = H_{ij}((\frac{\pi_{it}}{\pi_{jt}})^{\hat{\kappa}} -1) f_{ijt} &(&i,j) \in \mathcal{C}
        \label{eq:MINLP_compressorLoss2}\\
        &\hspace{0.5em} x_{ijt} \in \{0,1\}  &(&i,j) \in \mathcal{A}
        \label{eq:switchingActive}\\
        &\hspace{0.5em} u_{ijt}, v_{ijt} \in \{0,1\} &(&i,j) \in \mathcal{C}
        \label{eq:startUpDownCompressor}\\
        &\hspace{0.5em} \eqref{eq:activeElementConditions_1}, \eqref{eq:CompressorConditions_1}, \eqref{eq:activeTechnicalConditions}, \eqref{eq:weymouthPIpe}. \nonumber
\end{align} 
\end{subequations}
In this formulation, the objective function $g$ in \eqref{eq:Gas_objFunction} is some linear objective function, e.g.,  
the total fuel consumption and start-up costs, and the supply cost of gas stores. Constraint \eqref{eq:MINLP_flowBalance-gas} models natural gas flow balance by the conservation of mass at each node. Constraints \eqref{eq:MINLP_balance-store} and \eqref{eq:MINLP_endInventory} correspond to the inventory balance equation for the amount of gas stored and the ending storage levels at each store. Constraints \eqref{eq:pressureBounds_connectives} and \eqref{eq:gasFlowBounds_connectives} are the operational limit constraints for potential and gas flow variables, respectively. Constraints \eqref{eq:MINLP_activeElementRelative1} and \eqref{eq:MINLP_activeElementRelative2} are the operational constraints for the active elements whereas constraints \eqref{eq:MINLP_compressorLoss1} and \eqref{eq:MINLP_compressorLoss2} are the fuel consumption equations for compressors. Constraints \eqref{eq:switchingActive} and \eqref{eq:startUpDownCompressor} represent the binary restrictions for control decisions for active arcs, and start-up and shut-down variables for compressors, respectively. We note that problem \eqref{eq:MINLP} is (highly) nonlinear and nonconvex due to 
constraints \eqref{eq:weymouth1} and \eqref{eq:MINLP_activeElementRelative1}-\eqref{eq:MINLP_compressorLoss2}.

\section{Main Results} \label{section:mainResults}

In this section, we will present our main results. We first introduce the notation used in the remainder of this paper in Section \ref{subsection:notation}. We present our convexification results on the feasible set for pipes in Section \ref{ss:feasible_sets_pipes}, and the feasible set for compressors in Section \ref{ss:feasible_sets_compressors}. More precisely, we give the convex hull descriptions of these two nonconvex sets and show that they are second-order cone representable (Theorem \ref{theorem:convexHull_pipes}) and power cone representable (Theorem \ref{theorem:convexHull_compressors_union}), respectively. To compare our convexification approaches, we provide the polyhedral and second-order outer-approximations known from the literature for pipes. We also propose a second-order cone approximation of the feasible set of compressors (Theorem~\ref{theorem:compressorOuterApprox}) that effectively utilizes the state-of-the-art solvers.

\subsection{Notation} \label{subsection:notation} For a set $X \subseteq \mathbb{R}^n$, we denote its convex hull, set of extreme points and interior as $\text{conv}(X)$, $\textup{extr}(X)$ and $\text{int}(X)$, respectively. We also denote the projection of $X$ onto the space of $x$-variables as $\text{proj}_x(X)$. We define 
the $3$-dimensional power cone parameterized by $a \in [0,1]$ as
\[{\bf P}_{a} := \{x \in \mathbb{R}^3: x_1^a x_2^{1-a} \ge |x_3|, \ x_1, x_2 \ge 0\}.\]
For a convex set $Y \subseteq \mathbb{R}^n$, we say that $Y$ is conic representable if $Y = \{ x \in \mathbb{R}^n: \exists w \in \mathbb{R}^m: T(x, w) \succeq_{K} 0\}$, 
where $K$ is a proper cone and $T(.)$ is an affine mapping. In particular, we will call that $Y$ is polyhedrally representable/SOCr/POWr if $K$ is a polyhedral/second-order/power cone. For a compact conic representable set $P := \{x \in \mathbb{R}^n : Ax \preceq_K~b\}$ and $z \in \{0, 1\}$,  we  denote the set $\{x \in \mathbb{R}^n: Ax \preceq_K b z\}$ as $ P z$. 

Consider a set of the form $\{(w,x,y) \in \mathbb{R}^3: w= xy, l_x \le x \le u_x, l_y \le y \le u_y\}$ with finite variable bounds. We denote the convex hull of this set 
using McCormick envelopes (see, \cite{McCormick1976}) as follows: 
\begin{align*}
\mathcal{M}(l_x,  u_x, l_y, u_y) = \{ (w,x,y) \in \mathbb{R}^3: & w  \ge \max \{l_y x + l_x y - l_x l_y,  u_y x + u_x y - u_x u_y\} \\
\ & w \le \min \{u_y x + l_x y - l_x u_y, l_y x + u_x y - u_x l_y\} 
\}.
\end{align*} 
In the remainder of this paper, we omit $t$ indices for notational brevity.

\subsection{Convexification of the feasible set of a pipe} \label{ss:feasible_sets_pipes}
In this section, we will focus on the feasible set of a specific  pipe $(i,j) \in \mathcal{P}$, which satisfies the constraints  in~\eqref{eq:weymouthPIpe}.  Assume that gas flow bounds satisfy the relations $\munderbar{\alpha}_{ij} \le 0 \le \bar{\alpha}_{ij}$ and  $\munderbar{\beta}_{ij} \le 0 \le \bar{\beta}_{ij}$. 
We define the following polytope 
\[\mathcal{B}_{ij} := \{(f_{ij}, \pi_i, \pi_j) \in \mathbb{R}^3 : \munderbar{\alpha}_{ij} \le f_{ij} \le  \bar{\alpha}_{ij},  \ \munderbar{\beta}_{ij} \le \pi_i - \pi_j \le  \bar{\beta}_{ij},  \ \munderbar{\pi}_{i} \le \pi_i \le  \bar{\pi}_{i},  \ \munderbar{\pi}_{j} \le \pi_j \le  \bar{\pi}_{j}\},\]
and the   following (bounded) nonconvex set: 
\begin{align*} \label{set:WeymouthSet}
    \mathcal{X}_{ij} := \{(f_{ij}, \pi_{i}, \pi_{j}) \in \mathcal{B}_{ij}: \pi_{i} - \pi_{j}  = w_{ij} |f_{ij}|f_{ij} \}.
\end{align*}
We refer to the set $\mathcal{X}_{ij}$ as the \textit{Weymouth set}, and our aim is to find its convex hull description in this section.  
It will be convenient to introduce the following sets
\[
    \mathcal{X}^+_{ij} := \{(f_{ij}, \pi_{i}, \pi_{j}) \in \mathcal{B}^+_{ij}: \pi_{i} - \pi_{j}  = w_{ij} f_{ij}^2 \}
    \text{ and }
     \mathcal{X}^-_{ij} := \{(f_{ij}, \pi_{i}, \pi_{j}) \in \mathcal{B}^-_{ij}: \pi_{j} - \pi_{i}  = w_{ij} f_{ij}^2 \},
\]
where 
\begin{align*}
    &\mathcal{B}^+_{ij} := \{(f_{ij}, \pi_i, \pi_j) \in \mathcal{B}_{ij}: 0 \le f_{ij} \le  \bar{\alpha}_{ij},  \ 0 \le \pi_i - \pi_j \le  \bar{\beta}_{ij}\}
    \\
    &\mathcal{B}^-_{ij} := \{(f_{ij}, \pi_i, \pi_j) \in \mathcal{B}_{ij}: \munderbar{\alpha}_{ij} \le f_{ij} \le  0,  \ \munderbar{\beta}_{ij} \le \pi_i - \pi_j \le  0 \}.  
\end{align*}
Note that the relation $ \mathcal{X}_{ij} = \mathcal{X}^+_{ij} \cup \mathcal{X}^-_{ij}$ is trivially satisfied.
%
\subsubsection{Second-order cone representable convex hull} \label{ss:socr_convex_hull_pipes}
Our first  main result is the following theorem, which states that the convex hull of the Weymouth set is SOCr:
\begin{theorem} \label{theorem:convexHull_pipes}
 $\textup{conv}(\mathcal{X}_{ij})$ is SOCr.
\end{theorem}
The proof strategy of Theorem~\ref{theorem:convexHull_pipes} is to show that both $\textup{conv}(\mathcal{X}_{ij}^+)$ and $\textup{conv}(\mathcal{X}_{ij}^-)$ are SOCr sets.

\begin{proposition}\label{prop:convexHull_pipes_Pos}
     We have $\textup{conv}(\mathcal{X}_{ij}^+) = \{ (f_{ij}, \pi_{i}, \pi_{j}) \in \mathcal{ \check B}^+_{ij} :  \pi_{i} - \pi_{j}  \ge w_{ij}  f_{ij}^2 \} $, where $\mathcal{ \check B}^+_{ij}$ is a polytope. Hence, $\textup{conv}(\mathcal{X}_{ij}^+)$ is SOCr.
\end{proposition} 
\begin{proof}
    Due to \cite{tawarmalani2013decomposition}, we have that 
    $$\textup{conv}(\mathcal{X}_{ij}^+) = 
    \textup{conv}(\{ (f_{ij}, \pi_{i}, \pi_{j}) \in \mathcal{B}^+_{ij}: \pi_{i} - \pi_{j}  \ge w_{ij}  f_{ij}^2  \}) \cap
    \textup{conv}(\{ (f_{ij}, \pi_{i}, \pi_{j}) \in \mathcal{B}^+_{ij}: \pi_{i} - \pi_{j}  \le w_{ij}  f_{ij}^2  \}).$$
    Notice that the first set is the intersection of the polytope $\mathcal{B}^+_{ij}$ with the convex quadratic inequality $\pi_{i} - \pi_{j}  \ge w_{ij}  f_{ij}^2$, hence, it is already SOCr. 
    The second set is the convex hull of the intersection of the polytope $\mathcal{B}^+_{ij}$ with the reverse convex constraint $\pi_{i} - \pi_{j}  \le w_{ij}  f_{ij}^2$, which is known to be a polytope due to~\cite{hillestad1980linear}, hence it is also SOCr. Let us denote this polytope by $\mathcal{ \check B}^+_{ij} $ which is a subset of $\mathcal{ B}^+_{ij} $, by construction.
    Finally, as the intersection of two SOCr sets, $\textup{conv}(\mathcal{X}_{ij}^+)$ is SOCr.
\end{proof}
We note that the extreme points of the polytope $\mathcal{ \check B}^+_{ij}$ can be explicitly identified by a simple procedure (see Algorithm~\ref{alg:weymouthExtremesPlus} in Appendix~\ref{app:weymouth}).

 The proof of the following proposition is similar to that of Proposition~\ref{prop:convexHull_pipes_Pos}, hence, it is omitted.
\begin{proposition}\label{prop:convexHull_pipes_Neg}
     We have $\textup{conv}(\mathcal{X}_{ij}^-) = \{ (f_{ij}, \pi_{i}, \pi_{j}) \in \mathcal{ \check B}^-_{ij} :  \pi_{j} - \pi_{i}  \ge w_{ij}  f_{ij}^2 \} $, where $\mathcal{ \check B}^-_{ij}$ is a polytope. Hence, $\textup{conv}(\mathcal{X}_{ij}^-)$ is SOCr.
\end{proposition} 

We are now ready to prove Theorem~\ref{theorem:convexHull_pipes}.
\begin{proof}[Proof of Theorem~\ref{theorem:convexHull_pipes}]
    Due to Propositions~\ref{prop:convexHull_pipes_Pos} and~\ref{prop:convexHull_pipes_Neg}, we deduce that $\textup{conv}(\mathcal{X}_{ij}^+)$ and $\textup{conv}(\mathcal{X}_{ij}^-)$ are SOCr. Since these sets are also compact and  satisfy the relation $ \mathcal{X}_{ij} = \mathcal{X}^+_{ij} \cup \mathcal{X}^-_{ij}$, the convex hull of their union is also SOCr \cite{ben2001lectures}.
\end{proof}
 Since the proof of Theorem~\ref{theorem:convexHull_pipes} is constructive, we are able to give an extended formulation for $\textup{conv}(\mathcal{X}_{ij})$ as below:
 \[ \textup{conv}(\mathcal{X}_{ij}) =  \left \{  \begin{matrix*} (f_{ij}, \pi_i, \pi_j) \in  \mathbb{R}^3, \\ (f_{ij}^+, \pi_i^+, \pi_j^+) \in \mathbb{R}^3, \\ (f_{ij}^-, \pi_i^-, \pi_j^-) \in \mathbb{R}^3, \\ z_{ij} \in [0,1] \end{matrix*}  : \begin{matrix*}[l] 
(f_{ij}^+, \pi_i^+, \pi_j^+) \in  \mathcal{ \check B}^+_{ij} z_{ij}, \ z_{ij} (\pi^+_i - \pi^+_j) \ge w_{ij} (f_{ij}^+)^2 \\ 
(f_{ij}^-, \pi_i^-, \pi_j^-) \in  \mathcal{ \check B}^-_{ij} (1-z_{ij}), \ (1-z_{ij})(\pi^-_j - \pi^-_i) \ge w_{ij} (f_{ij}^-)^2 \\ 
 f_{ij} = f^+_{ij} + f^-_{ij}, \quad \pi_i = \pi^+_i +   \pi^-_i, \quad \pi_j = \pi^+_j + \pi^-_j \end{matrix*} \right \}.\]
Here, when we   have  the binary variable  $z_{ij}=1$, the  flow is positive and  the variables $(f_{ij}^+, \pi_i^+, \pi_j^+) $ belong to the set $\mathcal{ \check B}^+_{ij}$ while the variables $(f_{ij}^+, \pi_i^+, \pi_j^+) $  are all zero (the case with $z_{ij}=0$ is exactly the opposite). Although the above extended formulation gives the convex hull of the Weymouth set for a specific pipe, hence, the  variable $z_{ij}$ can be relaxed to a continuous variable, we need to enforce it to take a binary value in the gas network optimization problems.

\begin{figure}[H]
    \centering
 {\scalebox{0.90}{
\begin{tikzpicture}
 \draw[gray, opacity=0.5, dashdotted, xstep=0.5cm, ystep=0.5cm] (-2.75,-2.75) grid (2.75,2.75); 
    \draw[<->] (-2.75,0) -- (2.75,0) node[right] { $f_{ij}$}; 
    \draw[<->] (0,-2.75) -- (0,2.75) node[above] { $\pi_i - \pi_j$}; 
    \draw[scale=1, thick, domain=-2.25:2.25, smooth, variable=\x, color= blue] plot ({\x}, {0.4*\x^2}); 
    \node [circle,fill=blue!90,   dotted, scale=0.15, label={[label distance=-0.2cm]105:\scriptsize$ (\bar{\alpha}_{ij}, \bar{\beta}_{ij})$}] (a) at (2.25, 0.4*2.25*2.25)  {$\hspace{0.5em}$}; 
    \node [circle,fill=blue!90,   dotted, scale=0.15, label={[label distance=-0.2cm]-5:\scriptsize$ (\munderbar{\alpha}_{ij}, \munderbar{\beta}_{ij})$}] (a) at  (-2.25, -0.4*2.25*2.25) {$\hspace{0.5em}$}; 
    \draw[scale=1,   dotted,  domain=-2.25:-1.5, color=red] plot[id=x] function{2.706*x + 3.812}; 
    \draw[scale=1,   dotted,  domain=1.5:2.3, color=red] plot[id=x] function{2.235*x + -2.9505}; 
    \filldraw[blue, opacity=0.1,   dotted, xstep=0.5cm, ystep=0.5cm] (-2.25,-2.25*2.25*0.4) -- (-1.75, -1.75*1.75*0.4+0.3) -- (2.25, 2.25*2.25*0.4) -- (1.75, 1.75*1.75*0.4-0.3)  -- (-2.25,-2.25*2.25*0.4);
    \draw[red, densely dotted]  (-2.25,-2.25*2.25*0.4) -- (-1.75, -1.75*1.75*0.4+0.3) -- (2.25, 2.25*2.25*0.4) -- (1.75, 1.75*1.75*0.4-0.3)  -- (-2.25,-2.25*2.25*0.4);
\end{tikzpicture}}}
\caption{Polyhedral outer-approximation of  \(\mathcal{X}_{ij}\) in the $(f_{ij}, \pi_i - \pi_j)$-space.} \label{fig:polyOuterApprox}
\end{figure}
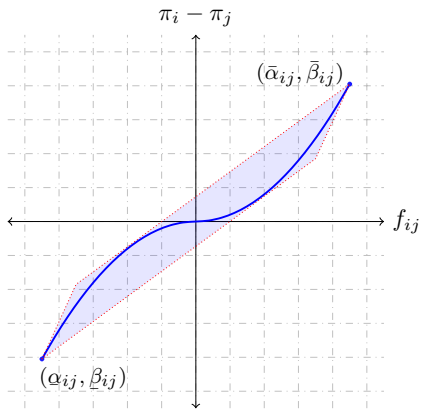

\subsubsection{Outer-approximations via different convex sets}  \label{ss:outer_approx_pipes}
Here, we introduce two convex outer-approximations of the Weymouth set taken from the literature that involve polyhedrally representable and SOCr sets. 

\paragraph{Polyhedral outer-approximation}  \label{ss:outer_approx_pipes_polyhedral}
We now provide a polyhedral approximation of the Weymouth set, which is adapted from~\cite{Pfetsch2015}. We use four hyperplanes to obtain the polyhedrally representable outer-approximation of this set. Our construction is geometrically illustrated in Figure \ref{fig:polyOuterApprox} in the space of $(f_{ij}, \pi_i - \pi_j)$-variables. To be precise, we first find the line passing through $(\munderbar{\alpha}_{ij}, \munderbar{\beta}_{ij})$ and tangent to the curve at point $(a_{ij}, w_{ij} a_{ij}^2)$ with $f_{ij}, \pi_i - \pi_j \ge 0$, where $a_{ij}$ is some constant. Similarly, we obtain the line passing through $(\bar{\alpha}_{ij}, \bar{\beta}_{ij})$ and tangent to this curve at point $(a_{ij}, -w_{ij} a_{ij}^2)$ with $f_{ij}, \pi_i - \pi_j \le 0$. For the remaining hyperplanes, we find two tangent lines to this curve at the points $(\munderbar{\alpha}_{ij}, \munderbar{\beta}_{ij})$ and  $(\bar{\alpha}_{ij}, \bar{\beta}_{ij})$, respectively. After some simple calculations, the polyhedral outer-approximation of the Weymouth set is formed by the set
\begin{align*}
     \mathcal{X}_{\text{poly}}:= \{(f_{ij}, \pi_i, \pi_j) \in \mathbb{R}^3:  &\pi_i - \pi_j \ge  w_{ij} a_{ij} (2f_{ij}-a_{ij}) , \ \pi_i - \pi_j \le  w_{ij} b_{ij} (b_{ij} - 2f_{ij}) \\ 
     &\pi_i - \pi_j \ge  w_{ij} \bar{\alpha}_{ij} (2f_{ij}-\bar{\alpha}_{ij}), \ \pi_i - \pi_j \le  w_{ij} \munderbar{\alpha}_{ij} (\munderbar{\alpha}_{ij}- 2f_{ij}) \},
\end{align*}
where $a_{ij} = (1- \sqrt{2})\munderbar{\alpha}_{ij}$ and  $b_{ij} = (1 - \sqrt{2}) \bar{\alpha}_{ij}$.

Note that other polyhedral outer-approximations of the Weymouth set such as piecewise linearizations (see, e.g., \cite{Pfetsch2015}) can also be applied. However, these formulations require using many binary variables for high-resolution solutions, even for a small-scale network. In this case, the global solver might fail to provide even locally optimal solutions as the size of the gas network increases due to the time limit. 
\paragraph{Second-order cone outer-approximation} \label{ss:outer_approx_pipes_socr}

A known SOCr relaxation of the Weymouth set is defined as follows \cite{Borraz2016SOCP} and used in the context of the gas expansion planning problem:
\begin{align*}
    \mathcal{X}_{\text{socr}}:= \{(f_{ij}, \pi_i, \pi_j) \in \mathbb{R}^3: \exists (\pi', z) \in \mathbb{R} \times [0,1] : \ & \munderbar{\alpha}_{ij} (1-z_{ij}) \le f_{ij} \le \bar{\alpha}_{ij} z_{ij} , \munderbar{\beta}_{ij} (1-z_{ij}) \le \pi_i - \pi_j \le \bar{\beta}_{ij} z_{ij} \\
    & \pi_{ij}' \ge w_{ij} f_{ij}^2, \ (\pi_{ij}', 2z_{ij}-1, \pi_i -\pi_j) \in \mathcal{M}(-1, 1, \munderbar{\beta}_{ij}, \bar{\beta}_{ij}) \}.
\end{align*}
Note that this set involves a second-order cone constraint and McCormick inequalities, but does not give the convex hull in general. 


\begin{remark} \label{remark:posFlows}
We note that in some GasLib test instances, flow directions are known to be positive. In such cases, the polyhedral and second-order cone outer-approximations can be formed by the sets $\mathcal{X}_{poly}~:=~\{(f_{ij}, \pi_i, \pi_j) \in \mathcal{B}_{ij}^+:  \pi_i - \pi_j \le \frac{\bar{\beta}_{ij}}{ \bar{\alpha}_{ij}}f_{ij}\}$ and $\mathcal{X}_{socr} := \{(f_{ij}, \pi_i, \pi_j) \in \check{\mathcal{B}}_{ij}^+: \pi_i - \pi_j \ge w_{ij} f_{ij}^2\}$, respectively. In our computational study, we will use these sets for the cases with positive gas flows. We also note that the set $\mathcal{X}_{socr}$ leads to the convex hull of the Weymouth set in the case of positive flows. Similar results also follow for both sets when the flow directions are negative.
\end{remark}

\subsection{Convexification of the feasible set of a compressor} \label{ss:feasible_sets_compressors}

In this section, we will focus on the feasible set of a specific compressor $(i,j) \in \mathcal{C}$, which satisfies the constraints  in~\eqref{eq:pressureBounds_connectives}-\eqref{eq:switchingActive}.   
We define the following polytopes with $1 < \munderbar{r}_{ij} < \bar{r}_{ij}$:
\begin{align*}
    \Pi_{ij} := \Big\{(\pi_i, \pi_j) \in [\munderbar{\pi}_{i} ,   \bar{\pi}_{i}]\times [ \munderbar{\pi}_{j} ,  \bar{\pi}_{j}] : \munderbar{r}_{ij} \le \frac{\pi_j}{\pi_i} \le \bar{r}_{ij} \Big\},  \ 
    \Pi_{ji} := \Big\{(\pi_i, \pi_j) \in [\munderbar{\pi}_{i} ,   \bar{\pi}_{i}]\times [ \munderbar{\pi}_{j} ,  \bar{\pi}_{j}] : \munderbar{r}_{ij} \le \frac{\pi_i}{\pi_j} \le \bar{r}_{ij}  \Big\},
\end{align*}
and the  (bounded) nonconvex set   $\mathcal{L}_{ij} := \mathcal{L}_{ij}^+ \cup \mathcal{L}_{ij}^- \cup \mathcal{L}_{ij}^0$, where
\begin{equation*}  \label{set:partitions}
\begin{split}
    &\hspace{0.5em} \mathcal{L}_{ij}^+ :=  \{ (\pi_i, \pi_j, f_{ij}, l_{ij}, x_{ij}) \in \mathbb{R}^4  \times \{1\}: l_{ij} = H_{ij}((\pi_{j}/\pi_{i})^{\hat{\kappa}} -1) f_{ij}, \ (\pi_i, \pi_j) \in \Pi_{ij},\ f_{ij} \in [0, \bar{f}_{ij}] \} \\
    &\hspace{0.5em} \mathcal{L}_{ij}^- :=  \{ (\pi_i, \pi_j, f_{ij}, l_{ij}, x_{ij}) \in \mathbb{R}^4  \times \{1\}: l_{ij} = H_{ij}((\pi_{i}/\pi_{j})^{\hat{\kappa}} -1) f_{ij}, \ (\pi_i, \pi_j) \in \Pi_{ji},\ f_{ij} \in [\munderbar{f}_{ij}, 0] \}\\
    &\hspace{0.5em} \mathcal{L}_{ij}^0 :=  \{ (\pi_i, \pi_j, f_{ij}, l_{ij}, x_{ij}) \in \mathbb{R}^2  \times \{0\}^3: \munderbar{\pi}_i \le \pi_i \le \bar{\pi}_i, \ \munderbar{\pi}_j \le \pi_j \le \bar{\pi}_j\}.
\end{split}
\end{equation*} 
We refer to the set $\mathcal{L}_{ij}$ as the \textit{compression set}, and our aim is to find its convex hull description in this section.



\subsubsection{Power cone representable convex hull} \label{ss:power_cone_convex_hull_compressors}
We state our second main result in the following theorem:
\begin{theorem} \label{theorem:convexHull_compressors_union}
$\textup{conv}(\mathcal{L}_{ij})$ is POWr. 
\end{theorem}
The proof strategy of Theorem~\ref{theorem:convexHull_compressors_union} is to show that the convex hull of extreme points of the sets  
$\mathcal{L}_{ij}^+ $, $\mathcal{L}_{ij}^- $ and $\mathcal{L}_{ij}^0 $
are POWr sets. Since the polytope $\mathcal{L}_{ij}^0$ is trivially POWr, we will only focus on sets $\mathcal{L}_{ij}^+ $ and $\mathcal{L}_{ij}^- $.
First, we present a preliminary result that will be used in the proof  
that will eventually give a complete extreme point characterization of the set    $\mathcal{L}_{ij}^+$ in Proposition~\ref{prop:extrPoints_compressors}. 
\begin{lemma} \label{lemma:extrPoints_1}
Let  $(\pi_i, \pi_j, f_{ij}, l_{ij}, x_{ij}) \in \textup{extr}(\mathcal{L}_{ij}^+)$. Then,  
\begin{enumerate}[(i)]
    \item $f_{ij} \in \{0, \bar{f}_{ij}\}$.
    \item  $(\pi_i, \pi_j) \notin \textup{int}(\Pi_{ij})$.
    \item Either $\pi_i \in\{\munderbar{\pi}_i, \bar{\pi}_i\}$  or $\pi_j \in\{\munderbar{\pi}_j, \bar{\pi}_j\}$.
\end{enumerate}
\end{lemma}
\begin{proof}
\

\begin{enumerate}[(i)]
\item
Suppose that there exists a point $y = (\pi_i, \pi_j, f_{ij}, l_{ij}, x_{ij}) \in \textup{extr}(\mathcal{L}^+_{ij})$ with $0 < f_{ij} < \bar{f}_{ij}$. Then, there exists a small enough $\epsilon > 0$  such that 
$    y'  = (\pi_i, \pi_j, (1+ \epsilon)f_{ij}, (1 + \epsilon)l_{ij}, x_{ij}) \in \mathcal{L}^+_{ij} $ and $
    y'' = (\pi_i, \pi_j, (1 - \epsilon)f_{ij}, (1 - \epsilon)l_{ij}, x_{ij}) \in \mathcal{L}^+_{ij}$.
This implies that $y \notin \textup{extr}(\mathcal{L}^+_{ij})$ since $y = \frac12 y' + \frac12 y''$, which leads to a contradiction.

\item
We first show that $(\pi_i, \pi_j) \notin \textup{int}(\Pi_{ij})$. In fact, suppose that there exists a point $y = (\pi_i, \pi_j, f_{ij}, l_{ij}, x_{ij}) \in \textup{extr}(\mathcal{L}^+_{ij})$ with  $(\pi_i, \pi_j) \in \textup{int}(\Pi_{ij})$.
Then, there exists a small enough $\epsilon > 0$  such that 
$
    y'  = ((1+\epsilon)\pi_i, (1+\epsilon)\pi_j, f_{ij}, l_{ij}, x_{ij}) \in \mathcal{L}^+_{ij} $ and $
    y''  = ((1-\epsilon)\pi_i, (1-\epsilon)\pi_j, f_{ij}, l_{ij}, x_{ij}) \in \mathcal{L}^+_{ij}$.
This implies that $y = \frac12 y' + \frac12 y''$, which leads to a contradiction.

\item
Suppose that there exists a point $y = (\pi_i, \pi_j, f_{ij}, l_{ij}, x_{ij}) \in \textup{extr}(\mathcal{L}^+_{ij})$ with  
$\munderbar{\pi}_{i} < \pi_i <  \bar{\pi}_{i}$ and $   \munderbar{\pi}_{j} < \pi_j <  \bar{\pi}_{j}$. Due to Part (ii), we must have $\frac{\pi_i}{\pi_j} \in \{\munderbar{r}_{ij}, \bar{r}_{ij}\}$. The rest of the proof is the same as the proof of Part (ii).
\end{enumerate}
\end{proof}
In Lemma \ref{lemma:extrPoints_1}, we identify three cases for an extreme point of the nonconvex set $\mathcal{L}^+_{ij}$. By using these results, we are able to characterize $\textup{extr}(\mathcal{L}_{ij}^+)$ through five disjunctions. We now present the following proposition on the characterization of the extreme points of the nonconvex set $\mathcal{L}_{ij}^+$.
\begin{proposition} \label{prop:extrPoints_compressors}
$\textup{extr}(\mathcal{L}_{ij}^+)$ is contained in  $\mathcal{U}_{ij}^{+1} \cup \mathcal{U}_{ij}^{+2} \cup \mathcal{U}_{ij}^{+3} \cup \mathcal{U}_{ij}^{+4}  \cup \mathcal{U}_{ij}^{+5} $, where
\begin{align*}
    & \mathcal{U}_{ij}^{+1} :=  \{ ( {\pi}_i, \pi_j, 0, 0, 1) \in \mathbb{R}^5:  (\pi_i,\pi_j) \in \Pi_{ij} \}\\
    & \mathcal{U}_{ij}^{+2} :=  \{ (\munderbar{\pi}_i, \pi_j, \bar{f}_{ij}, l_{ij}, 1) \in \mathbb{R}^5: l_{ij} = H_{ij} \bar{f}_{ij} ((\pi_{j}/\munderbar{\pi}_i)^{{\hat{\kappa}}} -1)  , \ \max\{\munderbar{\pi}_j ,  \munderbar{\pi}_i\munderbar{r}_{ij}  \} \le  \pi_j \le \min\{ \bar \pi_j,  \munderbar{\pi}_i\bar{r}_{ij}  \} \}\\
     & \mathcal{U}_{ij}^{+3} :=  \{ (\bar{\pi}_i, \pi_j, \bar{f}_{ij}, l_{ij}, 1) \in \mathbb{R}^5: l_{ij} = H_{ij} \bar{f}_{ij} ((\pi_{j}/\bar{\pi}_i)^{{\hat{\kappa}}} -1)  , \ \max\{\munderbar{\pi}_j , \bar{\pi}_i\munderbar{r}_{ij}   \} \le  \pi_j \le \min\{ \bar \pi_j,  \bar{\pi}_i\bar{r}_{ij}  \} \}\\
    & \mathcal{U}_{ij}^{+4} :=  \{ (\pi_i, \munderbar{\pi}_j, \bar{f}_{ij}, l_{ij}, 1) \in \mathbb{R}^5: l_{ij} = H_{ij}\bar{f}_{ij} ((\munderbar{\pi}_j/\pi_i)^{{\hat{\kappa}}} -1) , \ 
    \max\{\munderbar{\pi}_i , \munderbar{\pi}_j / \bar{r}_{ij}  \} \le  \pi_i \le \min\{ \bar \pi_i,  \munderbar{\pi}_j / \munderbar{r}_{ij}  \} \} \\
    & \mathcal{U}_{ij}^{+5} :=  \{ (\pi_i, \bar{\pi}_j, \bar{f}_{ij}, l_{ij}, 1) \in \mathbb{R}^5: l_{ij} = H_{ij}\bar{f}_{ij} ((\bar{\pi}_j/\pi_i)^{{\hat{\kappa}}} -1) , \ 
    \max\{\munderbar{\pi}_i , \bar{\pi}_j / \bar{r}_{ij}  \} \le  \pi_i \le \min\{ \bar \pi_i,  \bar{\pi}_j / \munderbar{r}_{ij}  \} \}.
\end{align*}
\end{proposition}

Proposition~\ref{prop:extrPoints_compressors} states that the set of extreme points of $\mathcal{L}_{ij}^+$ is given by the union of the polytope $\mathcal{U}_{ij}^{+1}$ and the bounded nonconvex sets $\mathcal{U}_{ij}^{+2}, \dots, \mathcal{U}_{ij}^{+5}$. Our aim is now to find the convex hull of each of these nonconvex sets. Let us continue by presenting the following lemmas that are useful in finding these convex hulls.

The following lemma follows from the fact that  for $n\in(0,1)$,  $u \le x^n$ is equivalent to $(x, 1, u) \in {\bf P}_{n}$ over $(u,x)\in\mathbb{R}_+^2$.
\begin{lemma} \label{lem:powerConeReprPlus}
    Consider the following nonconvex set $$\mathcal{V}^n(a,b,\munderbar{x}, \bar x) = \{(x,y) \in \mathbb{R}^2_+: y = a ( (x/b)^n -1), \ \munderbar{x} \le x \le \bar{x}\}, $$ where $n \in (0,1) $, $a >0$ and  $0<b\le\munderbar{x} \le \bar x$. 
    Then,
        \begin{align*}
        \textup{conv}( \mathcal{V}^n(a,b,\munderbar{x}, \bar x) ) = \{(x,y)\in \mathbb{R}^2_+: \exists u \in \mathbb{R}_+:\ & (x/b, 1, u) \in {\bf P}_{n},  y = a(u-1), \\
        &
        y \ge a\Big ( \frac{ (\bar{x}/b) ^n - (\munderbar{x}/b)^n}{\bar{x} - \munderbar{x}} \Big ) (x - \munderbar{x}) + a ( (\munderbar{x}/b)^n-1) \},
    \end{align*}
    is POWr.
\end{lemma}

The following lemma follows from the fact that  for $n\in(0,1)$,  $u \ge x^{-n}$ is equivalent to $(u,x,1) \in {\bf P}_{\frac{1}{1+n}}$ over $(u,x)\in\mathbb{R}_{++}^2$.
\begin{lemma} \label{lem:powerConeReprMinus}
    Consider the following nonconvex set $$\mathcal{W}^n(a,b,\munderbar{x}, \bar x) = \{(x,y) \in \mathbb{R}^2_+: y = a ( (b/x)^n -1), \ \munderbar{x} \le x \le \bar{x}\}, $$ where $n \in (0,1) $, $a >0$ and  $0\le\munderbar{x} \le \bar x<b$. 
    Then,
        \begin{align*}
        \textup{conv}( \mathcal{W}^n(a,b,\munderbar{x}, \bar x) ) = \{(x,y)\in \mathbb{R}^2_+: \exists u \in \mathbb{R}_+:\ & (u,x/b, 1) \in {\bf P}_{\frac{1}{1+n}},  y = a(u-1), \\
        &
        y \le a\Big ( \frac{ (b/\bar{x}) ^n - (b/\munderbar{x})^n}{\bar{x} - \munderbar{x}} \Big ) (x - \munderbar{x}) + a ( (b/\munderbar{x})^n-1) \},
    \end{align*}
    is POWr.
\end{lemma}

Now, we are ready to present the convex hull of the nonconvex sets $\mathcal{U}_{ij}^{+2}, \dots, \mathcal{U}_{ij}^{+5}$. In fact, the following proposition is an immediate consequence of Proposition~\ref{prop:extrPoints_compressors}, Lemma~\ref{lem:powerConeReprPlus} and Lemma~\ref{lem:powerConeReprMinus}:
\begin{proposition} \label{proposition:convexHull_compressors_U+_disj}
We have 
\begin{align*}
    & \textup{conv}(\mathcal{U}_{ij}^{+2}) =  \left\{ (\munderbar{\pi}_i, \pi_j, \bar{f}_{ij}, l_{ij}, 1) \in \mathbb{R}^5: 
    (\pi_i,l_{ij}) \in \textup{conv} \big( \mathcal{V}^{\hat \kappa}(  H_{ij} \bar f_{ij} , \munderbar{\pi}_i, 
    \max\{\munderbar{\pi}_j ,  \munderbar{\pi}_i\munderbar{r}_{ij}  \} , \min\{ \bar \pi_j,  \munderbar{\pi}_i\bar{r}_{ij}  \} ) \big) \right\}\\
     & \textup{conv}(\mathcal{U}_{ij}^{+3}) =   \left\{ (\bar{\pi}_i, \pi_j, \bar{f}_{ij}, l_{ij}, 1) \in \mathbb{R}^5: 
    (\pi_i,l_{ij}) \in \textup{conv} \big( \mathcal{V}^{\hat \kappa} ( 
     H_{ij} \bar{f}_{ij} , \bar{\pi}_i , \max\{\munderbar{\pi}_j , \bar{\pi}_i\munderbar{r}_{ij}   \} ,    \min\{ \bar \pi_j,  \bar{\pi}_i\bar{r}_{ij} \} ) \big) \right \}\\
    & \textup{conv}(\mathcal{U}_{ij}^{+4}) =   \left\{ (\pi_i, \munderbar{\pi}_j, \bar{f}_{ij}, l_{ij}, 1) \in \mathbb{R}^5: 
    (\pi_j,l_{ij}) \in \textup{conv} \big( \mathcal{W}^{\hat \kappa} (  H_{ij}\bar{f}_{ij} , \munderbar{\pi}_j  ,  
    \max\{\munderbar{\pi}_i , \munderbar{\pi}_j / \bar{r}_{ij}  \} , \min\{ \bar \pi_i,  \munderbar{\pi}_j / \munderbar{r}_{ij}  \} ) \big) \right\} \\
    & \textup{conv}(\mathcal{U}_{ij}^{+5}) =   \left\{ (\pi_i, \bar{\pi}_j, \bar{f}_{ij}, l_{ij}, 1) \in \mathbb{R}^5:
    (\pi_j,l_{ij}) \in \textup{conv} \big( \mathcal{W}^{\hat \kappa} (  H_{ij}\bar{f}_{ij} , \bar{\pi}_j  ,  
    \max\{\munderbar{\pi}_i , \bar{\pi}_j / \bar{r}_{ij}  \} , \min\{ \bar \pi_i,  \bar{\pi}_j / \munderbar{r}_{ij}  \} ) \big) \right\}.
\end{align*}
\end{proposition} 

Finally, we present our main result on the convex hull of the nonconvex set $\mathcal{L}_{ij}^+$:
\begin{proposition} \label{proposition:convexHull_compressors}
 $\textup{conv}(\mathcal{L}_{ij}^+)$ is POWr.
\end{proposition}
\begin{proof}
 Recall that  we have  $\textup{extr}(\mathcal{L}_{ij}^+) \subseteq \mathcal{U}_{ij}^{+1} \cup \mathcal{U}_{ij}^{+2} \cup \mathcal{U}_{ij}^{+3} \cup \mathcal{U}_{ij}^{+4}  \cup \mathcal{U}_{ij}^{+5} \subseteq \mathcal{L}_{ij}^+   $.
 Being a polytope, $\mathcal{U}_{ij}^1$ is trivially POWr. 
    On the other hand, $\textup{conv}(\mathcal{U}_{ij}^{+2}), \dots, \textup{conv}(\mathcal{U}_{ij}^{+5})$ are POWr sets due to  Proposition~\ref{proposition:convexHull_compressors_U+_disj}.  
Since $\textup{conv}(\mathcal{U}_{ij}^{+1}), \dots, \textup{conv}(\mathcal{U}_{ij}^{+5})$ are compact POWr sets and 
the set $\textup{conv}(\mathcal{L}_{ij}^+)$ is the convex hull of their union, the result follows due to~\cite{ben2001lectures}.
\end{proof}
Let us present a similar analysis for the extreme point characterization of the nonconvex set $\mathcal{L}_{ij}^-$. We start by presenting the following lemma. Note that the proof is omitted since it is similar to that of Lemma~\ref{lemma:extrPoints_1}. 

\begin{lemma} \label{lemma:extrPoints_1Neg}
Let $(\pi_i, \pi_j, f_{ij}, l_{ij}, x_{ij}) \in \textup{extr}(\mathcal{L}_{ij}^-)$. Then,
\begin{enumerate}[(i)]
    \item $f_{ij} \in \{\munderbar{f}_{ij}, 0\}$.
    \item  $(\pi_i, \pi_j) \notin \textup{int}(\Pi_{ji})$.
    \item Either $\pi_i \in\{\munderbar{\pi}_i, \bar{\pi}_i\}$  or $\pi_j \in\{\munderbar{\pi}_j, \bar{\pi}_j\}$.
\end{enumerate} 
\end{lemma}

Using Lemma \ref{lemma:extrPoints_1Neg}, we now present the following proposition on the characterization of the extreme points of the nonconvex set $\mathcal{L}_{ij}^-$.

\begin{proposition} \label{prop:extrPoints_compressorsNeg}
$\textup{extr}(\mathcal{L}_{ij}^-)$ is contained in  $\mathcal{U}_{ij}^{-1} \cup \mathcal{U}_{ij}^{-2} \cup \mathcal{U}_{ij}^{-3} \cup \mathcal{U}_{ij}^{-4}  \cup \mathcal{U}_{ij}^{-5} $, where 
\begin{align*}
    & \mathcal{U}_{ij}^{-1} :=  \{ ( {\pi}_i, \pi_j, 0, 0, 1) \in \mathbb{R}^5:  (\pi_i,\pi_j) \in \Pi_{ji} \}\\
    & \mathcal{U}_{ij}^{-2} :=  \{ ( \pi_i, \munderbar{\pi}_j,\munderbar{f}_{ij}, l_{ij}, 1) \in \mathbb{R}^5: l_{ij} = H_{ij} \munderbar{f}_{ij} ((\pi_{i}/\munderbar{\pi}_j)^{{\hat{\kappa}}} -1)  , \ \max\{\munderbar{\pi}_i ,  \munderbar{\pi}_j\munderbar{r}_{ij}  \} \le  \pi_i \le \min\{ \bar \pi_i,  \munderbar{\pi}_j\bar{r}_{ij}  \} \}\\
     & \mathcal{U}_{ij}^{-3} :=  \{ ( \pi_i, \bar{\pi}_j, \munderbar{f}_{ij}, l_{ij}, 1) \in \mathbb{R}^5: l_{ij} = H_{ij} \munderbar{f}_{ij} ((\pi_{i}/\bar{\pi}_j)^{{\hat{\kappa}}} -1)  , \ \max\{\munderbar{\pi}_i , \bar{\pi}_j\munderbar{r}_{ij}   \} \le  \pi_i \le \min\{ \bar \pi_i,  \bar{\pi}_j\bar{r}_{ij}  \} \}\\
    & \mathcal{U}_{ij}^{-4} :=  \{ ( \munderbar{\pi}_i, \pi_j, \munderbar{f}_{ij}, l_{ij}, 1) \in \mathbb{R}^5: l_{ij} = H_{ij}\munderbar{f}_{ij} ((\munderbar{\pi}_i/\pi_j)^{{\hat{\kappa}}} -1) , \ 
    \max\{\munderbar{\pi}_j , \munderbar{\pi}_i / \bar{r}_{ij}  \} \le  \pi_j \le \min\{ \bar \pi_j,  \munderbar{\pi}_i / \munderbar{r}_{ij}  \} \} \\
    & \mathcal{U}_{ij}^{-5} :=  \{ ( \bar{\pi}_i, \pi_j, \munderbar{f}_{ij}, l_{ij}, 1) \in \mathbb{R}^5: l_{ij} = H_{ij}\munderbar{f}_{ij} ((\bar{\pi}_i/\pi_j)^{{\hat{\kappa}}} -1) , \ 
    \max\{\munderbar{\pi}_j , \bar{\pi}_i / \bar{r}_{ij}  \} \le  \pi_j \le \min\{ \bar \pi_j,  \bar{\pi}_i / \munderbar{r}_{ij}  \} \}.
\end{align*}
\end{proposition}

Now, we are now ready to present the convex hull of the nonconvex sets $\mathcal{U}_{ij}^{-2}, \dots, \mathcal{U}_{ij}^{-5}$. In fact, the following proposition is an immediate consequence of Proposition~\ref{prop:extrPoints_compressorsNeg}, Lemma~\ref{lem:powerConeReprPlus} and Lemma~\ref{lem:powerConeReprMinus}:
\begin{proposition} \label{proposition:convexHull_compressors_U+_disjNeg}
We have 
\begin{align*}
    & \textup{conv}(\mathcal{U}_{ij}^{-2}) =  \left\{ ( \pi_i, \munderbar{\pi}_j,\munderbar{f}_{ij}, l_{ij}, 1) \in \mathbb{R}^5:
    (\pi_i,-l_{ij}) \in \textup{conv} \big( \mathcal{V}^{\hat \kappa}( - H_{ij} \munderbar f_{ij} , \munderbar{\pi}_j, 
  \max\{\munderbar{\pi}_i ,  \munderbar{\pi}_j\munderbar{r}_{ij}  \} , \min\{ \bar \pi_i,  \munderbar{\pi}_j\bar{r}_{ij}   \} ) \big) \right\}\\
     & \textup{conv}(\mathcal{U}_{ij}^{-3}) =   \left\{ ( \pi_i, \bar{\pi}_j, \munderbar{f}_{ij}, l_{ij}, 1) \in \mathbb{R}^5: 
    (\pi_i,-l_{ij}) \in \textup{conv} \big( \mathcal{V}^{\hat \kappa} ( 
     - H_{ij} \munderbar{f}_{ij} , \bar{\pi}_j ,\max\{\munderbar{\pi}_i , \bar{\pi}_j\munderbar{r}_{ij}   \} , \min\{ \bar \pi_i,  \bar{\pi}_j\bar{r}_{ij}   \} ) \big) \right \}\\
    & \textup{conv}(\mathcal{U}_{ij}^{-4}) =   \left\{ ( \munderbar{\pi}_i, \pi_j, \munderbar{f}_{ij}, l_{ij}, 1) \in \mathbb{R}^5:
    (\pi_j,-l_{ij}) \in \textup{conv} \big( \mathcal{W}^{\hat \kappa} (  -H_{ij}\munderbar{f}_{ij} , \munderbar{\pi}_i  ,  
    \max\{\munderbar{\pi}_j , \munderbar{\pi}_i / \bar{r}_{ij}  \}, \min\{ \bar \pi_j,  \munderbar{\pi}_i / \munderbar{r}_{ij}  \} ) \big) \right\} \\
    & \textup{conv}(\mathcal{U}_{ij}^{-5}) =   \left\{ ( \bar{\pi}_i, \pi_j, \munderbar{f}_{ij}, l_{ij}, 1) \in \mathbb{R}^5:
    (\pi_j,-l_{ij}) \in \textup{conv} \big( \mathcal{W}^{\hat \kappa} ( - H_{ij}\munderbar{f}_{ij} , \bar{\pi}_i  ,  
    \max\{\munderbar{\pi}_j , \bar{\pi}_i / \bar{r}_{ij}  \} , \min\{ \bar \pi_j,  \bar{\pi}_i / \munderbar{r}_{ij}  \}) \big) \right\}.
\end{align*}
\end{proposition} 
Finally, we formalize our main results on the convex hull of the nonconvex set $\mathcal{L}_{ij}^-$ and the polytope $\mathcal{L}_{ij}^0$. 
\begin{proposition} 
\label{proposition:convexHull_compressors_neg_zero}
 $\textup{conv}(\mathcal{L}_{ij}^-)$ 
and  $\textup{conv}(\mathcal{L}_{ij}^0)$ are POWr.
\end{proposition} 

We are now ready to prove Theorem~\ref{theorem:convexHull_compressors_union}.
\begin{proof}[Proof of Theorem~\ref{theorem:convexHull_compressors_union}]
 Due to Propositions~\ref{proposition:convexHull_compressors} and~\ref{proposition:convexHull_compressors_neg_zero}, we deduce that $\textup{conv}(\mathcal{L}_{ij}^+)$, $\textup{conv}(\mathcal{L}_{ij}^-)$ and $\textup{conv}(\mathcal{L}_{ij}^0)$  are POWr. Since these sets are also compact and  satisfy the relation $ \mathcal{L}_{ij} = \mathcal{L}^+_{ij} \cup \mathcal{L}^-_{ij} \cup \mathcal{L}^0_{ij} $, the convex hull of their union is also POWr \cite{ben2001lectures}.
\end{proof}
We note that although conic interior solvers allow for non-symmetric cones, e.g., power cones,  
their algorithmic frameworks for MINLPs with non-symmetric cones are not as efficient as those with symmetric cones, e.g., second-order cones. In fact, MISOCPs can be solved quite efficiently at the state-of-the-art, see, e.g., Mittelmann benchmarks \cite{Mittelmann}. To fully utilize the current solvers, we propose an SOCr outer-approximation of the compression set below. 

\subsubsection{Second-order cone representable outer-approximation}  \label{ss:outer_approx_compressors}

In practice, the isentropic exponent $\kappa$ can be approximated by the constant 0.228 \citep{Koch2015}. We also follow this common practice and set the exponent $\hat{\kappa} = \kappa / 2$ to 0.114. By doing so, we are also able to provide an SOCr outer-approximation to the power cone representable convex hull of $\mathcal{L}_{ij}$ for practicality. Moreover, our formulation can be readily given to the state-of-the-art solvers. Before we explain our construction, let us first present our third main result, which provides the mathematical formulation for the SOCr outer-approximation to the convex hull of the compression set in the extended space:

\begin{theorem} \label{theorem:compressorOuterApprox}
Let $n = 0.114$ and $a > 0$. Then, the following set provides an extended formulation for the SOCr outer-approximation of $\textup{conv}(\mathcal{L}_{ij})$:
\begin{align*} \mathcal{Y}_{ij} :=
 \left\{
 \begin{array}{cc}  
    (\pi_i, \pi_j, f_{ij}, l_{ij}, x_{ij}, x_{ij}^+, x_{ij}^-) \in \mathbb{R}^7\\[0.15cm]
    (\pi_i^{-5}, \pi_j^{-5}, f_{ij}^{-5}, l_{ij}^{-5}, \lambda_{ij}^{-5}) \in \mathbb{R}^5 \\[0.15cm]
    \vdots \\[0.15cm] 
    (\pi_i^{+5}, \pi_j^{+5}, f_{ij}^{+5}, l_{ij}^{+5}, \lambda_{ij}^{+5}) \in \mathbb{R}^5 
\end{array}
: 
\begin{array}{cc}
    & \sum_{k=1}^{11} \lambda_{ij}^k = 1, \ \lambda_{ij}^k \ge 0, \ k = 1, \dots, 11\\[0.15cm]
     &x_{ij}^+ = \sum_{k=+1}^{+5}  \lambda_{ij}^{k}, \ x_{ij}^- = \sum_{k=-5}^{-1} \lambda_{ij}^{k} \\[0.15cm]
     & \pi_i = \sum_{k=-5}^{+5} \pi_i^k, \ \pi_j = \sum_{k=-5}^{+5} \pi_j^k\\[0.15cm]
    & f_{ij} = \sum_{k=-5}^{+5}  f_{ij}^k, \ l_{ij} = \sum_{k=-5}^{+5} l_{ij}^k \\[0.15cm]
    & \munderbar{\pi}_i \le \pi_i \le \bar{\pi}_i, \ \munderbar{\pi}_j \le \pi_j \le \bar{\pi}_j, \ \munderbar{f}_{ij} x_{ij}^- \le f_{ij} \le \bar{f}_{ij} x_{ij}^+\\[0.15cm]
    &  x_{ij} = x_{ij}^+ + x_{ij}^-, \ x_{ij}^+ \in [0,1], \ x_{ij}^- \in [0,1] \\[0.15cm] 
\end{array}
\right. \\ & \hspace{-13cm} \left.
\begin{array}{cc}  
    &\munderbar{\pi}_i \lambda_{ij}^{0} \le \pi_i^{0} \le \bar{\pi}_i \lambda_{ij}^{0}, \ \munderbar{\pi}_j \lambda_{ij}^{0} \le \pi_j^{0} \le \bar{\pi}_j \lambda_{ij}^{0},  \ f_{ij}^{0} = l_{ij}^{0} = 0  \\[0.15cm]
    &\munderbar{\pi}_i \lambda_{ij}^{+1} \le \pi_i^{+1} \le \bar{\pi}_i \lambda_{ij}^{+1}, \ \munderbar{\pi}_j \lambda_{ij}^{+1} \le \pi_j^{+1} \le \bar{\pi}_j \lambda_{ij}^{+1}, \ \munderbar{r}_{ij} \pi_i^{+1} \le \pi_j^{+1} \le \bar{r}_{ij} \pi_i^{+1}, \ f_{ij}^{+1} = l_{ij}^{+1} = 0  \\[0.15cm]
    &(\pi_j^{+2},  l_{ij}^{+2}, \lambda_{ij}^{+2}) \in \mathcal{\check V} (  H_{ij} \bar f_{ij} , \munderbar{\pi}_i, 
    \max\{\munderbar{\pi}_j ,  \munderbar{\pi}_i\munderbar{r}_{ij}  \} , \min\{ \bar \pi_j,  \munderbar{\pi}_i\bar{r}_{ij}  \} ), \ \pi_i^{+2} = \munderbar{\pi}_i \lambda_{ij}^{+2}, \ f_{ij}^{+2} = \bar{f}_{ij} \lambda_{ij}^{+2} \\[0.15cm]
    &(\pi_j^{+3}, l_{ij}^{+3}, \lambda_{ij}^{+3}) \in \mathcal{\check V}  (H_{ij} \bar{f}_{ij} , \bar{\pi}_i , \max\{\munderbar{\pi}_j , \bar{\pi}_i\munderbar{r}_{ij}   \} ,    \min\{ \bar \pi_j,  \bar{\pi}_i\bar{r}_{ij} \} ), \ \pi_i^{+3} = \bar{\pi}_i \lambda_{ij}^{+3}, \ f_{ij}^{+3} = \bar{f}_{ij} \lambda_{ij}^{+3} \\[0.15cm]
    &(\pi_i^{+4}, l_{ij}^{+4}, \lambda_{ij}^{+4}) \in \mathcal{\check W}  (  H_{ij}\bar{f}_{ij} , \munderbar{\pi}_j  ,  
    \max\{\munderbar{\pi}_i , \munderbar{\pi}_j / \bar{r}_{ij}  \} , \min\{ \bar \pi_i,  \munderbar{\pi}_j / \munderbar{r}_{ij}  \} ) , \ \pi_j^{+4} = \munderbar{\pi}_j \lambda_{ij}^{+4}, \ f_{ij}^{+4} = \bar{f}_{ij} \lambda_{ij}^{+4} \\[0.15cm]
    &(\pi_i^{+5}, l_{ij}^{+5}, \lambda_{ij}^{+5}) \in \mathcal{\check W}   (  H_{ij}\bar{f}_{ij} , \bar{\pi}_j  ,  
    \max\{\munderbar{\pi}_i , \bar{\pi}_j / \bar{r}_{ij}  \} , \min\{ \bar \pi_i,  \bar{\pi}_j / \munderbar{r}_{ij}  \} ) , \ \pi_j^{+5} = \bar{\pi}_j \lambda_{ij}^{+5}, \ f_{ij}^{+5} = \bar{f}_{ij} \lambda_{ij}^{+5} \\[0.15cm]
    &\munderbar{\pi}_i \lambda_{ij}^{-1} \le \pi_i^{-1} \le \bar{\pi}_i \lambda_{ij}^{-1}, \ \munderbar{\pi}_j \lambda_{ij}^{-1} \le \pi_j^{-1} \le \bar{\pi}_j \lambda_{ij}^{-1}, \ \munderbar{r}_{ij} \pi_j^{-1} \le \pi_i^{-1} \le \bar{r}_{ij} \pi_j^{-1}, \ f_{ij}^{-1} = l_{ij}^{-1} = 0  \\[0.15cm]
    &(\pi_i^{-2},  -l_{ij}^{-2}, \lambda_{ij}^{-2}) \in \mathcal{\check V} ( - H_{ij} \munderbar f_{ij} , \munderbar{\pi}_j, 
    \max\{\munderbar{\pi}_i ,  \munderbar{\pi}_j\munderbar{r}_{ij}  \} , \min\{ \bar \pi_i,  \munderbar{\pi}_j\bar{r}_{ij}   \} ) , \ \pi_j^{-2} = \munderbar{\pi}_j \lambda_{ij}^{-2}, \ f_{ij}^{-2} = \munderbar{f}_{ij} \lambda_{ij}^{-2} \\[0.15cm]
    &(\pi_i^{-3}, -l_{ij}^{-3}, \lambda_{ij}^{-3}) \in \mathcal{\check V} ( 
     - H_{ij} \munderbar{f}_{ij} , \bar{\pi}_j ,\max\{\munderbar{\pi}_i , \bar{\pi}_j\munderbar{r}_{ij}   \} , \min\{ \bar \pi_i,  \bar{\pi}_j\bar{r}_{ij}   \} )  , \ \pi_j^{-3} = \bar{\pi}_j \lambda_{ij}^{-3}, \ f_{ij}^{-3} = \munderbar{f}_{ij} \lambda_{ij}^{-3} \\[0.15cm]
    &(\pi_j^{-4}, -l_{ij}^{-4}, \lambda_{ij}^{-4}) \in \mathcal{\check W}  (  -H_{ij}\munderbar{f}_{ij} , \munderbar{\pi}_i  ,   \max\{\munderbar{\pi}_j , \munderbar{\pi}_i / \bar{r}_{ij}  \}, \min\{ \bar \pi_j,  \munderbar{\pi}_i / \munderbar{r}_{ij}  \} ), \ \pi_i^{-4} = \munderbar{\pi}_i \lambda_{ij}^{-4}, \ f_{ij}^{-4} = \munderbar{f}_{ij} \lambda_{ij}^{-4} \\[0.15cm]
    &(\pi_j^{-5}, -l_{ij}^{-5}, \lambda_{ij}^{-5}) \in \mathcal{\check W}  ( - H_{ij}\munderbar{f}_{ij} , \bar{\pi}_i  ,  
    \max\{\munderbar{\pi}_j , \bar{\pi}_i / \bar{r}_{ij}  \} , \min\{ \bar \pi_j,  \bar{\pi}_i / \munderbar{r}_{ij}  \}) , \ \pi_i^{-5} = \bar{\pi}_i \lambda_{ij}^{-5}, \ f_{ij}^{-5} = \munderbar{f}_{ij} \lambda_{ij}^{-5} \\[0.15cm]
 \end{array}   
    \right\},
\end{align*} 
where
\begin{align*}
\mathcal{\check V}(a,b,\munderbar{x}, \bar x) = \{(x,y, \lambda)\in \mathbb{R}^3_+: \exists (u, v_1, v_2) \in \mathbb{R}^3_+:\ &  u^2 \le \lambda v_1, \  v_1^2 \le \lambda v_2, \ b v_2^2 \le  \lambda x \\ & y = a(u- \lambda), \ \munderbar{x} \lambda \le x \le \bar{x} \lambda \\
        &
        y \ge a\Big ( \frac{ (\bar{x}/b) ^n - (\munderbar{x}/b)^n}{\bar{x} - \munderbar{x}} \Big ) (x - \munderbar{x}\lambda) + a ( (\munderbar{x}/b)^n-1)\lambda \},
   \end{align*}
with $0<b\le\munderbar{x} \le \bar x$ and
    \begin{align*}
       \mathcal{\check W}(a,b,\munderbar{x}, \bar x) = \{(x,y, \lambda)\in \mathbb{R}^3_+: \exists (u, v_1, v_2, v_3) \in \mathbb{R}^4_+:\ & u v_1 \ge \lambda^2, \ \lambda v_2 \ge v_1^2, \ \lambda v_3 \ge v_2^2 , \ u x \ge b v_3^2
       \\ &y = a(u-\lambda), \ \munderbar{x}\lambda \le x \le \bar{x}\lambda\\
        &
        y \le a\Big ( \frac{ (b/\bar{x}) ^n - (b/\munderbar{x})^n}{\bar{x} - \munderbar{x}} \Big ) (x - \munderbar{x}\lambda) + a ( (b/\munderbar{x})^n-1) \lambda\},
    \end{align*}
    with $0\le\munderbar{x} \le \bar x<b$.

\end{theorem}

We note that the proof of Theorem~\ref{theorem:compressorOuterApprox} is constructive. First, we obtain an SOCr outer-approximation to each POWr convex hull of the sets $\mathcal{L}_{ij}^+ $, $\mathcal{L}_{ij}^- $ and $\mathcal{L}_{ij}^0$. The polytope $\mathcal{L}_{ij}^0$ is trivially SOCr. Hence, we will only focus on $\mathcal{L}_{ij}^+ $ and $\mathcal{L}_{ij}^- $. Here, we would like to remind the reader of Propositions~\ref{proposition:convexHull_compressors_U+_disj} and \ref{proposition:convexHull_compressors_U+_disjNeg}. Our task is exactly to find outer-approximations to the nontrivial POWr sets, in particular, finding outer-approximations to
$\mathcal{V}^{\hat{\kappa}}(a,b,\munderbar{x}, \bar x)$ and $\mathcal{W}^{\hat{\kappa}}(a,b,\munderbar{x}, \bar x)$ in Propositions~\ref{proposition:convexHull_compressors_U+_disj} and \ref{proposition:convexHull_compressors_U+_disjNeg}.
\begin{figure}[H]
\centering
\begin{subfigure}{0.45\textwidth}
    \centering
    \scalebox{1}{
    \begin{tikzpicture}
        \def\xmin{0.5}
        \def\xmax{4}

        \draw[gray, opacity=0.5, dashdotted, xstep=0.5cm, ystep=0.5cm] (0,0) grid (4.25, 2.75);
        \draw[->] (0,0) -- (4.25,0) node[right] {$x$};
        \draw[->] (0,0) -- (0,2.75) node[above] {$y$};
        
        \draw[domain=\xmin:\xmax, smooth, variable=\x, thick, blue]     plot ({\x}, { -2.4 + 3.3*pow( (\x ), 0.2)});

        \draw[domain=\xmin:\xmax, smooth, variable=\x, densely dotted, red]     plot ({\x}, { -2.2 + 3.3*pow( (1.01*\x ), 0.25)});

        \draw[red, densely dotted]
        (0.5,  -2.2 + 3.3*1.01^0.25*0.5^0.25) --
        (0.5, -2.4 + 3.3*0.5^0.2) -- 
        (4, -2.4 + 3.3*4^0.2) -- (4,  -2.2 + 3.3*1.01^0.25*4^0.25);

        \fill[blue, opacity=0.1]
            plot[domain=\xmin:\xmax] 
                ({\x}, {0.4233*\x + 0.2612}) -- 
            plot[domain=\xmax:\xmin] 
                ({\x}, { -2.2 + 3.3*pow( (1.01*\x ), 0.25)})
            -- cycle;

        \filldraw[blue] (\xmin, {({ -2.4 + 3.3*pow( (\xmin ), 0.2)})}) circle (1pt) node[below left] {$\munderbar{x}$};
        \filldraw[blue] (\xmax, {({ -2.4 + 3.3*pow( (\xmax ), 0.2)})}) circle (1pt) node[below right] {$\bar{x}$};

    \end{tikzpicture}
    } 
\caption{The outer-approximation to $\mathcal{V}^{\hat{\kappa}}$.}
 \label{fig:outerApproxV}
     \end{subfigure}
\begin{subfigure}{0.45\textwidth}
    \centering
    \scalebox{1}{
    \begin{tikzpicture}
        \def\a{1.1}
        \def\b{2}
        \def\n{0.114}
        \def\xmin{0.5}
        \def\xmax{4}

        \draw[gray, opacity=0.5, dashdotted, xstep=0.5cm, ystep=0.5cm] (0,0) grid (4.25, 2.75);
        \draw[->] (0,0) -- (4.25,0) node[right] {$x$};
        \draw[->] (0,0) -- (0,2.75) node[above] {$y$};
        
        \draw[domain=\xmin:\xmax, smooth, variable=\x, thick, blue]     plot ({\x}, { -1.4 + 3.3*pow( (\x ), -0.22)});

        \draw[domain=\xmin:\xmax, smooth, variable=\x, densely dotted, red]     plot ({\x}, { -1.7 + 3.3*pow( (1.01*\x ), -0.28)});

        \draw[red, densely dotted]
        (0.5,   2.296) --
        (0.5, -1.4 + 3.3*0.5^-0.22) --
        (4, -1.4 + 3.3*4^-0.22)--
        (4,  0.532);
    
        \fill[blue, opacity=0.1]
           plot[domain=\xmin:\xmax] 
               ({\x}, {-0.4032*\x + 2.6452}) -- 
           plot[domain=\xmax:\xmin] 
               ({\x}, { -1.7 + 3.3*pow( (1.01*\x ), -0.28)}) -- 
            cycle;

        \filldraw[blue] (\xmin, {({  -1.4 + 3.3*pow( (\xmin ), -0.22)})}) circle (1pt) node[below left] {$\munderbar{x}$};
        \filldraw[blue] (\xmax, {({  -1.4 + 3.3*pow( (\xmax), -0.22)})}) circle (1pt) node[below right] {$\bar{x}$};

    \end{tikzpicture}
    } 
\caption{The outer-approximation to $\mathcal{W}^{\hat{\kappa}}$.}
 \label{fig:outerApproxW}
     \end{subfigure}
\caption{The illustrations of the second-order cone representable outer-approximations.}
 \label{fig:outerApproxPOWr}
\end{figure}
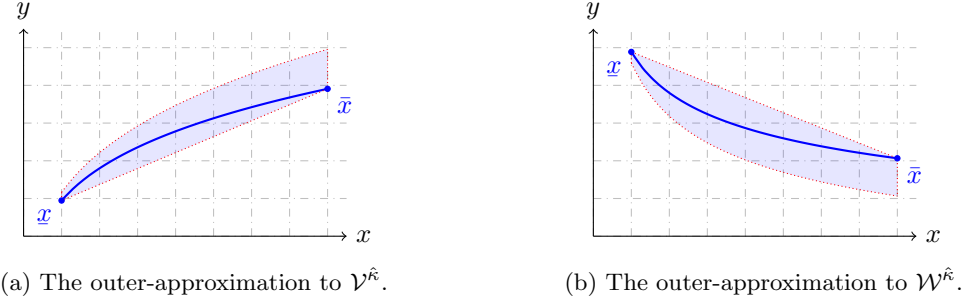
Next, we focus on a very special case in which we first approximate $\hat{\kappa}$ by the constants 1/8 and -1/9 for the POWr sets $\mathcal{V}^{\hat{\kappa}}$ and $\mathcal{W}^{\hat{\kappa}}$, respectively. In Figure~\ref{fig:outerApproxPOWr}, we geometrically summarize our construction: We overestimate the concave function $y = x^{0.114}$ with $y = x^{1/8}$ and find the lower envelope for $\mathcal{V}^{\hat{\kappa}}$ (Figure~\ref{fig:outerApproxV}). Similarly, we underestimate the convex function $y = x^{-0.114}$ with $y = x^{-1/9}$ and find the upper envelope for $\mathcal{W}^{\hat{\kappa}}$ (Figure~\ref{fig:outerApproxW}). Also, we observe that the hypograph of $y = x^{1/8}$ and the epigraph of $y = x^{-1/9}$ are both SOCr. These observations allow us to obtain SOCr outer-approximations for these POWr sets. Let us formalize them by presenting the following lemmas. 
\begin{lemma} \label{lem:secondOrderConeReprPlus}
    Consider the following convex set 
        \begin{align*}
     \mathcal{\hat V}(a,b,\munderbar{x}, \bar x) = \{(x,y)\in \mathbb{R}^2_+: \exists u \in \mathbb{R}_+:\ & u \le x^{1/8}, \  y = a(u-1), \ \munderbar{x} \le x \le \bar{x} \\
        &
        y \ge a\Big ( \frac{ (\bar{x}/b) ^n - (\munderbar{x}/b)^n}{\bar{x} - \munderbar{x}} \Big ) (x - \munderbar{x}) + a ( (\munderbar{x}/b)^n-1) \},
    \end{align*}
    where $n = 0.114$, $a >0$ and  $0<b\le\munderbar{x} \le \bar x$. Then, $\mathcal{\hat V}(a,b,\munderbar{x}, \bar x)$ is SOCr and it provides an outer-approximation to $ \textup{conv}(\mathcal{V}^n(a,b,\munderbar{x}, \bar x))$.
\end{lemma}
\begin{proof}
Suppose $n = 0.114$. Then, the proof immediately follows from the following relation 
    \begin{align*}
\{ (u,x)\in\mathbb{R}^2_+ :  u \le x^n\}   
\subseteq
&\{ (u,x)\in\mathbb{R}^2_+ :  u \le x^{1/8}\} \\
=
&\{ (u,x)\in\mathbb{R}^2_+ : \exists v\in\mathbb{R}^2_+ :  u^2 \le v_1, v_1^2 \le v_2, v_2^2 \le x \}.
\end{align*}
\end{proof}
\begin{lemma} \label{lem:secondOrderConeReprMinus}
    Consider the following convex set
    \begin{align*}
       \mathcal{\hat W}(a,b,\munderbar{x}, \bar x) = \{(x,y)\in \mathbb{R}^2_+: \exists u \in \mathbb{R}_+:\ &u \ge x^{-1/9}, \ y = a(u-1), \ \munderbar{x} \le x \le \bar{x}\\
        &
        y \le a\Big ( \frac{ (b/\bar{x}) ^n - (b/\munderbar{x})^n}{\bar{x} - \munderbar{x}} \Big ) (x - \munderbar{x}) + a ( (b/\munderbar{x})^n-1) \},
    \end{align*}
   where $n = 0.114 $, $a >0$ and $0\le\munderbar{x} \le \bar x<b$. Then, $\mathcal{\hat W}(a,b,\munderbar{x}, \bar x)$ is SOCr and it provides an outer-approximation to $\textup{conv}(\mathcal{W}^n(a,b,\munderbar{x}, \bar x))$.
\end{lemma}
\begin{proof}
    Suppose $n = 0.114$. Then, the proof immediately follows from the following relation 
    \begin{align*}
        \{ (u,x)\in\mathbb{R}^2_{++} :  u \ge x^{-n} \}   
        \subseteq
        &\{ (u,x)\in\mathbb{R}^2_{++} :  u \ge x^{-1/9}\} \\
        =
        &\{ (u,x)\in\mathbb{R}^2_{++} : \exists v\in\mathbb{R}^3_+ :  u v_1 \ge 1, v_2 \ge v_1^2, v_3 \ge v_2^2 , u x \ge v_3^2 \}.
    \end{align*}
\end{proof}
The following proposition is an immediate consequence of Proposition~\ref{proposition:convexHull_compressors_U+_disj}, Lemma~\ref{lem:secondOrderConeReprPlus} and Lemma~\ref{lem:secondOrderConeReprMinus}:
\begin{proposition} \label{prop:outerApprox_compressors_L+_disj}
$\textup{conv}(\mathcal{L}_{ij}^+)$ is contained in  $\mathcal{\hat L}_{ij}^+:=\mathcal{U}_{ij}^{+1} \cup \mathcal{\hat U}_{ij}^{+2} \cup \mathcal{\hat U}_{ij}^{+3} \cup \mathcal{\hat U}_{ij}^{+4}  \cup \mathcal{\hat U}_{ij}^{+5}$, where
\begin{align*}
    & \mathcal{\hat U}_{ij}^{+2} = \{ (\munderbar{\pi}_i, \pi_j, \bar{f}_{ij}, l_{ij}, 1) \in \mathbb{R}^5: 
    (\pi_i,l_{ij}) \in \mathcal{\hat V}(  H_{ij} \bar f_{ij} , \munderbar{\pi}_i, 
    \max\{\munderbar{\pi}_j ,  \munderbar{\pi}_i\munderbar{r}_{ij}  \} , \min\{ \bar \pi_j,  \munderbar{\pi}_i\bar{r}_{ij}  \} ) \}\\
     & \mathcal{\hat U}_{ij}^{+3} =   \{ (\bar{\pi}_i, \pi_j, \bar{f}_{ij}, l_{ij}, 1) \in \mathbb{R}^5: 
    (\pi_i,l_{ij}) \in  \mathcal{\hat V} ( 
     H_{ij} \bar{f}_{ij} , \bar{\pi}_i , \max\{\munderbar{\pi}_j , \bar{\pi}_i\munderbar{r}_{ij}   \} ,    \min\{ \bar \pi_j,  \bar{\pi}_i\bar{r}_{ij} \} ) \} \\
    & \mathcal{\hat U}_{ij}^{+4} =  \{ (\pi_i, \munderbar{\pi}_j, \bar{f}_{ij}, l_{ij}, 1) \in \mathbb{R}^5: 
    (\pi_j,l_{ij}) \in \mathcal{\hat W} (  H_{ij}\bar{f}_{ij} , \munderbar{\pi}_j  ,  
    \max\{\munderbar{\pi}_i , \munderbar{\pi}_j / \bar{r}_{ij}  \} , \min\{ \bar \pi_i,  \munderbar{\pi}_j / \munderbar{r}_{ij}  \} ) \} \\
    & \mathcal{\hat U}_{ij}^{+5} =  \{ (\pi_i, \bar{\pi}_j, \bar{f}_{ij}, l_{ij}, 1) \in \mathbb{R}^5:
    (\pi_j,l_{ij}) \in \mathcal{\hat W} (  H_{ij}\bar{f}_{ij} , \bar{\pi}_j  ,  
    \max\{\munderbar{\pi}_i , \bar{\pi}_j / \bar{r}_{ij}  \} , \min\{ \bar \pi_i,  \bar{\pi}_j / \munderbar{r}_{ij}  \} ) \}.
\end{align*}
Moreover, $\textup{conv}(\mathcal{\hat L}_{ij}^+)$ is SOCr.
\end{proposition}
Similarly, the following proposition is an immediate consequence of Proposition~\ref{proposition:convexHull_compressors_U+_disjNeg}, Lemma~\ref{lem:secondOrderConeReprPlus} and Lemma~\ref{lem:secondOrderConeReprMinus}:
\begin{proposition} \label{prop:outerApprox_compressors_L-_disj}
$\textup{conv}(\mathcal{L}_{ij}^-)$ is contained in  $\mathcal{\hat L}_{ij}^-:=\mathcal{U}_{ij}^{-1} \cup \mathcal{\hat U}_{ij}^{-2} \cup \mathcal{\hat U}_{ij}^{-3} \cup \mathcal{\hat U}_{ij}^{-4}  \cup \mathcal{\hat U}_{ij}^{-5}$, where
\begin{align*}
    & \mathcal{\hat U}_{ij}^{-2} =  \{ ( \pi_i, \munderbar{\pi}_j,\munderbar{f}_{ij}, l_{ij}, 1) \in \mathbb{R}^5:
    (\pi_i,-l_{ij}) \in \mathcal{\hat V}( - H_{ij} \munderbar f_{ij} , \munderbar{\pi}_j, 
  \max\{\munderbar{\pi}_i ,  \munderbar{\pi}_j\munderbar{r}_{ij}  \} , \min\{ \bar \pi_i,  \munderbar{\pi}_j\bar{r}_{ij}   \} ) \}\\
     & \mathcal{\hat U}_{ij}^{-3} =   \{ ( \pi_i, \bar{\pi}_j, \munderbar{f}_{ij}, l_{ij}, 1) \in \mathbb{R}^5: 
    (\pi_i,-l_{ij}) \in  \mathcal{ \hat V} ( 
     - H_{ij} \munderbar{f}_{ij} , \bar{\pi}_j ,\max\{\munderbar{\pi}_i , \bar{\pi}_j\munderbar{r}_{ij}   \} , \min\{ \bar \pi_i,  \bar{\pi}_j\bar{r}_{ij}   \} )   \}\\
    & \mathcal{\hat U}_{ij}^{-4} =   \{ ( \munderbar{\pi}_i, \pi_j, \munderbar{f}_{ij}, l_{ij}, 1) \in \mathbb{R}^5:
    (\pi_j,-l_{ij}) \in \mathcal{\hat W} (  -H_{ij}\munderbar{f}_{ij} , \munderbar{\pi}_i  ,  
    \max\{\munderbar{\pi}_j , \munderbar{\pi}_i / \bar{r}_{ij}  \}, \min\{ \bar \pi_j,  \munderbar{\pi}_i / \munderbar{r}_{ij}  \} ) \} \\
    & \mathcal{\hat U}_{ij}^{-5} =   \{ ( \bar{\pi}_i, \pi_j, \munderbar{f}_{ij}, l_{ij}, 1) \in \mathbb{R}^5:
    (\pi_j,-l_{ij}) \in \mathcal{\hat W} ( - H_{ij}\munderbar{f}_{ij} , \bar{\pi}_i  ,  
    \max\{\munderbar{\pi}_j , \bar{\pi}_i / \bar{r}_{ij}  \} , \min\{ \bar \pi_j,  \bar{\pi}_i / \munderbar{r}_{ij}  \})  \}.
\end{align*}
Moreover, $\textup{conv}(\mathcal{\hat L}_{ij}^-)$ is SOCr.
\end{proposition} 
Now, it is relatively easy to give a formal description of the SOCr outer-approximation of $\textup{conv}(\mathcal{\hat L}_{ij})$. In particular, it suffices to take the convex combination of sets $\mathcal{\hat L}_{ij}^+, \mathcal{\hat L}_{ij}^-$ and $\mathcal{L}_{ij}^0$. The extended formulation for this convex outer-approximation is given by Theorem~\ref{theorem:compressorOuterApprox}.

\section{Solution Methodology} \label{section:solutionAlgorithm}
In Section \ref{section:mainResults}, we presented our main results on the nonconvex feasible sets for pipes and compressors. In this section, we will describe how we can solve the multi-period gas transportation and storage problem \eqref{eq:MINLP} by using them. We summarize our two-step solution framework in Algorithm \ref{alg:MINLP_GlobalAlgorithm}. The input of our algorithm is the choice of the mixed-integer convex relaxation for the nonconvex MINLP problem \eqref{eq:MINLP}. The tolerance of the convex relaxation is denoted by $\epsilon$. We also denote the $\epsilon$-optimal values of the binary and continuous variables produced by this relaxation as $x^*$ and $y^*$, respectively. In Step \ref{MINLP_GlobalAlgorithm:solveRelax} of Algorithm \ref{alg:MINLP_GlobalAlgorithm}, a mixed-integer convex relaxation of our choice is solved to optimality. The dual bound for this relaxation provides a lower bound on the objective value of problem \eqref{eq:MINLP}. In order to obtain (locally) optimal solutions, all the optimal binary values $x^*$ produced by the relaxation are fixed, which reduces the nonconvex MINLP problem \eqref{eq:MINLP} to a nonlinear program (NLP). By passing the optimal continuous values $y^*$ produced by the relaxation as the initial point, we solve this program with a local interior point solver in Step \ref{MINLP_GlobalAlgorithm:solveNLP}. If the solver converges to a (locally) optimal solution, which is also feasible for problem \eqref{eq:MINLP}, we obtain an upper bound on its optimal objective value. By using these lower and upper bounds, the relative optimality gap percentage is then calculated.
%
%
%
%
%
\begin{algorithm}
\caption{MINLP}
\label{alg:MINLP_GlobalAlgorithm}
\begin{algorithmic}[1]
\REQUIRE \texttt{convexRelaxation}.
\ENSURE Primal and dual bounds for problem \eqref{eq:MINLP}, the relative optimality gap in percentages.
\STATE Solve \texttt{convexRelaxation} of problem \eqref{eq:MINLP}. \label{MINLP_GlobalAlgorithm:solveRelax}
\STATE Obtain $\epsilon$-optimal solutions $x^*$ and $y^*$ produced by \texttt{convexRelaxation}.
\STATE Assign the dual bound of \texttt{convexRelaxation} to $LB$.
\STATE Solve problem \eqref{eq:MINLP} by fixing $x^*$ and passing $y^*$ as an initial point. \label{MINLP_GlobalAlgorithm:solveNLP}
\STATE Assign the primal bound of problem \eqref{eq:MINLP} to $UB$.
\STATE Assign the relative optimality $(1 - LB/UB)\%$ to $Gap$.
 \RETURN $LB, UB, Gap$.
\end{algorithmic}
\end{algorithm}

We conclude this section by providing the mixed-integer convex relaxations that can be used within Algorithm \ref{alg:MINLP_GlobalAlgorithm} to solve the multi-period gas transportation and storage problem \eqref{eq:MINLP}. These relaxations differ based on the choice of the convexification approaches for the nonconvex Weymouth set $\mathcal{X}$, which are explained in Section~\ref{ss:feasible_sets_pipes}.
Note that we use the mixed-integer version of the SOCr set described as in Theorem~\ref{theorem:compressorOuterApprox} for the nonconvex feasible set for compressors. We consider the following cases for \texttt{convexRelaxation} in Algorithm~\ref{alg:MINLP_GlobalAlgorithm} with respect to the convexification of the Weymouth set:
\begin{itemize}    
    \item \texttt{R1}: The MISOCP relaxation with the convex hull  as described in Section~\ref{ss:socr_convex_hull_pipes}.
    \item \texttt{R2}: The MISOCP relaxation with the SOCr outer-approximation proposed by \cite{Borraz2016SOCP}  as described in Section~\ref{ss:outer_approx_pipes}.
    \item \texttt{R3}: The MISOCP relaxation with the polyhedral outer-approximation adapted from \cite{Pfetsch2015}
    as described in Section~\ref{ss:outer_approx_pipes}.
\end{itemize}

\section{Numerical Experiments} \label{section:computationalStudy}
In this section, we present our computational study conducted on various GasLib test instances from the literature \citep{Gaslib, Hennings2023}. We focus on short-term gas network operations and let $|\mathcal{T}| = 24$. Note that demand parameters in GasLib instances are given for a single period. Given this, we generate a new demand dataset by using the 24-hour gas demand data in \citet{Schwele2019SOCP} such that the single-period demand from GasLib instances is the average gas demand within the planning horizon. Table \ref{table:gaslibInstances} shows the cardinality of gas network elements for all the test instances used in the computational study.

\begin{table}[h!]
\begin{center}
{\scalebox{1.00}{
\begin{tabular}{cccccccc}
\hline
\multicolumn{1}{c}{Instance} & \multicolumn{1}{c}{$|\mathcal{N}|$} & \multicolumn{1}{c}{$|\mathcal{N}_{source}|$} & \multicolumn{1}{c}{$|\mathcal{N}_{sink}|$} & \multicolumn{1}{c}{$|\mathcal{P}|$}& \multicolumn{1}{c}{$|\mathcal{C}|$} & \multicolumn{1}{c}{$|\mathcal{V}|$} & \multicolumn{1}{c}{$|\mathcal{R}|$} \\ \hline
Gaslib-11 &11	&3	&3&	8&	2&	0&	1
  \\ \hline
Gaslib-24 &24	&3	 &5 &	21 &	3 &	1 &	0 
    \\ \hline
Gaslib-39  & 39&	2&	5	&28&	2	&4	&7       \\ \hline
Gaslib-40   &40	&3	&29	&39	&6	&0	&0
 \\ \hline
Gaslib-134 & 134&	3	&45&	131&	1	&1	&0
 \\ \hline
\end{tabular}}}
\caption{The network elements in GasLib instances.}
\label{table:gaslibInstances}
\end{center}
\end{table}
\subsection{Computational Setting}
All computational experiments are conducted on a 64-bit workstation with two Intel(R) Xeon(R) Gold 6248R CPU (3.00GHz) processors. The workstation runs on the Windows operations system and has 256 GB 2993 Mhz RAM. All mixed-integer convex relaxations for the multi-period gas storage optimization problem are solved by using the Gurobi Optimizer 11.0.1 \citep{Gurobi}. The $\texttt{BarHomogeneous}$ and $\texttt{BarConvTol}$ parameters are set to $1$ and $5 \times 10^{-6}$ in Gurobi, respectively. After fixing the binary variables obtained from these relaxations, the resulting NLP is solved by the interior point solver Ipopt 3.11.1 \citep{ipopt}.
For these NLPs, the $\texttt{bound$\_$relax$\_$factor}$ and $\texttt{max$\_$iter}$ parameters are set to $0$ and $5 \times 10^4$ in Ipopt, respectively.  We also compare the computational performance of our approaches by solving the nonconvex MINLP problem \eqref{eq:MINLP} with the global optimization solver BARON Solver 22.4.20 (hereafter referred to as \texttt{BARON}) \citep{baron}. The time limit within each solver is set to one hour. Time in all computational experiments is measured in seconds. 

In order to assess the scalability of our approaches, we also stress the gas network by scaling the gas demand from the set $\{0.1, 0.5, 1.0, 1.5 \}$ for each test instance. Moreover, we also conduct our experiments under different objectives that are used in the gas network optimization problems. In particular, we minimize the following objective functions with positive cost parameters:
\begin{itemize}
    \item $g_1 = \sum_{(i,j) \in \mathcal{C}} \sum_{t \in \mathcal{T}} ( C^{fc}_{ij} |l_{ijt}| + C^{up}_{ij} u_{ijt}) + \sum_{i \in \mathcal{N}} \sum_{j \in \mathcal{S}(i)} \sum_{t \in \mathcal{T}}  C^{wd}_j s_{kt}^-$. 
    \item $g_2 = \sum_{(i,j) \in \mathcal{C}} \sum_{t \in \mathcal{T}}  ( \gamma_1 |\pi_{it} - \pi_{jt}| + \gamma_2| (\pi_{it} -\pi_{i,t-1}) - (\pi_{jt} -\pi_{j,t-1}) |)$. 
    \item $g_3 = \sum_{(i,j) \in \mathcal{C}} \sum_{t \in \mathcal{T}}  C^{fc}_{ij} |l_{ijt}|$.
\end{itemize}
We note that the absolute values in these objective functions can be easily handled by common linearization approaches. The objective function $g_1$ and $g_3$ are the operational costs whereas $g_2$ is a simplified version of energy consumption by compressors, see, e.g., \citep{burlacu2019}. We present our computational results by presenting the following metrics that are obtained within Algorithm~\ref{alg:MINLP_GlobalAlgorithm}:
\begin{itemize}
    \item Dual bound: Lower bound obtained by \texttt{convexRelaxation}.
    \item Primal bound: Upper bound obtained by the NLP.
    \item Relaxation time: Wall-clock time in \texttt{convexRelaxation}.
    \item Total time: Total wall-clock time.
    \item Gap ($\%$): The relative optimality gap in percentages.
\end{itemize}
 The abbreviation ``inf'' (under column ``Dual bound'') is used if the mixed-integer convex relaxation is proven to be infeasible. Similarly, under ``Primal bound'' columns, we use the abbreviation ``inf'' if the local solver converges to a locally infeasible solution. Under ``Relaxation time'' column, we use the abbreviation ``TL'' if the one-hour time limit has been reached. 

Next, we explain the implementation of the disjunctive formulations of \eqref{eq:MINLP} for Gurobi and \texttt{BARON}. For each $(i,j) \in \mathcal{A}$, we define two binary variables $x_{ij}^+$ and $x_{ij}^-$ that satisfy the relation $x_{ij} = x_{ij}^+ + x_{ij}^-$ where $x_{ij} \in [0,1]$. For both Gurobi and \texttt{BARON}, we replace the disjunctive formulations given by constraints \eqref{eq:activeElementConditions_1}, \eqref{eq:MINLP_activeElementRelative1} and \eqref{eq:MINLP_activeElementRelative2} for each (control) valve $(i,j) \in \mathcal{A} \setminus \mathcal{C}$ by the following system of linear equations:
\begin{subequations} 
\begin{align*}
&\hspace{0.5em} \pi_j - \munderbar{r}_{ij} \pi_i  \ge (\munderbar{\pi}_j - \munderbar{r}_{ij} \bar{\pi}_i) (1-x_{ij}^+ ), \ \pi_j - \bar{r}_{ij} \pi_i  \le (\bar{\pi}_j - \bar{r}_{ij} \munderbar{\pi}_i) (1-x_{ij}^+ ) &(&i,j) \in \mathcal{A} \setminus \mathcal{C}\\
&\hspace{0.5em}  \pi_i - \munderbar{r}_{ij} \pi_j  \ge (\munderbar{\pi}_i - \munderbar{r}_{ij} \bar{\pi}_j) (1-x_{ij}^- ), \ \pi_i - \bar{r}_{ij} \pi_j  \le (\bar{\pi}_i - \bar{r}_{ij} \munderbar{\pi}_j) (1-x_{ij}^- ) &(&i,j) \in \mathcal{A} \setminus \mathcal{C}\\
&\hspace{0.5em} x_{ij}^- \munderbar{f}_{ij} \le f_{ij} \le \bar{f}_{ij} x_{ij}^+  &(&i,j) \in \mathcal{A} \setminus \mathcal{C}.
\end{align*}
\end{subequations}
For the global solver \texttt{BARON}, we also recast the nonconvex system given by constraints \eqref{eq:activeElementConditions_1}, \eqref{eq:CompressorConditions_1}, \eqref{eq:MINLP_activeElementRelative1}-\eqref{eq:MINLP_compressorLoss2} for each compressor $(i,j) \in \mathcal{C}$ by the following system of nonlinear equations:
\begin{subequations}  \label{eq:constraintsBaron}
\begin{align}
&\hspace{0.5em} l_{ij} = H_{ij} \big(({\pi_{i}} / {\pi_{j}})^{\hat{\kappa}}x_{ij}^+ + ({\pi_{j}} /{\pi_{i}})^{\hat{\kappa}} x_{ij}^- -1 \big) f_{ij} &(&i,j) \in \mathcal{C}\\
&\hspace{0.5em} \pi_j - \munderbar{r}_{ij} \pi_i  \ge (\munderbar{\pi}_j - \munderbar{r}_{ij} \bar{\pi}_i) (1-x_{ij}^+ ), \ \pi_j - \bar{r}_{ij} \pi_i  \le (\bar{\pi}_j - \bar{r}_{ij} \munderbar{\pi}_i) (1-x_{ij}^+ ) &(&i,j) \in \mathcal{C}\\
&\hspace{0.5em}  \pi_i - \munderbar{r}_{ij} \pi_j  \ge (\munderbar{\pi}_i - \munderbar{r}_{ij} \bar{\pi}_j) (1-x_{ij}^- ), \ \pi_i - \bar{r}_{ij} \pi_j  \le (\bar{\pi}_i - \bar{r}_{ij} \munderbar{\pi}_j) (1-x_{ij}^- ) &(&i,j) \in \mathcal{C}\\
&\hspace{0.5em} x_{ij}^- \munderbar{f}_{ij} \le f_{ij} \le \bar{f}_{ij} x_{ij}^+  &(&i,j) \in \mathcal{C}.
\end{align}
\end{subequations}
In Table~\ref{table:summaryFormulations}, we summarize the different formulations that are used in our computational experiments. We remind the reader that we use the MISOCP relaxation of each formulation for \texttt{R1}, \texttt{R2} and \texttt{R3} in Gurobi.
\begin{table}[h!] 
\centering
\begin{tabular}{ccccccccc}
\hline
              Nonconvex sets       &  & \texttt{R1}    &    &  \texttt{R2}   &    &  \texttt{R3}   &  &  \texttt{BARON} \\ \cline{1-1} \cline{3-3} \cline{5-5} \cline{7-7} \cline{9-9}  
Weymouth set ($\mathcal{X}_{ij}$) for pipes    &  & \textup{Conv}($\mathcal{X}_{ij}$) (Theorem \ref{theorem:convexHull_pipes})                     &  &             $\mathcal{X}_{\text{socr}}$        &  & $\mathcal{X}_{\text{poly}}$                       &  &   $\mathcal{X}_{ij}$               \\
\multicolumn{9}{c}{}                                                                                                      \\
Compression set  ($\mathcal{L}_{ij}$) for compressors   &  & \multicolumn{5}{c}{$\mathcal{Y}_{ij}$ (Theorem \ref{theorem:compressorOuterApprox})}                                         &  & \eqref{eq:constraintsBaron}                    \\ \hline
\end{tabular}
\caption{The summary of different formulations used in solvers.}
\label{table:summaryFormulations}
\end{table}
\subsection{Computational Results}
In this section, we present the results of our computational study. In particular, we consider small-scale gas networks with $11, 24, 39$, and $40$ nodes and medium-scale real-world gas networks with $134$ nodes. 
In Tables~\ref{tab:resultsGasLib11}, \ref{tab:resultsGasLib24}, \ref{tab:resultsGasLib39} and \ref{tab:resultsGasLib40}, we report our detailed results for small-scale networks for the single nomination taken from the literature. For GasLib-134, there are multiple nominations provided in three categories, namely, warm, mild and cold. We pick a single nomination from each category and present our computational results under different temperature conditions for GasLib-134. In particular, we pick nominations ``2015-08-23", ``2015-05-24" and ``2016-02-10" for mild, warm and cold temperatures, respectively. In Tables~\ref{tab:resultsGasLib134warm}, \ref{tab:resultsGasLib134mild}, and \ref{tab:resultsGasLib134cold}, we report our detailed results for GasLib-134 under warm, mild and cold nominations, respectively. In a preliminary study, we ran each computational experiment by omitting the MISOCP relaxation for the compression set. We observed that the dual bounds given by relaxation methods deteriorate although there are no significant differences in the dual bounds across relaxations. Moreover, we also observed that there are some instances in which the local solver was unable to produce primal bounds when initialized with starting points produced by those methods without the MISOCP relaxation for compressors. As a result, we include the MISOCP relaxation of the outer-approximation in Theorem~\ref{theorem:compressorOuterApprox} for relaxation methods \texttt{R1}, \texttt{R2} and \texttt{R3}.


For the GasLib-11 network, we observe that the performance of \texttt{R1} is consistently better than other methods in terms of optimality gaps for all instances. We also observe that \texttt{BARON} reaches the time limit and is unable to provide near-global solutions under the objectives $g_1$ and $g_3$ although the primal bounds of \texttt{BARON} are very close to the dual bounds produced by the convex relaxations. Although \texttt{BARON} provides better primal bounds under $g_2$ compared to the relaxation methods, its computational performance decreases and cannot prove optimality under the time limit as the congestion level increases. We see that although the differences in total time are mostly negligible among all convex relaxation methods, \texttt{R2} and \texttt{R3} are computationally faster than \texttt{R1}. However, our method consistently outperforms these methods in terms of both dual and primal bounds and is able to provide near-global solutions for all instances.
\afterpage{
\begin{landscape}
\begin{table}[p!]
    \centering
        \scalebox{0.75}{ 
          
\begin{tabular}{ccccccccccccccccccccc}
\hline
\multirow{2}{*}{Stress} &  & \multirow{2}{*}{Solver} &  & \multicolumn{5}{c}{$g_1$}                                              &  & \multicolumn{5}{c}{$g_2$}                                              &  & \multicolumn{5}{c}{$g_3$}                                              \\ \cline{5-9} \cline{11-15} \cline{17-21} 
                        &  &                         &  & Dual bound & Primal bound & Relaxation time & Total time & Gap (\%) &  & Dual bound & Primal bound & Relaxation time & Total time & Gap (\%) &  & Dual bound & Primal bound & Relaxation time & Total time & Gap (\%) \\ \cline{1-1} \cline{3-3} \cline{5-9} \cline{11-15} \cline{17-21} 
\multirow{4}{*}{0.1}    &  & \texttt{R1}               &  & \bf  11585.95   & 11585.95     & 7.34            & 21.13       & \bf 0.00     &  &  23.91      & 24.39        & 4.39            & 5.50       & 1.98     &  &\bf  332.10     & 332.10       & 2.27            & 2.88       & \bf  0.00     \\
                        &  & \texttt{R2}               &  & 11585.91   & 11663.14     & 5.89            & 7.78       & 0.66     &  & 23.19      & 24.39        & 3.23            & 4.22       & 4.92     &  & 332.07     & 334.99       & 1.41            & 2.74       & 0.87     \\
                        &  & \texttt{R3}               &  & 11585.91   & 11586.45     & 1.70            & 3.67       &   0.00     &  & 23.19      & inf          & 1.16            & 5.17       & -   &  & 332.07     & inf          & 1.22            & 3.22       & -   \\
                        &  & \texttt{BARON}            &  & 11393.64   &   11585.95     & -               & TL         & 1.66     &  & \bf 24.39      & 24.39        & -               & 33.77      & \bf 0.01     &  & 134.84     & 332.10       & -               & TL         & 59.40    \\ \cline{1-1} \cline{3-3} \cline{5-9} \cline{11-15} \cline{17-21} 
\multirow{4}{*}{0.5}    &  & \texttt{R1}               &  & \bf 12985.73   & 12986.74     & 7.22            & 8.91       & \bf 0.01     &  & 27.68      & 29.09        & 3.96            & 5.30       & 4.86     &  & \bf 1660.49    & 1660.49      & 1.94            & 3.39       &\bf  0.00     \\
                        &  & \texttt{R2}               &  & \bf 12985.73   & 12988.73     & 1.51            & 2.72       & 0.02     &  & 27.68      & 29.12        & 2.75            & 3.98       & 4.94     &  & 1660.11    & 1660.49      & 1.09            & 2.51       & 0.02     \\
                        &  & \texttt{R3}               &  & \bf 12985.73   & 12987.19     & 1.13            & 2.87       & 
  0.01     &  & 27.68      & 30.35        & 2.03            & 4.13       & 8.82     &  &\bf  1660.49    & 1660.49      & 1.09            & 2.17       & \bf 0.00     \\
                        &  & \texttt{BARON}            &  & 12009.43   & 12985.73     & -               & TL         & 7.52     &  & \bf  28.34      & 28.61        & -               & TL         & \bf  0.94     &  & 665.27     & 1660.49      & -               & TL         & 59.94    \\ \cline{1-1} \cline{3-3} \cline{5-9} \cline{11-15} \cline{17-21} 
\multirow{4}{*}{1.0}    &  & \texttt{R1}               &  & \bf  14735.46   & 14735.46     & 4.30            & 6.39       & \bf 0.00     &  & \bf 33.81      & 35.83        & 65.99           & 67.97      & \bf 5.62     &  & \bf 3320.98    & 3320.98      & 1.99            & 4.72       & \bf 0.00     \\
                        &  & \texttt{R2}               &  & \bf  14735.46   & 14740.06     & 1.43            & 3.30       & 0.03     &  & 33.32      & 35.66        & 58.55           & 62.75      & 6.57     &  & \bf 3320.98    & 3320.98      & 1.42            & 2.73       & \bf 0.00     \\
                        &  & \texttt{R3}               &  & \bf  14735.46   & 14735.84     & 1.10            & 2.42       & 0.00     &  & 32.42      & 35.52        & 20.45           & 21.67      & 8.74     &  & \bf 3320.98    & inf          & 1.04            & 2.13       & -   \\
                        &  & \texttt{BARON}            &  & 12743.96   & 14738.38     & -               & TL         & 13.53    &  & 33.11      & 35.67        & -               & TL         & 7.16     &  & 1324.65    & 3320.98      & -               & TL         & 60.11    \\ \cline{1-1} \cline{3-3} \cline{5-9} \cline{11-15} \cline{17-21} 
\multirow{4}{*}{1.5}    &  & \texttt{R1}               &  & 16485.19   & 16493.52
     & 4.18            & 7.12       & 0.05     &  & \bf 45.65      & 46.39        & 73.85           & 75.18      &\bf  1.58     &  & \bf 4981.47    & 4981.47     & 1.73            & 3.05       & \bf 0.00     \\
                        &  & \texttt{R2}               &  & \bf 16485.26   & 16508.25     & 1.19            & 3.16       & 0.14     &  & \bf 45.65      & 46.39        & 47.70           & 50.12      &\bf  1.58     &  & \bf 4981.47    & 4981.47      & 1.31            & 2.41       & \bf 0.00     \\
                        &  & \texttt{R3}               &  & 16485.19   & 16492.95
     & 1.37            & 4.22       & \bf  0.05     &  & 43.36      & 46.39        & 16.04           & 17.25      & 6.52     &  & \bf 4981.47    & 4981.47      & 1.37            & 2.91       &\bf  0.00     \\
                        &  & \texttt{BARON}            &  & 13816.16   & 16485.18     & -               & TL         & 16.19    &  & 44.30      & 45.70        & -               & TL         & 3.06     &  & 2273.25    & 4981.47      & -               & TL         & 54.37    \\ \hline
\end{tabular}

}
    \caption{Computational results on GasLib-11.}
    \label{tab:resultsGasLib11}
\end{table}
\begin{table}[p!]
    \centering
        \scalebox{0.75}{ 
          
\begin{tabular}{ccccccccccccccccccccc}
\hline
\multirow{2}{*}{Stress} &  & \multirow{2}{*}{Solver} &  & \multicolumn{5}{c}{$g_1$}                                              &  & \multicolumn{5}{c}{$g_2$}                                              &  & \multicolumn{5}{c}{$g_3$}                                              \\ \cline{5-9} \cline{11-15} \cline{17-21} 
                        &  &                         &  & Dual bound & Primal bound & Relaxation time & Total time & Gap (\%) &  & Dual bound & Primal bound & Relaxation time & Total time & Gap (\%) &  & Dual bound & Primal bound & Relaxation time & Total time & Gap (\%) \\ \cline{1-1} \cline{3-3} \cline{5-9} \cline{11-15} \cline{17-21} 
\multirow{3}{*}{0.1}    &  & \texttt{R1}               &  & \bf 17934.98   & 17934.98     & 1.97            & 5.71       & \bf 0.00     &  & \bf 31.28      & 35.04        & 2.04            & 3.36       & 10.74    &  & \bf 1048.60    & 1048.60      & 1.97            & 15.30      & \bf 0.00     \\
                        &  & \texttt{R3}               &  & \bf  17934.98   & 17934.98     & 2.03            & 11.08      & \bf 0.00     &  & \bf  31.28      & 35.04        & 2.02            & 3.37       & 10.74    &  & \bf 1048.60    & 1048.60      & 2.00            & 15.48      & \bf 0.00     \\
                        &  & \texttt{BARON}            &  & 17040.26   & 17934.98     & -               & TL         & 4.99     &  & \bf 31.28      & 31.28        & -               & 2.45       &\bf  0.00     &  & 119.68     & 1048.60      & -               & TL         & 88.59    \\ \cline{1-1} \cline{3-3} \cline{5-9} \cline{11-15} \cline{17-21} 
\multirow{3}{*}{0.5}    &  & \texttt{R1}               &  & \bf  22258.90   & 22258.90     & 1.98            & 5.84       & \bf  0.00     &  & \bf 30.04      & 32.13        & 2.24            & 5.63       & 6.51     &  &\bf  5242.98    & 5242.98      & 1.99            & 5.80       & \bf 0.00     \\
                        &  & \texttt{R3}               &  & \bf  22258.90   & 22258.90     & 1.95            & 7.08       & \bf  0.00     &  & \bf 30.04      & 32.13        & 2.26            & 6.38       & 6.51     &  & \bf 5242.98    & 5242.98      & 1.92            & 5.62       &\bf  0.00     \\
                        &  & \texttt{BARON}            &  & 17380.53   & 22258.90     & -               & TL         & 21.92    &  &\bf  30.04      & 30.04        & -               & 17.94      & \bf  0.01     &  & 389.67     & 5242.98      & -               & TL         & 92.57    \\ \cline{1-1} \cline{3-3} \cline{5-9} \cline{11-15} \cline{17-21} 
\multirow{3}{*}{1.0}    &  & \texttt{R1}               &  &\bf  27663.79   & 27663.79     & 2.02            & 5.94       & \bf 0.00     &  & \bf 31.60      & 31.85        & 2.02            & 5.86       & 0.78     &  &\bf  10485.97   & 10485.97     & 1.95            & 6.73       & \bf 0.00     \\
                        &  & \texttt{R3}               &  & \bf  27663.79   & 27663.79     & 2.05            & 5.97       & \bf 0.00     &  & \bf 31.60      & 31.85        & 2.27            & 5.45       & 0.78     &  & \bf 10485.97   & 10485.97     & 1.92            & 6.77       & \bf 0.00     \\
                        &  & \texttt{BARON}            &  & 17729.56   & 27663.79     & -               & TL         & 35.91    &  & 31.59      & 31.60        & -               & 11.98      & \bf 0.01     &  & 526.13     & 10485.97     & -               & TL         & 94.98    \\ \cline{1-1} \cline{3-3} \cline{5-9} \cline{11-15} \cline{17-21} 
\multirow{3}{*}{1.5}    &  & \texttt{R1}               &  & \bf 33068.69   & 33068.69     & 1.99            & 4.72       &\bf  0.00     &  &\bf  35.61      & 35.96        & 2.29            & 9.23       & 0.98     &  & \bf 15728.95   & 15728.95     & 1.92            & 5.29       & \bf 0.00     \\
                        &  & \texttt{R3}               &  &\bf  33068.69   & 33068.69     & 1.95            & 5.80       & \bf 0.00     &  & \bf 35.61      & 35.96        & 2.57            & 6.08       & 0.98     &  & \bf 15728.95   & 15728.95     & 1.94            & 5.45       & \bf 0.00     \\
                        &  & \texttt{BARON}            &  & 18133.35   & 33068.69     & -               & TL         & 45.16    &  & 35.60      & 35.61        & -               & 15.99      & \bf 0.01     &  & 792.86     & 15728.95     & -               & TL         & 94.96    \\ \hline
\end{tabular}

}
    \caption{Computational results on GasLib-24.}
    \label{tab:resultsGasLib24}
\end{table}
\end{landscape}
}
\afterpage{
\begin{landscape}
\begin{table}[p!]
    \centering
        \scalebox{0.75}{ 
          
\begin{tabular}{ccccccccccccccccccccc}
\hline
\multirow{2}{*}{Stress} &  & \multirow{2}{*}{Solver} &  & \multicolumn{5}{c}{$g_1$}                                              &  & \multicolumn{5}{c}{$g_2$}                                              &  & \multicolumn{5}{c}{$g_3$}                                              \\ \cline{5-9} \cline{11-15} \cline{17-21} 
                        &  &                         &  & Dual bound & Primal bound & Relaxation time & Total time & Gap (\%) &  & Dual bound & Primal bound & Relaxation time & Total time & Gap (\%) &  & Dual bound & Primal bound & Relaxation time & Total time & Gap (\%) \\ \cline{1-1} \cline{3-3} \cline{5-9} \cline{11-15} \cline{17-21} 
\multirow{4}{*}{0.1}    &  & \texttt{R1}               &  & \bf 13018.55   & 13018.55     & 5.02            & 13.54      & \bf  0.00     &  & \bf 13.09      & 15.99        & 20.63           & 37.92      & 18.12    &  & \bf 1673.82    & 1673.82      & 4.43            & 14.79      &\bf  0.00     \\
                        &  & \texttt{R2}               &  & \bf 13018.55   & 13031.37     & 2.33            & 11.80      & 0.10     &  & \bf 13.09      & 32.27        & 4.58            & 360.75     & 59.44    &  & \bf 1673.82    & 1673.82      & 2.22            & 30.17      & \bf 0.00     \\
                        &  & \texttt{R3}               &  & \bf 13018.55   & inf          & 1.72            & 13.49      & -   &  & \bf 13.09      & inf          & 1.89            & 19.26      & -   &  & 1673.81    & 1673.82      & 1.71            & 17.90      & 0.00     \\
                        &  & \texttt{BARON}            &  & 11344.73   & 13018.55     &     -            & TL         & 12.86    &  & \bf 13.09      & 13.12        &            -     & TL         & \bf  0.26     &  &      -      & inf          &         -        & 81.95      &        -  \\ \cline{1-1} \cline{3-3} \cline{5-9} \cline{11-15} \cline{17-21} 
\multirow{4}{*}{0.5}    &  & \texttt{R1}               &  & \bf 20148.75   & 20148.75     & 4.50            & 10.74      &\bf  0.00     &  & \bf 13.09      & 27.34        & 2619.95         & 2642.07    & 52.12    &  & \bf 8369.08    & 8369.08      & 4.40            & 14.03      & \bf 0.00     \\
                        &  & \texttt{R2}               &  & \bf 20148.75   & inf          & 6.79            & 19.66      & -  &  & \bf  13.09      & 34.86        & TL              & 3663.08    & 62.45    &  & \bf 8369.08    & 8369.08      & 2.54            & 38.86      & \bf 0.00     \\
                        &  & \texttt{R3}               &  & \bf 20148.75   & 20268.14     & 1.71            & 38.30      & 0.59     &  & \bf 13.09      & 28.78        & 2.09            & 1133.24    & 54.52    &  & \bf 8369.08    & 8369.08      & 1.74            & 20.60      & \bf 0.00     \\
                        &  & \texttt{BARON}            &  & 11779.67   & 20148.75     &      -           & TL         & 41.54    &  & \bf  13.09      & 14.68        &      -           & TL         & \bf 10.81    &  &     -       & inf          &           -      & 42.66      &      -    \\ \cline{1-1} \cline{3-3} \cline{5-9} \cline{11-15} \cline{17-21} 
\multirow{4}{*}{1.0}    &  & \texttt{R1}               &  & \bf 29061.50   & 29080.17     & 5.53            & 26.19      &\bf  0.06     &  &   15.03      & 16.29        & TL              & 3662.87    & \bf 7.70     &  & \bf 16738.15   & 16738.15     & 5.08            & 16.00      & \bf 0.00     \\
                        &  & \texttt{R2}               &  &\bf  29061.50   & 29138.14     & 5.63            & 37.73      & 0.26     &  & 14.98      & 29.90        & TL              & 3645.26    & 49.90    &  &\bf  16738.15   & inf          & 2.76            & 14.31      & -   \\
                        &  & \texttt{R3}               &  & \bf 29061.50   & inf          & 1.81            & 12.27      & -   &  & 13.80      & 34.95        & 2.75            & 292.89     & 60.52    &  & \bf 16738.15   & inf          & 1.80            & 19.05      & -  \\
                        &  & \texttt{BARON}            &  & 12323.35   & 29069.24     &       -          & TL         & 57.61    &  & \bf 16.13      & 18.93        &     -            & TL         & 14.80    &  & 0.00       & 16738.15     &      -           & TL         & 100.00   \\ \cline{1-1} \cline{3-3} \cline{5-9} \cline{11-15} \cline{17-21} 
\multirow{4}{*}{1.5}    &  & \texttt{R1}               &  &\bf  37974.25   & 37985.86     & TL              & 3630.18    & \bf 0.03     &  & \bf 20.16      & 23.91        & TL              & 3786.89    & \bf 15.68    &  & \bf 25107.23   & 25107.23     & 114.36          & 139.28     & \bf 0.00     \\
                        &  & \texttt{R2}               &  &\bf  37974.25   & inf          & 177.45          & 189.00     & - &  & 19.48      & 32.70        & TL              & 3729.74    & 40.43    &  & 25106.97   & inf          & 2.43            & 25.62      & -   \\
                        &  & \texttt{R3}               &  & \bf 37974.25   & 38607.72     & 1.59            & 17.36      & 1.64     &  & 18.81      & 32.30        & 324.08          & 353.57     & 41.78    &  & \bf 25107.23   & 25107.23     & 1.58            & 23.36      & \bf 0.00     \\
                        &  & \texttt{BARON}            &  &            & inf          &           -      & 174.49     &       -   &  &      -      & inf          &          -       & 13.90      &  -        &  & 724.00     & 25107.22     &      -           & TL         & 97.12    \\ \hline
\end{tabular}

}
    \caption{Computational results on GasLib-39.}
    \label{tab:resultsGasLib39}
\end{table}
\begin{table}[p!]
    \centering
        \scalebox{0.75}{ 
          
\begin{tabular}{ccccccccccccccccccccc}
\hline
\multirow{2}{*}{Stress} &  & \multirow{2}{*}{Solver} &  & \multicolumn{5}{c}{$g_1$}                                              &  & \multicolumn{5}{c}{$g_2$}                                              &  & \multicolumn{5}{c}{$g_3$}                                              \\ \cline{5-9} \cline{11-15} \cline{17-21} 
                        &  &                         &  & Dual bound & Primal bound & Relaxation time & Total time & Gap (\%) &  & Dual bound & Primal bound & Relaxation time & Total time & Gap (\%) &  & Dual bound & Primal bound & Relaxation time & Total time & Gap (\%) \\ \cline{1-1} \cline{3-3} \cline{5-9} \cline{11-15} \cline{17-21} 
\multirow{4}{*}{0.1}    &  & \texttt{R1}               &  & 30402.54   & 30402.76     & TL              & 3652.87    & 0.00     &  & \bf 0.57       & 10.55        & TL              & 3626.27    & \bf  94.63    &  & 2183.14    & 2197.27      & 970.01          & 1006.11    & \bf 0.64     \\
                        &  & \texttt{R2}               &  & 30402.49   & 30557.82     & 1224.97         & 1272.33    & 0.51     &  & 0.51       & 66.50        & TL              & 3633.55    & 99.24    &  & 2178.19    & 2945.06      & TL              & 3690.41    & 26.04    \\
                        &  & \texttt{R3}               &  & \bf 30402.62   & 30402.76     & 2.95            & 49.94      & \bf  0.00     &  & 0.05       & 3.98         & 4.31            & 59.93      & 98.65    &  & \bf  2183.22    & 3164.83      & 3.22            & 60.60      & 31.02    \\
                        &  & \texttt{BARON}            &  & -          & inf          & -               & 533.76     & -        &  & -          & inf          & -               & 381.03     & -        &  & -          & inf          & -               & 104.43     & -        \\ \cline{1-1} \cline{3-3} \cline{5-9} \cline{11-15} \cline{17-21} 
\multirow{4}{*}{0.5}    &  & \texttt{R1}               &  & \bf  39653.82   & 39653.82     & 8.17            & 19.97      &\bf  0.00     &  & \bf  19.30      & 188.20       & TL              & 3648.39    & 89.74    &  & \bf 10916.85   & 10916.85     & 630.30          & 684.10     &\bf  0.00     \\
                        &  & \texttt{R2}               &  & \bf 39653.82   & 39653.82     & 101.83          & 255.37     & \bf  0.00     &  & 12.46      & 69.18        & TL              & 3681.64    & \bf  81.98    &  & 10862.79   & 12343.13     & TL              & 3706.31    & 11.99    \\
                        &  & \texttt{R3}               &  & 39651.30   & 39653.82     & 3.18            & 9.17       & 0.01     &  & 1.37       & 67.79        & 3.07            & 91.59      & 97.97    &  & 10914.33   & 10916.85     & 3.03            & 9.06       & 0.02     \\
                        &  & \texttt{BARON}            &  & -          & inf          & -               & 69.95      & -        &  & -          & inf          & -               & 82.79      & -        &  & -          & inf          & -               & 202.82     & -        \\ \cline{1-1} \cline{3-3} \cline{5-9} \cline{11-15} \cline{17-21} 
\multirow{4}{*}{1.0}    &  & \texttt{R1}               &  & 51217.52   & 51217.64     & 107.06          & 114.30     & 0.00     &  & -          & -            & TL              & 3656.69    & -        &  & \bf 21833.70   & 22997.54     & TL              & 4277.62    & 5.06     \\
                        &  & \texttt{R2}               &  & \bf 51217.62   & 51217.64     & TL              & 3631.30    & \bf  0.00     &  & -          & -            & TL              & 3618.00    & -        &  & 21833.61   & 23425.75     & 1873.25         & 1964.20    & 6.80     \\
                        &  & \texttt{R3}               &  & 51217.55   & 51217.64     & 3.04            & 11.64      & 0.00     &  & \bf  59.31      & 271.03       & 4.16            & 44.50      & \bf  78.12    &  & 21833.61   & 22296.61     & 3.26            & 65.98      & \bf  2.08     \\
                        &  & \texttt{BARON}            &  & -          & inf          & -               & 127.04     & -        &  & -          & inf          & -               & 100.29     & -        &  & -          & inf          & -               & 86.15      & -        \\ \cline{1-1} \cline{3-3} \cline{5-9} \cline{11-15} \cline{17-21} 
\multirow{4}{*}{1.5}    &  & \texttt{R1}               &  & inf        & -            & 9.23            & 9.23       & -        &  & inf        & -            & 13.32           & 13.32      & -        &  & inf        &           -   & 10.56           & 10.56      & -        \\
                        &  & \texttt{R2}               &  & -          & -            & TL              & 3612.88    & -        &  & inf        & -            & 161.81          & 161.81     & -        &  & inf        &     -         & 258.62          & 258.62     & -        \\
                        &  & \texttt{R3}               &  & inf        & -            & 2.92            & 2.92       & -        &  & inf        & -            & 2.84            & 2.84       & -        &  & inf        &       -       & 2.74            & 2.74       & -        \\
                        &  & \texttt{BARON}            &  & -          & inf          & -               & 37.15      & -        &  & -          & inf          & -               & 230.32     & -        &  &     -       & inf          &        -         & 59.18      & -        \\ \hline
\end{tabular}

}
    \caption{Computational results on GasLib-40.}
    \label{tab:resultsGasLib40}
\end{table}
\end{landscape}


\begin{landscape}

\begin{table}[p!]
    \centering
        \scalebox{0.9}{ 
          
\begin{tabular}{ccccccccccccccc}
\hline
\multirow{2}{*}{Stress} &  & \multirow{2}{*}{Solver} &  & \multicolumn{5}{c}{$g_1$}                                              &  & \multicolumn{5}{c}{$g_3$}                                              \\ \cline{5-9} \cline{11-15} 
                        &  &                         &  & Dual bound & Primal bound & Relaxation time & Total time & Gap (\%) &  & Dual bound & Primal bound & Relaxation time & Total time & Gap (\%) \\ \cline{1-1} \cline{3-3} \cline{5-9} \cline{11-15} 
\multirow{4}{*}{0.1}    &  &  \texttt{R1}                &  & \bf 5788.02    & 5788.02      & 22.49           & 103.82     & \bf  0.00     &  & \bf 158.35     & 158.35       & 12.21           & 61.41      &\bf  0.00     \\
                        &  &  \texttt{R2}                &  & \bf 5788.02    & 5788.02      & 80.14           & 440.38     &\bf  0.00     &  &\bf  158.35     & 158.35       & 4.22            & 50.31      & \bf 0.00     \\
                        &  &  \texttt{R3}                &  & \bf 5788.02    & 5788.31      & 2.91            & 29.48      & 0.01     &  & \bf 158.35     & inf          & 2.81            & 158.39     & -        \\
                        &  &  \texttt{BARON}             &  & 5706.12    & 5788.02      & -               & TL         & 1.41     &  & 37.77      & 158.35       & -               & TL         & 76.15    \\ \cline{1-1} \cline{3-3} \cline{5-9} \cline{11-15} 
\multirow{4}{*}{0.5}    &  &  \texttt{R1}                &  & \bf 6468.08    & 6468.08      & 33.04           & 108.32     & \bf 0.00     &  & \bf 791.74     & 791.74       & 11.08           & 57.16      & \bf 0.00     \\
                        &  &  \texttt{R2}                &  & \bf 6468.08    & 6468.08      & 83.02           & 326.09     &\bf  0.00     &  & \bf 791.74     & 791.74       & 4.83            & 35.11      & \bf 0.00     \\
                        &  &  \texttt{R3}                &  & \bf 6468.08    & 6468.08      & 3.13            & 302.65     &\bf  0.00     &  & \bf 791.74     & inf          & 2.80            & 30.95      & -        \\
                        &  &  \texttt{BARON}             &  & 5959.13    & 6468.08      & -               & TL         & 7.87     &  & 213.41     & 791.74       & -               & TL         & 73.05    \\ \cline{1-1} \cline{3-3} \cline{5-9} \cline{11-15} 
\multirow{4}{*}{1.0}    &  &  \texttt{R1}                &  & \bf 7318.16    & 7318.17      & 13.27           & 22.91      & \bf 0.00     &  & \bf 1583.49    & 1583.49      & 4.43            & 76.89      & \bf 0.00     \\
                        &  &  \texttt{R2}                &  & \bf 7318.16    & 7318.17      & 4.25            & 63.86      & \bf 0.00     &  & \bf 1583.49    & inf          & 2.69            & 66.61      & -        \\
                        &  &  \texttt{R3}                &  & \bf 7318.16    & 7318.17      & 2.80            & 42.07      &\bf  0.00     &  &\bf  1583.49    & 1583.49      & 13.01           & 47.25      & \bf 0.00     \\
                        &  &  \texttt{BARON}             &  & 6282.99    & 7318.17      & -               & TL         & 14.15    &  & 468.40     & 1583.49      & -               & TL         & 70.42    \\ \cline{1-1} \cline{3-3} \cline{5-9} \cline{11-15} 
\multirow{4}{*}{1.5}    &  &  \texttt{R1}                &  & \bf 8168.25    & 8168.25      & 14.61           & 57.90      & \bf 0.00     &  & \bf 2375.23    & 2375.23      & 15.38           & 132.75     & \bf  0.00     \\
                        &  &  \texttt{R2}                &  & \bf 8168.25    & 8168.25      & 5.20            & 50.44      & \bf 0.00     &  & \bf 2375.23    & 2375.23      & 5.11            & 75.13      & \bf 0.00     \\
                        &  &  \texttt{R3}                &  & \bf 8168.25    & 8168.25      & 2.96            & 49.44      & \bf 0.00     &  & \bf 
 2375.23    & inf          & 2.81            & 39.43      & -        \\
                        &  &  \texttt{BARON}             &  & 6620.45    & 8168.25      & -               & TL         & 18.95    &  & 742.09     & 2375.23      & -               & TL         & 68.76    \\ \hline
\end{tabular}

}
    \caption{Computational results on GasLib-134 under warm nomination.}
    \label{tab:resultsGasLib134warm}
\end{table}
\begin{table}[p!]
    \centering
        \scalebox{0.90}{ 
          

\begin{tabular}{ccccccccccccccc}
\hline
\multirow{2}{*}{Stress} &  & \multirow{2}{*}{Solver} &  & \multicolumn{5}{c}{$g_1$}                                              &  & \multicolumn{5}{c}{$g_3$}                                              \\ \cline{5-9} \cline{11-15} 
                        &  &                         &  & Dual bound & Primal bound & Relaxation time & Total time & Gap (\%) &  & Dual bound & Primal bound & Relaxation time & Total time & Gap (\%) \\ \cline{1-1} \cline{3-3} \cline{5-9} \cline{11-15} 
\multirow{4}{*}{0.1}    &  &  \texttt{R1}                &  & \bf 5768.87    & 5768.87      & 90.78           & 222.79     & \bf 0.00     &  & \bf 140.89     & 140.89       & 12.40           & 37.25      & \bf 0.00     \\
                        &  &  \texttt{R2}                &  & \bf 5768.87    & inf          & 88.90           & 197.99     & -        &  & \bf 140.89     & 140.91       & 4.67            & 59.23      & 0.01     \\
                        &  &  \texttt{R3}                &  & \bf 5768.87    & 5768.87      & 2.90            & 60.08      &\bf  0.00     &  & \bf 140.89     & inf          & 2.68            & 95.96      & -        \\
                        &  &  \texttt{BARON}             &  & 5694.70    & 5768.87      &                 & TL         & 1.29     &  & 48.18      & 140.89       & -               & TL         & 65.80    \\ \cline{1-1} \cline{3-3} \cline{5-9} \cline{11-15} 
\multirow{4}{*}{0.5}    &  &  \texttt{R1}                &  & \bf 6372.35    & 6372.35      & 47.97           & 59.59      & \bf 0.00     &  & \bf 704.47     & 704.47       & 11.40           & 88.23      & \bf 0.00     \\
                        &  &  \texttt{R2}                &  & \bf 6372.35    & 6372.35      & 1360.48         & 1524.30    & \bf 0.00     &  & \bf 704.47     & 706.76       & 4.15            & 50.06      & 0.32     \\
                        &  &  \texttt{R3}                &  & \bf 6372.35    & 6372.35      & 3.24            & 204.03     &\bf  0.00     &  & \bf 704.47     & inf          & 2.78            & 25.18      & -        \\
                        &  &  \texttt{BARON}             &  & 5909.59    & 6372.35      & -               & TL         & 7.26     &  & 190.81     & 704.47       & -               & TL         & 72.91    \\ \cline{1-1} \cline{3-3} \cline{5-9} \cline{11-15} 
\multirow{4}{*}{1.0}    &  &  \texttt{R1}                &  & \bf 7126.70    & 7126.71      & 12.23           & 60.37      & \bf 0.00     &  & \bf 1408.94    & 1408.94      & 13.56           & 74.83      &\bf  0.00     \\
                        &  &  \texttt{R2}                &  & \bf 7126.70    & 7126.71      & 5.01            & 16.16      &\bf  0.00     &  & \bf 1408.94    & 1408.95      & 4.40            & 44.09      & 0.00     \\
                        &  &  \texttt{R3}                &  & \bf 7126.70    & 7128.98      & 2.69            & 87.34      & 0.03     &  & \bf 1408.94    & inf          & 2.80            & 17.25      & -        \\
                        &  &  \texttt{BARON}             &  & 6200.02    & 7126.71      & -               & TL         & 13.00    &  & 382.52     & 1408.94      & -               & TL         & 72.85    \\ \cline{1-1} \cline{3-3} \cline{5-9} \cline{11-15} 
\multirow{4}{*}{1.5}    &  &  \texttt{R1}                &  & \bf 7881.06    & 7881.06      & 13.09           & 28.21      & \bf 0.00     &  & \bf 2113.41    & 2113.41      & 14.88           & 92.64      & \bf 0.00     \\
                        &  &  \texttt{R2}                &  & \bf 7881.06    & 7881.06      & 5.04            & 32.70      & \bf 0.00     &  & \bf  2113.41    & 2113.41      & 4.96            & 45.44      & \bf 0.00     \\
                        &  &  \texttt{R3}                &  & \bf 7881.06    & 7881.06      & 2.78            & 27.38      &\bf  0.00     &  &\bf  2113.41    & inf          & 2.93            & 37.68      & -        \\
                        &  &  \texttt{BARON}             &  & 6491.19    & 7881.06      & -               & TL         & 17.64    &  & 521.14     & 2113.41      & -               & TL         & 75.34    \\ \hline
\end{tabular}

}
    \caption{Computational results on GasLib-134 under mild nomination.}
    \label{tab:resultsGasLib134mild}
\end{table}
\end{landscape}
}
\afterpage{
\begin{landscape}
\begin{table}[p!]
    \centering
        \scalebox{0.90}{ 
          
\begin{tabular}{ccccccccccccccc}
\hline
\multirow{2}{*}{Stress} &  & \multirow{2}{*}{Solver} &  & \multicolumn{5}{c}{$g_1$}                                              &  & \multicolumn{5}{c}{$g_3$}                                              \\ \cline{5-9} \cline{11-15} 
                        &  &                         &  & Dual bound & Primal bound & Relaxation time & Total time & Gap (\%) &  & Dual bound & Primal bound & Relaxation time & Total time & Gap (\%) \\ \cline{1-1} \cline{3-3} \cline{5-9} \cline{11-15} 
\multirow{4}{*}{0.1}    &  &  \texttt{R1}                &  & \bf 5871.64    & 5871.64      & 127.08          & 207.25     & \bf 0.00     &  & \bf  223.62     & 223.62       & 12.32           & 73.45      &\bf  0.00     \\
                        &  &  \texttt{R2}                &  & \bf 5871.64    & inf          & 36.27           & 236.52     & -        &  & \bf 223.62     & inf          & 4.26            & 49.70      & -        \\
                        &  &  \texttt{R3}                &  & \bf 5871.64    & 5871.64      & 2.97            & 23.85      & \bf 0.00     &  & \bf 223.62     & inf          & 2.69            & 402.43     & -        \\
                        &  &  \texttt{BARON}             &  & -          & inf          & -               & 13.67      & -        &  & 143.99     & 223.62       & -               & TL         & 35.61    \\ \cline{1-1} \cline{3-3} \cline{5-9} \cline{11-15} 
\multirow{4}{*}{0.5}    &  &  \texttt{R1}                &  & 6886.19    & 6886.20      & 12.32           & 102.78     & 0.00     &  & \bf 1118.09    & 1118.09      & 12.73           & 128.06     & \bf 0.00     \\
                        &  &  \texttt{R2}                &  & \bf 6886.20    & 6886.20      & 67.93           & 149.98     &\bf  0.00     &  & \bf 1118.09    & inf          & 4.17            & 299.52     & -        \\
                        &  &  \texttt{R3}                &  & \bf 6886.20    & 6886.20      & 3.03            & 30.02      &\bf  0.00     &  & \bf 
 1118.09    & inf          & 2.80            & 124.10     & -        \\
                        &  &  \texttt{BARON}             &  & -          & inf          & -               & 157.98     & -        &  & 766.03     & 1118.09      & -               & TL         & 31.49    \\ \cline{1-1} \cline{3-3} \cline{5-9} \cline{11-15} 
\multirow{4}{*}{1.0}    &  &  \texttt{R1}                &  & \bf 8154.39    & 8154.39      & 13.44           & 106.25     & \bf 0.00     &  & \bf 2236.18    & 2236.18      & 12.62           & 134.56     & \bf 0.00     \\
                        &  &  \texttt{R2}                &  & 8154.38    & 8154.39      & 5.34            & 33.20      & 0.00     &  & \bf 2236.18    & inf          & 4.82            & 234.28     & -        \\
                        &  &  \texttt{R3}                &  & 8154.38    & 8154.39      & 2.85            & 31.31      & 0.00     &  & \bf 2236.18    & inf          & 2.75            & 37.82      & -        \\
                        &  &  \texttt{BARON}             &  & -          & inf          & -               & 55.34      & -        &  & 1356.40    & 2236.18      & -               & TL         & 39.34    \\ \cline{1-1} \cline{3-3} \cline{5-9} \cline{11-15} 
\multirow{4}{*}{1.5}    &  &  \texttt{R1}                &  & \bf 9422.59    & 9422.59      & 190.59          & 223.44     & \bf 0.00     &  &\bf  3354.27    & 3354.27      & 13.92           & 74.51      &\bf  0.00     \\
                        &  &  \texttt{R2}                &  & \bf 9422.59    & 11413.30     & 60.23           & 91.12      & 17.44    &  & \bf 3354.27    & inf          & 5.66            & 194.27     & -        \\
                        &  &  \texttt{R3}                &  & 9422.58    & 9422.59      & 2.97            & 61.81      & 0.00     &  & \bf 3354.27    & inf          & 2.85            & 117.10     & -        \\
                        &  &  \texttt{BARON}             &  & -          & inf          & -               & 77.63      & -        &  & 1512.13    & 12014.83     & -               & TL         & 87.41    \\ \hline
\end{tabular}

}
    \caption{Computational results on GasLib-134 under cold nomination.}
    \label{tab:resultsGasLib134cold}
\end{table}
\end{landscape}
}
For the GasLib-24 network, we restrict ourselves to the positive gas flows for pipes and exclude \texttt{R2} (see, Remark \ref{remark:posFlows} for unidirectional flows). We see that convex relaxation methods outperform \texttt{BARON} and are able to produce near-global solutions under the objectives $g_1$ and $g_3$. Moreover, our method \texttt{R1} is computationally faster than \texttt{R3} for all instances whereas \texttt{BARON} reaches the time limit. Under the objective $g_2$, we observe that \texttt{BARON} is able to terminate within the time limit and produce near-global solutions. On the other hand, convex relaxation methods are also able to produce the same solutions as \texttt{BARON}; however, they cannot prove global optimality due to the locally optimal solutions they produce. Still, the optimality gaps produced by these relaxations decrease considerably as the stress level in the network increases.

For the GasLib-39 network, we first note that we encounter some persistent issues with \texttt{BARON}. In particular, \texttt{BARON} gives inconsistent results for all instances. For example, \texttt{BARON} claims that the instances under stress levels $0.1$ and $0.5$ are infeasible under the objective $g_3$ whereas it is able to produce primal bounds under $g_1$ and $g_2$ for the same instances. The same infeasibility issue also occurs as the network congestion increases under the objectives $g_1$ and $g_2$. In terms of dual bounds, the convex relaxation methods produce similar results. In terms of primal bounds, \texttt{R1} consistently produces high-quality feasible solutions and better optimality gaps compared to other relaxation methods. In the local solver, we also observe that there are cases where both \texttt{R2} and \texttt{R3} converge to a point of local infeasibility when initiated with the starting solutions produced by these relaxations. Still, \texttt{R1} offers high-quality warm-starting points for local solvers for all instances. 

For the GasLib-40 network, we encounter the same issues with \texttt{BARON}. Interestingly, it claims false infeasibility for all feasible instances. We try to overcome this by solving the underlying problem in a shorter planning horizon and decreasing the stress levels. Still, \texttt{BARON} falsely states upon termination that these instances are infeasible. On the other hand, it is able to claim infeasibility within the time limit under the stress level of $1.5$. For most of the instances, \texttt{R1} is able to produce considerably better optimality gaps compared to other methods within the time limit. Moreover, it also produces considerably tighter dual bounds compared to other relaxation methods under the objective $g_2$ for low congested networks. All convex relaxations detect infeasibility within the time limit; however, \texttt{R1} and \texttt{R3} are computationally faster than \texttt{R2}.

For the GasLib-134 network, we present our results under different nominations. We note that we encounter inconsistency issues with \texttt{BARON}. Interestingly, it produces inconsistent primal and dual bounds under the objective $g_2$. Thus, we omit them from our computational results. Although \texttt{BARON} reaches the time limit under warm and mild nominations, it is able to produce consistent results. However, 
 it falsely claims infeasibility under the objective $g_1$ in Table~\ref{tab:resultsGasLib134cold} although it is able to give consistent primal bounds under the objective $g_3$ for cold nomination. All convex relaxation methods terminate within the time limit. Among these methods, \texttt{R1} consistently produces near-global solutions for all instances under each nomination. Also, we observe that there are cases where both \texttt{R2} and \texttt{R3} converge to a point of infeasibility in the local solver whereas \texttt{R1} provides high-quality warm-starting points for all instances.

We summarize our computational results as follows: \texttt{BARON} is able to provide consistent dual and primal bounds for GasLib-11 and GasLib-24; however, we encounter some persistent inconsistency and infeasibility issues with \texttt{BARON} as the gas network gets more complex. Under the objective functions $g_1$ and $g_3$, we see that the performance of the convex relaxations in terms of the optimality gap and total time is better than \texttt{BARON} for all instances. In particular, \texttt{BARON} returns considerably larger gaps as the stress level of the network increases. Moreover, \texttt{R1} consistently provides tighter dual bounds for all instances, and it also offers high-quality warm-starting points for the local solver Ipopt compared to \texttt{R2} and \texttt{R3}. We also observe that there are cases where both \texttt{R2} and \texttt{R3} converge to a point of local infeasibility when initiated with the starting solutions produced by these relaxations. We conclude that our relaxation method \texttt{R1} consistently leads to numerically more stable results in comparison with those methods from the literature and the global~solver.

Finally, we present the performance ratios in terms of the best optimality gap of different methods across three objectives and overall performance.  Whenever two or more methods have the same optimality gap for an instance, we distribute the ratios equally. Let us start by the GasLib-24 instance with known gas (positive) flow directions. Under $g_1$ and $g_3$, both \texttt{R1} and \texttt{R3} have the same performance ratio (50\%) whereas \texttt{BARON} outperforms these methods under $g_2$ for each instance. For the GasLib-24 instance, the overall performance for each of these three methods is equal (33.33\%). For the remaining GasLib instances, we calculate the performance ratio of each method for 24, 12 and 24 instances under the objective functions $g_1$, $g_2$ and $g_3$, respectively. We present them along with the overall performance out of all instances in Table~\ref{table:summaryPerformance}. Our proposed method \texttt{R1} consistently achieves the highest performance ratio, with 52\%, 40\% and 66\% under $g_1$, $g_2$ and $g_3$, respectively. Moreover, \texttt{R1} outperforms the relaxation methods \texttt{R2} and \texttt{R3} for all cases. Although our method \texttt{R1} also outperforms \texttt{BARON}, its performance becomes comparable and improves considerably under $g_2$ (35\%). However, we would like to note that the overall performance of \texttt{BARON} still remains low (8\%). \texttt{R2} and \texttt{R3} yield comparable and moderate results, with a total performance of 17\% and 20\%, respectively. Still, our proposed method \texttt{R1} shows the highest overall performance. These results show that \texttt{R1} not only provides tighter dual bounds but also provides numerically more stable and reliable solutions than both existing relaxations and the global solver \texttt{BARON}.

\begin{table}[h!] 
\centering
\scalebox{1.2}{
\begin{tabular}{ccccccccc}
\hline
      &  & \multirow{2}{*}{\begin{tabular}[c]{@{}c@{}}$g_1$\\ (out of 24)\end{tabular}} &  & \multirow{2}{*}{\begin{tabular}[c]{@{}c@{}}$g_2$\\ (out of 12)\end{tabular}} &  & \multirow{2}{*}{\begin{tabular}[c]{@{}c@{}}$g_3$\\ (out of 24)\end{tabular}} &  & \multirow{2}{*}{\begin{tabular}[c]{@{}c@{}}Total\\ (out of  60)\end{tabular}} \\
      &  &                                                                           &  &                                                                           &  &                                                                           &  &                                                                               \\ \hline
\texttt{R1}    &  & \bf 52\%                                                                        &  & \bf  40\%                                                                         &  &  \bf  66\%                                                                         &  &  \bf 55\%                                                                             \\
\texttt{R2}    &  & 18\%                                                                         &  & 15\%                                                                         &  & 18\%                                                                         &  & 17\%                                                                             \\
\texttt{R3}    &  & 29\%                                                                         &  & 10\%                                                                         &  & 15\%                                                                         &  & 20\%                                                                             \\
\texttt{BARON} &  & 1\%                                                                          &  & 35\%                                                                         &  & 1\%                                                                          &  & 8\%                                                                              \\ \hline
\end{tabular}
}
\caption{The summary of the performance ratios of different methods for GasLib-11, 39, 40 and 134 instances.}
\label{table:summaryPerformance}
\end{table}

\section{Conclusion} \label{section:conclusion}
In this paper, we studied the gas transportation and storage optimization problem, in which we explicitly considered the highly nonlinear and nonconvex inherent aspects due to gas physics as well as the discrete aspects due to the control decisions of active network elements. We studied the nonconvex sets induced by gas physics in detail and proposed MISOCP relaxations for this problem. In particular, our relaxation approach is based on the convex hull representations of the nonconvex sets as follows: We derived the exact convex hull formulation for the feasible (nonconvex) set for pipes and show that it is SOCr. We also provided the convex hull of the extreme points of the feasible (nonconvex) set for compressors and show that it is POWr. To fully utilize the state-of-the-art solver, we also proposed a second-order cone outer-approximation for this POWr set. Based on these relaxations, we also proposed an algorithmic framework to solve the problem. We tested our framework on multiple GasLib instances under different settings. Our computational results show that our method \texttt{R1} outperforms the state-of-the-art global solver and other relaxation methods from the literature. Our approach consistently provides tight dual bounds and produces (near) global solutions as well as high-quality warm-starting points for local solvers ensuring numerical consistency and scalability.

There are several promising directions for future research. Extending our approach to gas integrated energy systems for electricity and gas as well as incorporating uncertainty in demand and supply conditions could be an important area of study. Additionally, decomposition approaches for multi-period gas network optimization problems could improve computational efficiency and scalability.



\section*{Declarations}

\subsection*{Ethics approval and consent to participate}
Not applicable

\subsection*{Consent for publication}
Not applicable

\subsection*{Funding}

The Scientific and Technological Research Council of Turkey  (project number: 119M855).

\subsection*{Availability of data and materials}
Data is provided within the manuscript and references therein.

\subsection*{Competing interests}
The authors declare that they have no competing interests.

\subsection*{Authors' contributions}
B.C.O. worked on the methodology, implementation, visualization, analysis, writing and B.K. worked on the methodology, implementation, analysis, writing.

\subsection*{Acknowledgements}


This work was supported by the  Scientific and Technological Research Council of Turkey under grant number 119M855.
\bibliography{sn-bibliography}


\begin{thebibliography}{46}
\ifx \bisbn   \undefined \def \bisbn  #1{ISBN #1}\fi
\ifx \binits  \undefined \def \binits#1{#1}\fi
\ifx \bauthor  \undefined \def \bauthor#1{#1}\fi
\ifx \batitle  \undefined \def \batitle#1{#1}\fi
\ifx \bjtitle  \undefined \def \bjtitle#1{#1}\fi
\ifx \bvolume  \undefined \def \bvolume#1{\textbf{#1}}\fi
\ifx \byear  \undefined \def \byear#1{#1}\fi
\ifx \bissue  \undefined \def \bissue#1{#1}\fi
\ifx \bfpage  \undefined \def \bfpage#1{#1}\fi
\ifx \blpage  \undefined \def \blpage #1{#1}\fi
\ifx \burl  \undefined \def \burl#1{\textsf{#1}}\fi
\ifx \doiurl  \undefined \def \doiurl#1{\url{https://doi.org/#1}}\fi
\ifx \betal  \undefined \def \betal{\textit{et al.}}\fi
\ifx \binstitute  \undefined \def \binstitute#1{#1}\fi
\ifx \binstitutionaled  \undefined \def \binstitutionaled#1{#1}\fi
\ifx \bctitle  \undefined \def \bctitle#1{#1}\fi
\ifx \beditor  \undefined \def \beditor#1{#1}\fi
\ifx \bpublisher  \undefined \def \bpublisher#1{#1}\fi
\ifx \bbtitle  \undefined \def \bbtitle#1{#1}\fi
\ifx \bedition  \undefined \def \bedition#1{#1}\fi
\ifx \bseriesno  \undefined \def \bseriesno#1{#1}\fi
\ifx \blocation  \undefined \def \blocation#1{#1}\fi
\ifx \bsertitle  \undefined \def \bsertitle#1{#1}\fi
\ifx \bsnm \undefined \def \bsnm#1{#1}\fi
\ifx \bsuffix \undefined \def \bsuffix#1{#1}\fi
\ifx \bparticle \undefined \def \bparticle#1{#1}\fi
\ifx \barticle \undefined \def \barticle#1{#1}\fi
\bibcommenthead
\ifx \bconfdate \undefined \def \bconfdate #1{#1}\fi
\ifx \botherref \undefined \def \botherref #1{#1}\fi
\ifx \url \undefined \def \url#1{\textsf{#1}}\fi
\ifx \bchapter \undefined \def \bchapter#1{#1}\fi
\ifx \bbook \undefined \def \bbook#1{#1}\fi
\ifx \bcomment \undefined \def \bcomment#1{#1}\fi
\ifx \oauthor \undefined \def \oauthor#1{#1}\fi
\ifx \citeauthoryear \undefined \def \citeauthoryear#1{#1}\fi
\ifx \endbibitem  \undefined \def \endbibitem {}\fi
\ifx \bconflocation  \undefined \def \bconflocation#1{#1}\fi
\ifx \arxivurl  \undefined \def \arxivurl#1{\textsf{#1}}\fi
\csname PreBibitemsHook\endcsname

\bibitem[\protect\citeauthoryear{{IEA}}{2023}]{IEA2023}
\begin{botherref}
\oauthor{\bsnm{{IEA}}}:
{Global Gas Security Review 2023}.
\url{https://www.iea.org/reports/global-gas-security-review-2023}
(2023)
\end{botherref}
\endbibitem

\bibitem[\protect\citeauthoryear{R{\'\i}os-Mercado and Borraz-S{\'a}nchez}{2015}]{Rios2015}
\begin{barticle}
\bauthor{\bsnm{R{\'\i}os-Mercado}, \binits{R.Z.}},
\bauthor{\bsnm{Borraz-S{\'a}nchez}, \binits{C.}}:
\batitle{Optimization problems in natural gas transportation systems: A state-of-the-art review}.
\bjtitle{Applied Energy}
\bvolume{147},
\bfpage{536}--\blpage{555}
(\byear{2015})
\end{barticle}
\endbibitem

\bibitem[\protect\citeauthoryear{Gross et~al.}{2019}]{Gross2019}
\begin{barticle}
\bauthor{\bsnm{Gross}, \binits{M.}},
\bauthor{\bsnm{Pfetsch}, \binits{M.E.}},
\bauthor{\bsnm{Schewe}, \binits{L.}},
\bauthor{\bsnm{Schmidt}, \binits{M.}},
\bauthor{\bsnm{Skutella}, \binits{M.}}:
\batitle{Algorithmic results for potential-based flows: Easy and hard cases}.
\bjtitle{Networks}
\bvolume{73}(\bissue{3}),
\bfpage{306}--\blpage{324}
(\byear{2019})
\end{barticle}
\endbibitem

\bibitem[\protect\citeauthoryear{Labb{\'e} et~al.}{2021}]{labbe2021polynomial}
\begin{barticle}
\bauthor{\bsnm{Labb{\'e}}, \binits{M.}},
\bauthor{\bsnm{Plein}, \binits{F.}},
\bauthor{\bsnm{Schmidt}, \binits{M.}},
\bauthor{\bsnm{Th{\"u}rauf}, \binits{J.}}:
\batitle{{Deciding feasibility of a booking in the European gas market on a cycle is in P for the case of passive networks}}.
\bjtitle{Networks}
\bvolume{78}(\bissue{2}),
\bfpage{128}--\blpage{152}
(\byear{2021})
\end{barticle}
\endbibitem

\bibitem[\protect\citeauthoryear{Humpola}{2014}]{Humpola2014}
\begin{botherref}
\oauthor{\bsnm{Humpola}, \binits{J.}}:
{Gas Network Optimization by MINLP}.
PhD thesis,
Technische Universität Berlin
(2014)
\end{botherref}
\endbibitem

\bibitem[\protect\citeauthoryear{de~Wolf and Smeers}{2000}]{Wolf2000}
\begin{barticle}
\bauthor{\bsnm{Wolf}, \binits{D.}},
\bauthor{\bsnm{Smeers}, \binits{Y.}}:
\batitle{The gas transmission problem solved by an extension of the simplex algorithm}.
\bjtitle{Management Science}
\bvolume{46}(\bissue{11}),
\bfpage{1454}--\blpage{1465}
(\byear{2000})
\end{barticle}
\endbibitem

\bibitem[\protect\citeauthoryear{Ehrhardt and Steinbach}{2005}]{Ehrhardt2005}
\begin{bchapter}
\bauthor{\bsnm{Ehrhardt}, \binits{K.}},
\bauthor{\bsnm{Steinbach}, \binits{M.C.}}:
\bctitle{Nonlinear optimization in gas networks}.
In: \beditor{\bsnm{Bock}, \binits{H.G.}},
\beditor{\bsnm{Phu}, \binits{H.X.}},
\beditor{\bsnm{Kostina}, \binits{E.}},
\beditor{\bsnm{Rannacher}, \binits{R.}} (eds.)
\bbtitle{Modeling, Simulation and Optimization of Complex Processes},
pp. \bfpage{139}--\blpage{148}.
\bpublisher{Springer},
\blocation{Berlin, Heidelberg}
(\byear{2005})
\end{bchapter}
\endbibitem

\bibitem[\protect\citeauthoryear{Steinbach}{2007}]{Steinbach2007}
\begin{barticle}
\bauthor{\bsnm{Steinbach}, \binits{M.C.}}:
\batitle{{On PDE solution in transient optimization of gas networks}}.
\bjtitle{Journal of Computational and Applied Mathematics}
\bvolume{203}(\bissue{2}),
\bfpage{345}--\blpage{361}
(\byear{2007}).
\bcomment{Special Issue: The first Indo-German Conference on PDE, Scientific Computing and Optimization in Applications}
\end{barticle}
\endbibitem

\bibitem[\protect\citeauthoryear{Schmidt}{2015}]{Schmidt2015}
\begin{barticle}
\bauthor{\bsnm{Schmidt}, \binits{M.}}:
\batitle{An interior-point method for nonlinear optimization problems with locatable and separable nonsmoothness}.
\bjtitle{EURO Journal on Computational Optimization}
\bvolume{3}(\bissue{4}),
\bfpage{309}--\blpage{348}
(\byear{2015})
\end{barticle}
\endbibitem

\bibitem[\protect\citeauthoryear{Wu et~al.}{2007}]{wu2007gas}
\begin{barticle}
\bauthor{\bsnm{Wu}, \binits{Y.}},
\bauthor{\bsnm{Lai}, \binits{K.K.}},
\bauthor{\bsnm{Liu}, \binits{Y.}}:
\batitle{Deterministic global optimization approach to steady-state distribution gas pipeline networks}.
\bjtitle{Optimization and Engineering}
\bvolume{8}(\bissue{3}),
\bfpage{259}--\blpage{275}
(\byear{2007})
\end{barticle}
\endbibitem

\bibitem[\protect\citeauthoryear{Schmidt}{2013}]{Schmidt2013}
\begin{botherref}
\oauthor{\bsnm{Schmidt}, \binits{M.}}:
{A Generic Interior-Point Framework for Nonsmooth and Complementarity Constrained Nonlinear Optimization}.
PhD thesis,
Gottfried Wilhelm Leibniz Universität Hannover
(2013)
\end{botherref}
\endbibitem

\bibitem[\protect\citeauthoryear{Pfetsch et~al.}{2015}]{Pfetsch2015}
\begin{barticle}
\bauthor{\bsnm{Pfetsch}, \binits{M.E.}},
\bauthor{\bsnm{F{\"u}genschuh}, \binits{A.}},
\bauthor{\bsnm{Gei{\ss}ler}, \binits{B.}},
\bauthor{\bsnm{Gei{\ss}ler}, \binits{N.}},
\bauthor{\bsnm{Gollmer}, \binits{R.}},
\bauthor{\bsnm{Hiller}, \binits{B.}},
\bauthor{\bsnm{Humpola}, \binits{J.}},
\bauthor{\bsnm{Koch}, \binits{T.}},
\bauthor{\bsnm{Lehmann}, \binits{T.}},
\bauthor{\bsnm{Martin}, \binits{A.}},
\bauthor{\bsnm{Morsi}, \binits{A.}},
\bauthor{\bsnm{R{\"o}vekamp}, \binits{J.}},
\bauthor{\bsnm{Schewe}, \binits{L.}},
\bauthor{\bsnm{Schmidt}, \binits{M.}},
\bauthor{\bsnm{Schultz}, \binits{R.}},
\bauthor{\bsnm{Schwarz}, \binits{R.}},
\bauthor{\bsnm{Schweiger}, \binits{J.}},
\bauthor{\bsnm{Stangl}, \binits{C.}},
\bauthor{\bsnm{Steinbach}, \binits{M.C.}},
\bauthor{\bsnm{Vigerske}, \binits{S.}},
\bauthor{\bsnm{Willert}, \binits{B.M.}}:
\batitle{Validation of nominations in gas network optimization: models, methods, and solutions}.
\bjtitle{Optimization Methods and Software}
\bvolume{30}(\bissue{1}),
\bfpage{15}--\blpage{53}
(\byear{2015})
\end{barticle}
\endbibitem

\bibitem[\protect\citeauthoryear{Zhang and Zhu}{1996}]{Zhang1996}
\begin{barticle}
\bauthor{\bsnm{Zhang}, \binits{J.}},
\bauthor{\bsnm{Zhu}, \binits{D.}}:
\batitle{A bilevel programming method for pipe network optimization}.
\bjtitle{SIAM Journal on Optimization}
\bvolume{6}(\bissue{3}),
\bfpage{838}--\blpage{857}
(\byear{1996})
\end{barticle}
\endbibitem

\bibitem[\protect\citeauthoryear{Wu et~al.}{2000}]{Wu2000}
\begin{barticle}
\bauthor{\bsnm{Wu}, \binits{S.}},
\bauthor{\bsnm{R{\'\i}os-Mercado}, \binits{R.Z.}},
\bauthor{\bsnm{Boyd}, \binits{E.A.}},
\bauthor{\bsnm{Scott}, \binits{L.R.}}:
\batitle{Model relaxations for the fuel cost minimization of steady-state gas pipeline networks}.
\bjtitle{Mathematical and Computer Modelling}
\bvolume{31}(\bissue{2}),
\bfpage{197}--\blpage{220}
(\byear{2000})
\end{barticle}
\endbibitem

\bibitem[\protect\citeauthoryear{Andre et~al.}{2009}]{Andre2009}
\begin{barticle}
\bauthor{\bsnm{Andre}, \binits{J.}},
\bauthor{\bsnm{Bonnans}, \binits{F.}},
\bauthor{\bsnm{Cornibert}, \binits{L.}}:
\batitle{Optimization of capacity expansion planning for gas transportation networks}.
\bjtitle{European Journal of Operational Research}
\bvolume{197}(\bissue{3}),
\bfpage{1019}--\blpage{1027}
(\byear{2009})
\end{barticle}
\endbibitem

\bibitem[\protect\citeauthoryear{Babonneau et~al.}{2012}]{Babonneau2012}
\begin{barticle}
\bauthor{\bsnm{Babonneau}, \binits{F.}},
\bauthor{\bsnm{Nesterov}, \binits{Y.}},
\bauthor{\bsnm{Vial}, \binits{J.-P.}}:
\batitle{{Design and Operations of Gas Transmission Networks}}.
\bjtitle{Operations Research}
\bvolume{60}(\bissue{1}),
\bfpage{34}--\blpage{47}
(\byear{2012})
\end{barticle}
\endbibitem

\bibitem[\protect\citeauthoryear{Fügenschuh and Humpola}{2013}]{Fugenschuh2013}
\begin{botherref}
\oauthor{\bsnm{Fügenschuh}, \binits{A.}},
\oauthor{\bsnm{Humpola}, \binits{J.}}:
{A Unified View on Relaxations for a Nonlinear Network Flow Problem}.
Technical Report 13-31,
ZIB,
Takustr. 7, 14195 Berlin
(2013)
\end{botherref}
\endbibitem

\bibitem[\protect\citeauthoryear{Wu et~al.}{2017}]{Wu2017Mono}
\begin{bchapter}
\bauthor{\bsnm{Wu}, \binits{F.}},
\bauthor{\bsnm{Nagarajan}, \binits{H.}},
\bauthor{\bsnm{Zlotnik}, \binits{A.}},
\bauthor{\bsnm{Sioshansi}, \binits{R.}},
\bauthor{\bsnm{Rudkevich}, \binits{A.M.}}:
\bctitle{Adaptive convex relaxations for gas pipeline network optimization}.
In: \bbtitle{2017 American Control Conference (ACC)},
pp. \bfpage{4710}--\blpage{4716}
(\byear{2017})
\end{bchapter}
\endbibitem

\bibitem[\protect\citeauthoryear{Ordoudis et~al.}{2019}]{Ordoudis2019}
\begin{barticle}
\bauthor{\bsnm{Ordoudis}, \binits{C.}},
\bauthor{\bsnm{Pinson}, \binits{P.}},
\bauthor{\bsnm{Morales}, \binits{J.M.}}:
\batitle{An integrated market for electricity and natural gas systems with stochastic power producers}.
\bjtitle{European Journal of Operational Research}
\bvolume{272}(\bissue{2}),
\bfpage{642}--\blpage{654}
(\byear{2019})
\end{barticle}
\endbibitem

\bibitem[\protect\citeauthoryear{Borraz-S{\'a}nchez et~al.}{2016}]{Borraz2016SOCP}
\begin{barticle}
\bauthor{\bsnm{Borraz-S{\'a}nchez}, \binits{C.}},
\bauthor{\bsnm{Bent}, \binits{R.}},
\bauthor{\bsnm{Backhaus}, \binits{S.}},
\bauthor{\bsnm{Hijazi}, \binits{H.}},
\bauthor{\bsnm{Hentenryck}, \binits{P.V.}}:
\batitle{Convex relaxations for gas expansion planning}.
\bjtitle{INFORMS Journal on Computing}
\bvolume{28}(\bissue{4}),
\bfpage{645}--\blpage{656}
(\byear{2016})
\end{barticle}
\endbibitem

\bibitem[\protect\citeauthoryear{Martin et~al.}{2006}]{Martin2006}
\begin{barticle}
\bauthor{\bsnm{Martin}, \binits{A.}},
\bauthor{\bsnm{Möller}, \binits{M.}},
\bauthor{\bsnm{Moritz}, \binits{S.}}:
\batitle{Mixed integer models for the stationary case of gas network optimization}.
\bjtitle{Mathematical Programming}
\bvolume{105}(\bissue{2-3}),
\bfpage{563}
(\byear{2006})
\end{barticle}
\endbibitem

\bibitem[\protect\citeauthoryear{Zheng et~al.}{2010}]{Zheng2010}
\begin{bchapter}
\bauthor{\bsnm{Zheng}, \binits{Q.P.}},
\bauthor{\bsnm{Rebennack}, \binits{S.}},
\bauthor{\bsnm{Iliadis}, \binits{N.A.}},
\bauthor{\bsnm{Pardalos}, \binits{P.M.}}:
\bctitle{{Optimization Models in the Natural Gas Industry}}.
In: \beditor{\bsnm{Pardalos}, \binits{P.}},
\beditor{\bsnm{Rebennack}, \binits{S.}},
\beditor{\bsnm{Pereira}, \binits{M.}},
\beditor{\bsnm{Iliadis}, \binits{N.}} (eds.)
\bbtitle{Handbook of Power Systems I. Energy Systems}.
\bpublisher{{Springer}},
\blocation{{Berlin, Heidelberg}}
(\byear{2010})
\end{bchapter}
\endbibitem

\bibitem[\protect\citeauthoryear{Wang et~al.}{2018}]{Wang2018}
\begin{barticle}
\bauthor{\bsnm{Wang}, \binits{B.}},
\bauthor{\bsnm{Yuan}, \binits{M.}},
\bauthor{\bsnm{Zhang}, \binits{H.}},
\bauthor{\bsnm{Zhao}, \binits{W.}},
\bauthor{\bsnm{Liang}, \binits{Y.}}:
\batitle{{An MILP model for optimal design of multi-period natural gas transmission network}}.
\bjtitle{Chemical Engineering Research and Design}
\bvolume{129},
\bfpage{122}--\blpage{131}
(\byear{2018})
\end{barticle}
\endbibitem

\bibitem[\protect\citeauthoryear{Ojha et~al.}{2017}]{Ohja2017SDP}
\begin{bchapter}
\bauthor{\bsnm{Ojha}, \binits{A.}},
\bauthor{\bsnm{Kekatos}, \binits{V.}},
\bauthor{\bsnm{Baldick}, \binits{R.}}:
\bctitle{Solving the natural gas flow problem using semidefinite program relaxation}.
In: \bbtitle{2017 IEEE Power Energy Society General Meeting},
pp. \bfpage{1}--\blpage{5}
(\byear{2017})
\end{bchapter}
\endbibitem

\bibitem[\protect\citeauthoryear{Borraz-S{\'a}nchez et~al.}{2016}]{BorrazSanchez2016-2}
\begin{bchapter}
\bauthor{\bsnm{Borraz-S{\'a}nchez}, \binits{C.}},
\bauthor{\bsnm{Bent}, \binits{R.}},
\bauthor{\bsnm{Backhaus}, \binits{S.}},
\bauthor{\bsnm{Blumsack}, \binits{S.}},
\bauthor{\bsnm{Hijazi}, \binits{H.}},
\bauthor{\bsnm{Hentenryck}, \binits{P.}}:
\bctitle{Convex optimization for joint expansion planning of natural gas and power systems}.
In: \bbtitle{2016 49th Hawaii International Conference on System Sciences (HICSS)},
pp. \bfpage{2536}--\blpage{2545}
(\byear{2016})
\end{bchapter}
\endbibitem

\bibitem[\protect\citeauthoryear{Wen et~al.}{2018}]{Wen2018}
\begin{barticle}
\bauthor{\bsnm{Wen}, \binits{Y.}},
\bauthor{\bsnm{Qu}, \binits{X.}},
\bauthor{\bsnm{Li}, \binits{W.}},
\bauthor{\bsnm{Liu}, \binits{X.}},
\bauthor{\bsnm{Ye}, \binits{X.}}:
\batitle{Synergistic operation of electricity and natural gas networks via admm}.
\bjtitle{IEEE Transactions on Smart Grid}
\bvolume{9}(\bissue{5}),
\bfpage{4555}--\blpage{4565}
(\byear{2018})
\end{barticle}
\endbibitem

\bibitem[\protect\citeauthoryear{He et~al.}{2018}]{He2018}
\begin{barticle}
\bauthor{\bsnm{He}, \binits{Y.}},
\bauthor{\bsnm{Shahidehpour}, \binits{M.}},
\bauthor{\bsnm{Li}, \binits{Z.}},
\bauthor{\bsnm{Guo}, \binits{C.}},
\bauthor{\bsnm{Zhu}, \binits{B.}}:
\batitle{Robust constrained operation of integrated electricity-natural gas system considering distributed natural gas storage}.
\bjtitle{IEEE Transactions on Sustainable Energy}
\bvolume{9}(\bissue{3}),
\bfpage{1061}--\blpage{1071}
(\byear{2018})
\end{barticle}
\endbibitem

\bibitem[\protect\citeauthoryear{Singh and Kekatos}{2019}]{Singh2019}
\begin{bchapter}
\bauthor{\bsnm{Singh}, \binits{M.K.}},
\bauthor{\bsnm{Kekatos}, \binits{V.}}:
\bctitle{{Natural Gas Flow Equations: Uniqueness and an MI-SOCP Solver}}.
In: \bbtitle{2019 American Control Conference (ACC)},
pp. \bfpage{2114}--\blpage{2120}
(\byear{2019})
\end{bchapter}
\endbibitem

\bibitem[\protect\citeauthoryear{Schwele et~al.}{2019}]{Schwele2019SOCP}
\begin{bchapter}
\bauthor{\bsnm{Schwele}, \binits{A.}},
\bauthor{\bsnm{Ordoudis}, \binits{C.}},
\bauthor{\bsnm{Kazempour}, \binits{J.}},
\bauthor{\bsnm{Pinson}, \binits{P.}}:
\bctitle{Coordination of power and natural gas systems: Convexification approaches for linepack modeling}.
In: \bbtitle{2019 IEEE Milan PowerTech},
pp. \bfpage{1}--\blpage{6}
(\byear{2019})
\end{bchapter}
\endbibitem

\bibitem[\protect\citeauthoryear{Li et~al.}{2024}]{li2024misocp}
\begin{barticle}
\bauthor{\bsnm{Li}, \binits{Y.}},
\bauthor{\bsnm{Dey}, \binits{S.S.}},
\bauthor{\bsnm{Sahinidis}, \binits{N.V.}}:
\batitle{A reformulation-enumeration minlp algorithm for gas network design}.
\bjtitle{Journal of Global Optimization}
\bvolume{90}(\bissue{4}),
\bfpage{931}--\blpage{963}
(\byear{2024})
\end{barticle}
\endbibitem

\bibitem[\protect\citeauthoryear{Burlacu et~al.}{2019}]{burlacu2019}
\begin{barticle}
\bauthor{\bsnm{Burlacu}, \binits{R.}},
\bauthor{\bsnm{Egger}, \binits{H.}},
\bauthor{\bsnm{Gro{\ss}}, \binits{M.}},
\bauthor{\bsnm{Martin}, \binits{A.}},
\bauthor{\bsnm{Pfetsch}, \binits{M.E.}},
\bauthor{\bsnm{Schewe}, \binits{L.}},
\bauthor{\bsnm{Sirvent}, \binits{M.}},
\bauthor{\bsnm{Skutella}, \binits{M.}}:
\batitle{Maximizing the storage capacity of gas networks: a global minlp approach}.
\bjtitle{Optimization and Engineering}
\bvolume{20}(\bissue{2}),
\bfpage{543}--\blpage{573}
(\byear{2019})
\end{barticle}
\endbibitem

\bibitem[\protect\citeauthoryear{Borraz-S{\'a}nchez and Ríos-Mercado}{2009}]{Sanchez2009Tabu}
\begin{barticle}
\bauthor{\bsnm{Borraz-S{\'a}nchez}, \binits{C.}},
\bauthor{\bsnm{Ríos-Mercado}, \binits{R.Z.}}:
\batitle{Improving the operation of pipeline systems on cyclic structures by tabu search}.
\bjtitle{Computers \& Chemical Engineering}
\bvolume{33}(\bissue{1}),
\bfpage{58}--\blpage{64}
(\byear{2009})
\end{barticle}
\endbibitem

\bibitem[\protect\citeauthoryear{Koch et~al.}{2015}]{Koch2015}
\begin{bbook}
\bauthor{\bsnm{Koch}, \binits{T.}},
\bauthor{\bsnm{Pfetsch}, \binits{M.E.}},
\bauthor{\bsnm{Schewe}, \binits{L.}}:
\bbtitle{{Evaluating Gas Network Capacities}}.
\bpublisher{{MOS-SIAM Series on Optimization}},
\blocation{SIAM, Philadelphia}
(\byear{2015})
\end{bbook}
\endbibitem

\bibitem[\protect\citeauthoryear{Osiadacz}{1987}]{Osiadacz1987}
\begin{bbook}
\bauthor{\bsnm{Osiadacz}, \binits{A.}}:
\bbtitle{Simulation and Analysis of Gas Networks}.
\bpublisher{Gulf Publishing Company},
\blocation{Houston, TX}
(\byear{1987})
\end{bbook}
\endbibitem

\bibitem[\protect\citeauthoryear{Tuncer and Kocuk}{2023}]{Tuncer2023}
\begin{barticle}
\bauthor{\bsnm{Tuncer}, \binits{D.}},
\bauthor{\bsnm{Kocuk}, \binits{B.}}:
\batitle{{An MISOCP-Based Decomposition Approach for the Unit Commitment Problem With AC Power Flows}}.
\bjtitle{IEEE Transactions on Power Systems}
\bvolume{38}(\bissue{4}),
\bfpage{3388}--\blpage{3400}
(\byear{2023})
\end{barticle}
\endbibitem

\bibitem[\protect\citeauthoryear{Okumuşoğlu et~al.}{2024}]{Okumusoglu2024}
\begin{barticle}
\bauthor{\bsnm{Okumuşoğlu}, \binits{B.C.}},
\bauthor{\bsnm{Basciftci}, \binits{B.}},
\bauthor{\bsnm{Kocuk}, \binits{B.}}:
\batitle{{An Integrated Predictive Maintenance and Operations Scheduling Framework for Power Systems Under Failure Uncertainty}}.
\bjtitle{INFORMS Journal on Computing}
\bvolume{36}(\bissue{5}),
\bfpage{1335}--\blpage{1358}
(\byear{2024})
\doiurl{10.1287/ijoc.2022.0154}
\end{barticle}
\endbibitem

\bibitem[\protect\citeauthoryear{McCormick}{1976}]{McCormick1976}
\begin{barticle}
\bauthor{\bsnm{McCormick}, \binits{G.P.}}:
\batitle{Computability of global solutions to factorable nonconvex programs: Part i ---convex underestimating problems}.
\bjtitle{Mathematical Programming}
\bvolume{10}(\bissue{1}),
\bfpage{147}--\blpage{175}
(\byear{1976})
\end{barticle}
\endbibitem

\bibitem[\protect\citeauthoryear{Tawarmalani and Richard}{2013}]{tawarmalani2013decomposition}
\begin{botherref}
\oauthor{\bsnm{Tawarmalani}, \binits{M.}},
\oauthor{\bsnm{Richard}, \binits{J.-P.P.}}:
Decomposition techniques in convexification of inequalities.
Technical report,
Technical report
(2013)
\end{botherref}
\endbibitem

\bibitem[\protect\citeauthoryear{Hillestad and Jacobsen}{1980}]{hillestad1980linear}
\begin{barticle}
\bauthor{\bsnm{Hillestad}, \binits{R.J.}},
\bauthor{\bsnm{Jacobsen}, \binits{S.E.}}:
\batitle{Linear programs with an additional reverse convex constraint}.
\bjtitle{Applied Mathematics and Optimization}
\bvolume{6},
\bfpage{257}--\blpage{269}
(\byear{1980})
\end{barticle}
\endbibitem

\bibitem[\protect\citeauthoryear{Ben-Tal and Nemirovski}{2001}]{ben2001lectures}
\begin{bbook}
\bauthor{\bsnm{Ben-Tal}, \binits{A.}},
\bauthor{\bsnm{Nemirovski}, \binits{A.}}:
\bbtitle{Lectures on Modern Convex Optimization: Analysis, Algorithms, and Engineering Applications}.
\bpublisher{{SIAM-MOS Series on Optimization}},
\blocation{SIAM, Philadelphia}
(\byear{2001})
\end{bbook}
\endbibitem

\bibitem[\protect\citeauthoryear{Matthias}{2024}]{Mittelmann}
\begin{botherref}
\oauthor{\bsnm{Matthias}, \binits{M.}}:
{mittelmann-plots - Visualizations of Mittelmann benchmarks}.
\url{https://mattmilten.github.io/mittelmann-plots/}
(2024)
\end{botherref}
\endbibitem

\bibitem[\protect\citeauthoryear{Schmidt et~al.}{2017}]{Gaslib}
\begin{botherref}
\oauthor{\bsnm{Schmidt}, \binits{M.}},
\oauthor{\bsnm{Aßmann}, \binits{D.}},
\oauthor{\bsnm{Burlacu}, \binits{R.}},
\oauthor{\bsnm{Humpola}, \binits{J.}},
\oauthor{\bsnm{Joormann}, \binits{I.}},
\oauthor{\bsnm{Kanelakis}, \binits{N.}},
\oauthor{\bsnm{Koch}, \binits{T.}},
\oauthor{\bsnm{Oucherif}, \binits{D.}},
\oauthor{\bsnm{Pfetsch}, \binits{M.E.}},
\oauthor{\bsnm{Schewe}, \binits{L.}},
\oauthor{\bsnm{Schwarz}, \binits{R.}},
\oauthor{\bsnm{Sirvent}, \binits{M.}}:
Gaslib—a library of gas network instances.
Data
\textbf{2}(4)
(2017)
\doiurl{10.3390/data2040040}
\end{botherref}
\endbibitem

\bibitem[\protect\citeauthoryear{Hennings}{2023}]{Hennings2023}
\begin{botherref}
\oauthor{\bsnm{Hennings}, \binits{F.}}:
{Modeling and solving real-world transient gas network transport problems using mathematical programming}.
PhD thesis,
Technische Universität Berlin
(2023)
\end{botherref}
\endbibitem

\bibitem[\protect\citeauthoryear{Gurobi~Optimization}{2024}]{Gurobi}
\begin{botherref}
\oauthor{\bsnm{Gurobi~Optimization}, \binits{L.}}:
{Gurobi Optimizer Reference Manual}.
\url{www.gurobi.com/documentation/}
(2024)
\end{botherref}
\endbibitem

\bibitem[\protect\citeauthoryear{W{\"a}chter and Biegler}{2006}]{ipopt}
\begin{barticle}
\bauthor{\bsnm{W{\"a}chter}, \binits{A.}},
\bauthor{\bsnm{Biegler}, \binits{L.T.}}:
\batitle{On the implementation of an interior-point filter line-search algorithm for large-scale nonlinear programming}.
\bjtitle{Mathematical Programming}
\bvolume{106}(\bissue{1}),
\bfpage{25}--\blpage{57}
(\byear{2006})
\end{barticle}
\endbibitem

\bibitem[\protect\citeauthoryear{Tawarmalani and Sahinidis}{2005}]{baron}
\begin{barticle}
\bauthor{\bsnm{Tawarmalani}, \binits{M.}},
\bauthor{\bsnm{Sahinidis}, \binits{N.V.}}:
\batitle{A polyhedral branch-and-cut approach to global optimization}.
\bjtitle{Mathematical Programming}
\bvolume{103}(\bissue{2}),
\bfpage{225}--\blpage{249}
(\byear{2005})
\end{barticle}
\endbibitem

\end{thebibliography}

\begin{appendices}

\section{Extreme Points of the Weymouth Set}  \label{app:weymouth}

For simplicity of notation, consider the set
\(
\mathcal{X}^+ = \{ (x,y,z) \in \mathcal{B}^+ : x-y = w z^2   \},
\)
where
$\mathcal{B}^+ = \{ (x,y,z)\in  [\munderbar{x}, \bar{x}]\times [\munderbar{y}, \bar{y}]\times[0, {\alpha}] : 0 \le x-y \le {\beta} \} $ is a polytope such that ${\alpha}\ge 0$ and ${\beta} \ge 0$. 
We have already established in Proposition~\ref{prop:convexHull_pipes_Pos} that $\textup{conv}(\mathcal{X}^+) = \{ (x,y,z)\in \mathcal{\check B}^+ : x-y \ge w z^2 \}$ for some polytope $\mathcal{\check B}^+$. 
In order to construct  $\mathcal{\check B}^+$, we follow 
  \cite{hillestad1980linear}: We first find the extreme points of the intersection of each edge of the polytope $\mathcal{B}^+$ with the reverse convex constraint $ x-y \le  w z^2$ (Algorithm~\ref{alg:weymouthExtremesPlus} provides the set of these points denoted as $\mathcal{Q^+}$). Then, we obtain   $ \mathcal{\check B^+} = \textup{conv}(\mathcal{Q^+})$.
  
\begin{algorithm}[h!]
\caption{Extreme points of $\mathcal{\check B^+}$.}
\label{alg:weymouthExtremesPlus}
\begin{algorithmic}[1]
\REQUIRE $ \munderbar{x}, \bar{x}, \munderbar{y}, \bar{y},  {\alpha} ,  {\beta} $ such that $ 0 < \munderbar{x} \le \bar{x}$, $0 < \munderbar{y} \le \bar{y}$, ${\alpha}\ge 0$,  ${\beta} \ge 0$.
\ENSURE $\mathcal{Q^+}$ 
\STATE $\mathcal{Q^+} = \emptyset$
%
\FOR{$(\tilde x, \tilde y) \in \textup{extr}(\{ (x,y)\in  [\munderbar{x}, \bar{x}]\times [\munderbar{y}, \bar{y}]: 0 \le x-y \le {\beta} \})$}  
\STATE Compute $z' = \sqrt{(\tilde x-\tilde y)/w}$
\IF{$\alpha \le  z'$} 
\STATE $\mathcal{Q^+}  =\mathcal{Q^+} \cup \{ (\tilde x, \tilde y, z'), (\tilde x, \tilde y, \alpha) \}$
\ENDIF
\ENDFOR

\FOR{$  \tilde z  \in\{0, \alpha\}  $}
%
\FOR{$  \tilde x  \in \{\munderbar{x}, \bar{x}\}$}
\STATE Compute $y_1 = \max \{ \munderbar{y}, \tilde x - \beta, \tilde x -w \tilde z^2\}$ and $y_2 = \min\{ \bar y, \tilde x\}$
\IF{$ y_1 \le y_2 $}
\STATE $\mathcal{Q^+}  =\mathcal{Q^+} \cup \{ (\tilde x, y_1, \tilde z), (\tilde x, y_2, \tilde z) \}$
\ENDIF    
\ENDFOR
%
\FOR{$  \tilde y  \in \{\munderbar{y}, \bar{y}\}$}
\STATE Compute $x_1 = \max \{ \munderbar{x}, \tilde y\}$ and $y_2 = \min\{ \bar x, \tilde y+\beta, \tilde y+w\tilde z^2\}$
\IF{$ x_1 \le x_2 $}
\STATE $\mathcal{Q^+}  =\mathcal{Q^+} \cup \{ (x_1,\tilde y,  \tilde z), ( x_2, \tilde y,\tilde z) \}$
\ENDIF
\ENDFOR
%
\FOR{$ d  \in \{ 0, \beta\}$}
\IF{$d \le w \tilde z^2$}
    \FOR{$  \tilde x  \in \{\munderbar{x}, \bar{x}\}$}
        \STATE  $y' = \tilde x - d $
        \IF{$\munderbar{y} \le y' \le  \bar{y}$}
            \STATE $\mathcal{Q^+}  =\mathcal{Q^+} \cup \{( \tilde x, y',  \tilde z)  \}$
        \ENDIF
    \ENDFOR
    
    \FOR{$  \tilde y  \in \{\munderbar{y}, \bar{y}\}$}
            \STATE  $x' = \tilde y + d $
            \IF{$\munderbar{x} \le x' \le  \bar{x}$}
            \STATE $\mathcal{Q^+}  =\mathcal{Q^+} \cup \{ (x', \tilde  y,  \tilde z)  \}$
            \ENDIF
    \ENDFOR
\ENDIF   
 \ENDFOR
 \ENDFOR
 \RETURN $\mathcal{Q^+} $ 
\end{algorithmic}
\end{algorithm}

Algorithm~\ref{alg:weymouthExtremesPlus} can also be used to find the convex hull of the  set
\(
\mathcal{X}^- = \{ (x,y,z) \in \mathcal{B}^- : y-x = w z^2   \}
\), 
where
$\mathcal{B}^- = \{ (x,y,z)\in  [\munderbar{x}, \bar{x}]\times [\munderbar{y}, \bar{y}]\times[ {\alpha}', 0] : {\beta}' \le x-y \le 0 \} $ such that ${\alpha}'\le 0$ and ${\beta}' \le 0$.  
For this purpose, we first   run Algorithm~\ref{alg:weymouthExtremesPlus} with the input $ \munderbar{y}, \bar{y},\munderbar{x}, \bar{x},   -{\alpha}' ,  -{\beta}' $ to obtain $\mathcal{Q^+}$. Then, for any element $(y,x,z)\in\mathcal{Q^+}$, we have $(x,y,-z)\in \mathcal{Q^-}$ such that $ \mathcal{\check B^-} = \textup{conv}(\mathcal{Q^-})$.
 
\end{appendices}

\end{document}